\documentclass[10pt]{article}

\usepackage{latexsym,color,amsmath,amsthm,amssymb,amscd,amsfonts,graphicx}
\usepackage{dsfont}
\usepackage{yfonts}

\usepackage{booktabs}

\usepackage{dsfont}

\newtheorem{theorem}{Theorem}
\newtheorem{corollary}{Corollary}

\newtheorem{lemma}{Lemma}

\newtheorem{proposition}{Proposition}

\newtheorem{example}{Example}

\newtheorem{remark}{Remarks}
        
              \def\iff{\Longleftrightarrow}

\def\R{{\mathbb R}}

\def\iff{\Longleftrightarrow}

\def\vol{{\rm vol}}

\def\1{\mathds{1}}

\def\be{\begin{equation}}
\def\ee{\end{equation}}

\def\grad{\nabla}

\DeclareMathOperator{\tr}{tr}

\newcommand{\norm}[1]{\lVert #1\rVert}

\newcommand{\dd}{\,d}

\title{Gravitational Illumination Bodies}
\author{Victor Glasgo and Elisabeth  M. Werner}

\date{}

\begin{document}
 
 \maketitle

\begin{abstract}
For a given convex body $K$ in $\mathbb R ^n$ and $\delta>0$, we introduce a new class of  bodies, the  {\em gravitational illumination bodies} $I^\delta(K)$. We show some of their properties and explore their 
relation to floating bodies, illumination bodies and Ulam floating bodies. Unexpected phenomena happen and 
the gravitational illumination bodies exhibit very different behavior 
than the 
floating bodies, illumination bodies and Ulam floating bodies. While the classical bodies are all convex and either contained in $K$ or containing $K$, this is no longer the case for the new bodies.
They need not be convex, but they are  star shaped and it may happen that $I^\delta(K)$ has non-empty intersection with both, $K$ and its complement. Moreover, the containment behavior may change with the dimension.
\par
Most surprising is  the behavior of a volume derivative, $$\lim_{\delta \to 0} \frac{ |\vol_n(K)-\vol_n(I^\delta(K)) |}{\delta^\frac{2}{n+1}}.$$ 
In all dimensions the corresponding volume derivative for the classical bodies results in the equi-affine invariant affine surface area of the body $K$.
In the case of the gravitational illumination body, this is also the case in all dimensions, except in dimension $3$.
 In dimension $3$ the limit  is equal to $0$ for all convex bodies.
A non-zero derivative is obtained in dimension $3$ if one normalizes by  $\delta$ instead of $\delta^\frac{1}{2}$. The non-zero limit is a new  affine invariant involving the {\em  affine mean curvature } and the {\em Pick invariant} of the convex body.
\end{abstract}

    

\vskip 3mm
\section { Introduction}
The illumination body of a convex body $K \subset \mathbb R^n$ was introduced in \cite{Werner94} as 
$$
K^\delta=\{x \in \mathbb R^n:   \vol_n\left( [x,K] \setminus K\right) \leq\delta\},
$$
where $[x,K]$ denotes the convex hull of $x$ and $K$ and $\vol_n(K)$ is the volume of $K$. 
We use the illumination body to 
 define  a new class of  bodies, the {\em gravitational illumination bodies} $I^\delta(K)$: For $\delta \geq 0$ given,  $x$ is in the boundary $ \partial I^\delta(K)$ of $I^\delta(K)$ if there exists $x^\delta$ in the boundary of $K^\delta$ such that $x$ is the center of gravity $g$ of $[x^\delta,K] \backslash K$, or
\begin{equation}\label{gib}
\partial (I^\delta(K))= \{ g([x^\delta,K] \backslash K) :  \vol_n\left([x^\delta,K] \backslash K\right) = \delta\}.
\end{equation}
The illumination body  can be viewed as a dual concept to the (convex) floating body introduced in \cite{BaranyLarman1988,  SchuettWerner1990},  the intersection of all halfspaces whose defining hyperplanes cut off a set of  fixed volume $\delta$ from $K$. The duality is made precise in \cite{MordhorstWerner2020}.  
The gravitational illumination body can be viewed as a dual concept to the {\em metronoid} \cite{HuangSlomka}. The metronoid  is also called   {\em Ulam floating body} \cite{HuangSlomkaWerner}, as it is related to Ulam's problem \cite{Ryabogin2022, Ulam},  or  {\em surface of centers}, e.g., \cite{DannHerscoviciMyrosh,FlorentinSchuettWZhang, Ryabogin2023,  Zawalski} as its boundary consists of the centers of gravity of the caps of  $K$ determined by hyperplanes that cut off $\delta$. 
\par
There has been much research in recent years related to those notions.
One example is the  work devoted to the  homothety problem of Sch\"utt and Werner \cite{AlfonsecaNazarovRyaboginStancuYaskin2025, SchuettWerner1994, Stancu2006, WernerYe2011, Zawalski}. Other examples are the extenions of 
floating bodies and illumination bodies to spherical and 
hyperbolic space \cite{AssoulineBesauWerner, BesauWerner2016,BesauWerner2018}, to Riemannian manifolds \cite{AssoulineSchuettWerner},  and to the class of $R$-ball convex bodies \cite{OlivaWernerYalikun, SchuettWernerYalikun2025}. 
We refer to e.g., \cite {ArtsteinFlorentin,  BezdekLangiNazodi, ChernovDrachTatarko, DrachTatarko, DrachTatarko2026, FodorKeveiVigh, FodorPapvari, FodorGrunfelder, LangiNazodiTalata} for more on  this class.
\vskip 2mm
Floating bodies, illumination bodies  and Ulam floating bodies are always convex. 
This is also the case for gravitational illumination bodies {\bf in the plane}. But in higher dimensions 
a first surprising behavior happens. There are examples of polytopes in $\mathbb{R}^3$ such that the gravitational illumination body is not convex for all
small enough $\delta$. However, for large enough $\delta$ the gravitational illumination body of a polytope is star shaped and even convex (depending on $\delta$),  also in dimension higher than $2$,  (see Theorem \ref{polytopes} for the precise formulation).
The gravitational illumination bodies of $C^2_+$ convex bodies 
are star shaped for small $\delta$.
\par
Floating bodies $K_\delta$,   illumination bodies $K^\delta$ and Ulam floating bodies $M_\delta(K)$ 
exhibit consistent behavior in all dimensions and for all $\delta$ with respect to containment with respect to $K$. Floating bodies and Ulam floating bodies are always contained in $K$ and illumination bodies always contain $K$. 
This is no longer the case for the gravitational illumination bodies. 
We show that for the Euclidean ball there  is a change in behavior from dimension $3$ to dimension $4$. 
Even ``worse":  Examples in dimension $2$ and $3$ illustrate that the gravitational illumination body can have non-empty intersection with both, $K$ and $K^C$, the complement of $K$. 
When we compare relations between floating bodies, Ulam floating bodies, illumination bodies and gravitational illumination bodies,  further interesting dimension-dependent features appear.
This is illustrated in the next table  in the case of the Euclidean ball, where we have for  sufficiently small positive $\delta$
depending on $n$, 
\begin{center}
\small
\begin{tabular}{@{}ll@{}}
\toprule
Dimension & Inclusion chain for the Euclidean ball $K=B_2^n$\\
\midrule
$n=2,3$
 & $K_\delta\subset M_\delta(K)\subset K
       \subset I^\delta(K)\subset K^\delta$\\[3pt]
$4\le n\le11$
 & $K_\delta\subset M_\delta(K)\subset I^\delta(K)
       \subset K\subset K^\delta$\\[3pt]
$12\le n\le23$
 & $K_\delta\subset I^\delta(K)\subset M_\delta(K)
       \subset K\subset K^\delta$\\[3pt]
$n\ge24$
 & $I^\delta(K)\subset K_\delta\subset M_\delta(K)
       \subset K\subset K^\delta$\\
\bottomrule
\end{tabular}
\end{center}
\par
We show that  a suitable derivative of volume  of a smooth convex body and its gravitational illumination body,
\begin{equation}\label{limit-intro}
\lim_{\delta \to 0} \frac{ |\vol_n(K)-\vol_n(I^\delta(K)) |}{\delta^\frac{2}{n+1}}
\end{equation}
 gives rise, in {\bf all dimensions except in dimension $3$},  to the {\em affine surface area}, an important equi-affine invariant in
convex geometry and  differential geometry with  applications to many other areas.  Similar results have been established  in all dimensions,  {\bf including dimension $3$},  for  floating bodies \cite{SchuettWerner1990}, illumination bodies \cite{Werner94}, Ulam floating bodies 
 \cite{HuangSlomkaWerner}, and more \cite{SchuettWerner2004}, in the decades long  effort to extend affine surface area from smooth bodies to all convex bodies.
\par
Originally, the affine surface area $as(K)$ of a smooth convex body $K$ was introduced by Blaschke  \cite{Blaschke:1923}
as an integral over the boundary $\partial K$ of  $K$ of a power of the Gauss curvature $\kappa$, namely
\begin{equation}\label{asa}
as(K) = \int_{\partial K} \kappa(x)^\frac{1}{n+1} d\mu_K(x),
\end{equation}
where $\mu_K$ is the usual surface area measure on $\partial K$.
Due to its important properties, which make it an effective and powerful tool, it
is omnipresent in geometry, e.g., 
\cite{GHSW:2020, Hoehner:2022, LR:2010, Lutwak:1996, LYZ:2000}).
There are numerous  applications and extensions of affine surface area 
such as, 
$L_p$-versions \cite{HaddadLangharstPuttermanRoysdonYe, Lutwak:1996, SchuettWerner2004, SchuettWerner2023, TW:2019, TW:2023, WernerYalikun, WernerYe:2008, Ye:2015a, Ye:2015, Ye:2016}, 
extremal problems \cite{BasitHoehnerLangiLedford, GHSW:2020, Hoehner:2022}, 
the approximation theory of convex bodies by polytopes 
\cite{BLW:2018, BGT, Boe1, 
Hoehneretal2024, ProchnoSchuettWerner2022, 
Reitzner2002, Reitzner, SchuettWerner2003}, Minkowski problems \cite{GuoXiZhao, HuangLutwakYangZhang}, affine curvature flows, information theory, and differential equations \cite{Andrews:1999, BesauWerner1, BesauWerner, TW1, TW2, TW4}.
\vskip 2mm
In dimension $3$ yet again an extremely surprising phenomenon occurs. The derivative of volume corresponding to (\ref{limit-intro}) is then {\bf equal to $0$ for all} convex bodies. 
The asymptotic expansion of the volume difference shows that a non-zero derivative is obtained if one uses the normalization $\delta$ instead of $\delta^\frac{1}{2}$ and the non-zero limit  
then gives a quantity that is not just  {\bf equi-affine invariant}, but {\bf affine invariant} involving  the  {\em affine mean curvature} and the {\em Pick invariant} of the convex body, see Theorem \ref{main3} for the precise formulation.
We investigate this new invariant in higher dimensions in a forthcoming work \cite{BesauSchuettW2027}.
\vskip 4mm
\noindent
{\bf Further notation.} 
The  closed Euclidean ball centered at $a$ with radius $r$ is $B^n_2(a,r)$. We write in short $B^n_2=B^n_2(0,1)$ and $S^{n-1}= \partial B^n_2$. The Euclidean norm on $\mathbb{R}^n$
is $\| \cdot\|$.
For $u \in S^{n-1}$, we denote by  $H\left(x, u\right)$ the hyperplane through $x$ orthogonal to the vector $u$. 
$H^-(x,u)$ and  $H^+(x,u)$ are the closed halfspaces determined by $H\left(x, u\right)$.
The line segment joining $x$ and $y$ is $[x,y]$. 
For a convex body $K$ in $\mathbb{R}^n$, $\vol_n(K)$ is its $n$-dimensional volume and 
$\mu_K$ is  the usual surface area measure on $\partial K$.
We write $\sigma= \mu_{B^n_2}$.
By $N(x)$ or $N_K(x)$ we denote the unit outer normal to $K$ in the boundary point $x$ and 
$\kappa(x)=\kappa_K(x)$ is the Gauss curvature in $x$. 
A convex body $K$ is $C^k_+$ if it has $k$-times  continuously differentiable boundary with everywhere strictly positive Gauss curvature. 
\par
\noindent
For more information and details on convex bodies we refer to e.g., the books  \cite{GardnerBook,  SchneiderBook}.
\par
\noindent
Finally, without loss of generality, we assume throughout the paper that $0 \in \text{int} (K)$, the interior of $K$,  and 
 $c$, $c_1$, $d$,  $d_1$, etc.,  are absolute constants that may change from line to line.

\vskip 4mm
\noindent
{\bf Acknowledgements.}
Elisabeth Werner was supported by NSF grant DMS-2506790 and by Simons Fellowship SFI-MPS-SFM-00020853. 
\vskip 3mm
\noindent

\section{The main theorems}

We now state  our main theorems.  We need the notion of a star shaped set.
\vskip 2mm
\noindent
Let $L\subseteq\mathbb{R}^n$ with $0 \in \text{int}  (L)$. The set $L$ is
\emph{star shaped with respect to $0$} if
\[
  [0,x]
  \subseteq L
  \qquad\text{for every }x\in L.
\]
\vskip 3mm
\noindent
We first address convexity  properties of the gravitational illumination bodies.
Again surprising phenomena happen, depending on the dimension.  
\vskip 2mm
\noindent
\begin{theorem}\label{polytopes}
\par
\noindent
(i) Let $P$ be a convex polygon in $\R^2$. 
Then $I^\delta(P)$ is convex.
\vskip1mm
\noindent
(ii) Let $K$ be a convex body in $\R^2$. Then $I^\delta(K)$ is convex.
\vskip 1mm
\noindent
(iii)  Let $P$ be a convex polytope in $\R^n$, $ n\geq 3$, with the center of gravity $g(P)=0$. Then there exists $\delta_1=\delta_1(P)$ such that for all $\delta \geq \delta_1$, $I^\delta(P)$ is starshaped.
\vskip 1mm
\noindent
(iv)  Let $P$ be a convex polytope in $\R^n$, $ n\geq 3$,  with the center of gravity $g(P)=0$. Then there exists $\delta_2=\delta_2(P)$ such that for all $\delta \geq  \delta_2$, $I^\delta(P)$ is convex.
\end{theorem}
\vskip 2mm
\noindent
The gravitational illumination body is  not  convex in general.  In the Appendix we give examples of polytopes in $\mathbb{R}^3$ that are not convex for all
small enough $\delta$. 
\vskip 3mm
\noindent
On the other hand in the smooth case we have the following behavior.
\vskip 2mm
\begin{theorem}\label{starshaped}
Let $K\subset\R^n$ be of class $C^2_+$ and suppose that the center of gravity of $K$, $g(K)=0$. 
Then there is $\delta_0=\delta_0(K)>0$ such that for all $\delta < \delta_0$, 
$I^\delta(K)$ is  starshaped with respect to the origin.
\end{theorem}
\vskip 3mm
\noindent
For the asymptotic behavior we recall the symmetric difference metric $d_s$ which for two measurable sets $A$ and $B$ is the volume of the symmetric difference,
$$
d_s(A,B) = \vol_n(A\setminus B) + \vol_n(B\setminus A).
$$
\begin{theorem}\label{main}
Let $K$ be a convex body in $\mathbb{R}^n$ that is $C^2_+$.	 Then for  $n \neq 3$,
\begin{align*}
		\lim_{\delta \to 0} \frac{ \vol_n(K) - \vol_n( I ^\delta(K)) }{\delta^\frac{2}{n+1}}
		= \alpha_n  \int_{\partial K}  \kappa(x) ^\frac{1}{n+1}  d\mu_K (x)= \alpha_n \, as(K), 
	\end{align*}
where $\alpha_n=\frac{1}{2} \  \frac{n-3}{n+3} \left( \frac{n(n+1)}{\vol_{n-1}(B^{n-1}_2)}\right)^\frac{2}{n+1}$.
\par
\noindent
For $n=3$, 
\begin{align*}
		\lim_{\delta \to 0} \frac{ d_s(K, I ^\delta(K)) }{\delta^\frac{1}{2}}
		= 0.
	\end{align*}
\end{theorem}
\vskip 2mm
\noindent
{\bf Remarks.}
\par
\noindent
(i) The case $n=3$ is extremely surprising. For one, we have to consider the symmetric difference as in dimension $3$ even in the smooth case 
we have in general neither $I^\delta(K) \subset K$ nor $K \subset I^\delta(K)$.  An example is given in Remark \ref{Beispiele}, 3.
And then,  the corresponding limit is equal to $0$ for {\bf all} convex bodies.
\par
\noindent
Note also that  for $n=2$ the right hand side of Theorem \ref{main} is strictly smaller than $0$. That the left had side is also  strictly smaller than $0$
follows from  Corollary \ref{cor2}.
\par
\noindent
Let $K$ be a $C^2_+$ convex body in $\mathbb{R}^n$. If we consider  $F(n) = \lim_{\delta \to 0} \frac{d_s(K, I ^\delta(K)) }{\delta^\frac{2}{n+1}}$,  as a function of the dimension $n$, then $F$ has a minimum at $n=3$ and is strictly positive for $n\neq3$.
\par
\begin{figure}[ht]
    \centering
    \includegraphics[width=2in]{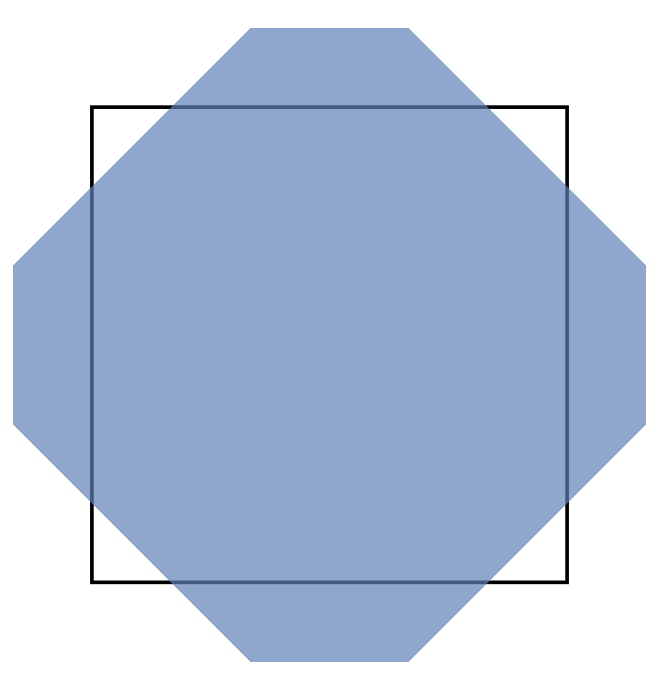}
    \caption{The gravitational illumination body of the square in $\R^2$ with $\delta$ sufficiently small neither contains the square nor is contained in the square.}
    \label{FigureGIBsquareexample}
\end{figure}

\noindent
(ii) We cannot expect that Theorem \ref{main} holds without any assumptions on $K$, even in dimension $2$. For instance, when $P$ is a polytope, for $\delta >0$ small enough, 
we have in general neither $I^\delta(P) \subset P$ nor $P \subset I^\delta(P)$. 
An example is the square in $\mathbb R^2$  centered at $0$
with side length $2$, see Figure \ref{FigureGIBsquareexample}.
\par
\noindent
The boundary of the gravitational illumination body of the square in the first quadrant is given by 
\[
\begin{cases}
y= 1 +\frac{\delta}{3},  & 0 \leq x \leq \frac{1}{3}   \\
y = -x  +\frac{4+\delta}{3} & \frac{1}{3} \leq x \leq 1 + \frac{\delta}{3}\\
x=1 +\frac{\delta}{3},  & 0 \leq y\leq \frac{1}{3}.
\end{cases}
\]
Furthermore, we observe  a surprising limiting behavior in the case of polytopes, namely that $\lim_{\delta\to 0}I^\delta(P)\neq P$. The limiting shape in the first quadrant is 
\[
\begin{cases}
y= 1,  & 0 \leq x \leq \frac{1}{3}   \\
y = -x  +\frac{4}{3} & \frac{1}{3} \leq x \leq 1 \\
x=1,  & 0 \leq y\leq \frac{1}{3}.
\end{cases}
\]
\vskip 3mm
\noindent
To determine the limit in dimension $3$, we need the {\em affine mean curvature 
 \(H_{\mathrm{aff}}\)}  and  the {\em Pick invariant $J$} which we introduce next.  For details we refer to \cite{NomizuSasaki}.
\par
\noindent
 Let \(x\in\partial K\). Let $T_x(\partial K)$ be the tangent space to $\partial K$ in $x$ and let \(\xi:\partial K\to\mathbb R^n\) be the inward affine  normal. Let $S_{\mathrm{aff}}$
 be the affine shape operator, i.e., 
\[
S_{\mathrm{aff}}(x)w= -D_w \xi
\qquad\text{for }w\in T_x(\partial K), 
\]
where \(D_w\xi\) differentiates \(\xi\) along the surface in the direction \(w\).
The eigenvalues of the affine shape operator are the affine principal curvatures. Then 
\[
H_{\mathrm{aff}}=\frac{1}{n-1} \operatorname{tr}S_{\mathrm{aff}}.
\]
The Blaschke connection \(\nabla^{\mathrm{aff}}\) and Blaschke metric \(h\) are defined by
\[
D_vw=\nabla^{\mathrm{aff}}_vw+h(v,w)\xi
\]
for tangent vectors  \(v,w\), where \(D\) is the ambient derivative. The affine cubic form at $x$ is the covariant tensor \(\mathcal C_x=\nabla^{\mathrm{aff}}h_x\), so \(\mathcal C_x(v,w,z)=(\nabla^{\mathrm{aff}}_vh)(w,z)\).
For $n=2$, the Pick invariant is $0$ for all $x \in \partial K$. For $n >2$, the Pick invariant at \(x\in \partial K\) is
\[
J(x)=\frac{\|\mathcal C_x\|_{h_x}^2}{4(n-1)(n-2)}, 
\]
where
\(\|\cdot\|_{h_x}\) is the norm induced by the Blaschke metric at \(x\). For a tangent vector \(v\in T_x(\partial K)\), it is \(\|v\|_{h_x}=\sqrt{h_x(v,v)}\).

\vskip 2mm
\noindent
Then, in dimension $3$,  the limiting behavior  is as follows.
\vskip 2mm
\begin{theorem}\label{main3}
	Let $K$ be a convex body in $\mathbb{R}^3$ that is $C^4_+$. Then 
\[
\lim_{\delta \to 0} \frac{d_s(K, I ^\delta(K)) }{\delta} = \frac{1}{12 \pi}  \int _{\partial K}  |J(x) -9\,  H_{\mathrm{aff}}(x)| \,  \kappa(x)^{1/4} d\mu_K (x)   \]
In particular, for the Euclidean ball $r B^n_2$ with radius $r$ we get
\[\lim_{\delta \to 0} \frac{\vol_3( I ^\delta(r B^3_2)) - \vol_3(r B^3_2)}{\delta} = 3. \]
\end{theorem}
\vskip 3mm
\noindent
{\bf Remarks.}
\vskip 2mm
\noindent
(i) The measure $\kappa^{1/(n+1)} d\mu_K$ is sometimes called {\em equiaffine surface area measure}  $dV_{aff}$ and ${\cal P}(K)=\int _{\partial K}  J V_{aff}$ is the {\em total Pick invariant}.
By the Pick-Berwald theorem \cite{NomizuSasaki},   ${\cal P}(K)=0$ if and only if $K$ is an ellipsoid.
The total mass of the {\em equiaffine mean curvature} $H_{aff}$, $\int _{\partial K}  H_{aff} V_{aff}$,  has appeared in  \cite{StancuWerner}.
\vskip 2mm
\noindent
(ii) We can define the  analog of the new  invariant on the right hand side of Theorem \ref{main3} for convex bodies  in $\mathbb{R}^n$. This is investigated in the forthcoming work \cite{BesauSchuettW2027}.

\vskip 4mm

\section{Properties of  the gravitational illumination bodies}

For $n \neq 3$, the volume derivative of Theorem \ref{main} results in the {\em equiaffine invariant} affine surface area,  i.e. 
$as(T K) =  as (K)$ for all affine maps $T$ on $\mathbb{R}^n$  such that $|\text{det}\  T|=1 $.
The next proposition shows that for $n=3$, the volume derivative of Theorem \ref{main3} gives a quantity that is affine  invariant,
i.e. gives the same quantity for {\bf all} affine maps independent of their determinant.
\vskip 2mm
\noindent

\begin{proposition} \label{prop3} 
\par
\noindent
(i) Let $K$ be a convex body in $\mathbb{R}^n$. Let $T:  \mathbb{R}^n  \rightarrow \mathbb{R}^n$ be an affine map such that $\text{det}\  T \neq 0$. Then 
$$ I^\delta( TK) = T \left( I^\frac{\delta}{|\text{det} \, T|}(K)\right).$$
\par
\noindent
(ii) Let  $K$ be a convex body in $\mathbb{R}^3$ that is $C^4_+$. Let $T:  \mathbb{R}^n  \rightarrow \mathbb{R}^n$ be an affine map such that $\text{det}\  T \neq 0$. Then
$$
\lim_{\delta \to 0} \frac{d_s(TK,  I ^\delta(T K))}{\delta}= \lim_{\delta \to 0} \frac{d_s(K, I ^\delta(K))}{\delta}.
$$
\end{proposition}
\vskip 2mm
\noindent
\begin{proof} 

(i) Let $y_0 \in \partial I^\delta( TK)$. This means that there exists $y^\delta$ such that $\int_{[y^\delta, TK] \setminus TK} dz  = \delta$ and
\begin{equation*}\label{y0}
y_0  =\frac{ \int_{[y^\delta, TK] \setminus TK} z  dz}{\int_{[y^\delta, TK] \setminus TK} dz} = \frac{1}{\delta}  \int_{[y^\delta, TK] \setminus TK} z  dz.
\end{equation*}
Then
\begin{eqnarray*}
T^{-1} \left( [y^\delta, TK] \right) &=& T^{-1} \left( \left\{  \alpha y^\delta + (1-\alpha) Tx: x \in K, 0 \leq \alpha \leq 1\right\}\right)\\ &=& 
\left\{  \alpha T^{-1} (y^\delta)  + (1-\alpha) x: x \in K, 0 \leq \alpha \leq 1\right\} = [x ^\delta, K],
\end{eqnarray*}
or, equivalently, $[y^\delta, TK]= T \left( [x ^\delta, K]\right)$,
where we have put $x^\delta = T^{-1} (y^\delta)$.  This implies that 
$$
[y^\delta, TK] \setminus TK = T \left( [x ^\delta, K] \setminus K \right)
$$
and 
\begin{equation}\label{x0}
\delta = \vol _n\left([y^\delta, TK] \setminus TK \right) =  | \det T | \    \vol _n  \left( [x ^\delta, K] \setminus K \right).
\end{equation}
Thus,
\begin{eqnarray*}
y_0 &= &\frac{ \int_{[y^\delta, TK] \setminus TK} z  dz}{\delta} = 
\frac{ \int_{T  \left(   [x ^\delta, K] \setminus K\right) } z  dz}{ | \det T | \    \vol _n  \left( [x ^\delta, K] \setminus K \right)} =
\frac{ | \det T |  \  \int_{  [x ^\delta, K] \setminus K\ } Tx   dx}{ | \det T | \    \vol _n  \left( [x ^\delta, K] \setminus K \right)}\\
& =& 
T \left(  \frac{\int_{  [x ^\delta, K] \setminus K\ } x  dx}{ \vol _n  \left( [x ^\delta, K] \setminus K \right)}\right) =  
T \left(  \frac{\int_{  [x ^\delta, K] \setminus K\ } x  dx}{\delta / | \det T | }\right)
\end{eqnarray*}
Together with (\ref{x0}), this means that
$$
I^\delta( TK) =T\left( I^\frac{\delta}{|\text{det} \ T|}(K)\right).
$$
\vskip 3mm
\noindent
(ii) We use (i) and get
\begin{eqnarray*}
\lim_{\delta \to 0} \frac{d_s(TK,  I ^\delta(T K))}{\delta}&=& \lim_{\delta \to 0} \frac{\vol_3( I ^{\delta}(T K)\setminus TK) + \vol_3(TK)\setminus I ^\delta(T K))}{\delta}\\
&=& \lim_{\delta \to 0} \frac{\vol_3(  T(I ^{\delta/| \det T |}(K)\setminus TK) + \vol_3(TK)\setminus T(I ^{\delta/| \det T |}(K))}{\delta}\\
&=& \lim_{\delta \to 0}|\det T | \frac{\vol_3(  (I ^{\delta/| \det T |}(K)\setminus K) + \vol_3(K)\setminus (I ^{\delta/| \det T |}(K))}{\delta}\\
&=& \lim_{\delta \to 0}\frac{\vol_3(  (I ^{\delta/| \det T |}(K)\setminus K) + \vol_3(K)\setminus (I ^{\delta/| \det T |}(K))}{\frac{\delta}{| \det T |}}\\
&= &  \lim_{s \to 0} \frac{d_s(K,  I ^s( K))}{s} .
\end{eqnarray*}
\end{proof}
\vskip 3mm
\noindent
The following proposition shows that the  gravitational illumination bodies behave nicely under convergence in Hausdorff distance.
We recall the Hausdorff distance between nonempty sets $A,B\subset \R^n$,
$$d_H(A,B)=\max\{\sup_{a\in A}d(a,B), \sup_{b\in B}d(A,b)\},$$
where $d(x,A)=\inf_{a\in A}d(x,a)$ is the distance from any point $x\in\R^n$ to the subset $A\subset\R^n$.
\vskip 3mm
\noindent
\begin{proposition}\label{seqpoly}
Let $K$ be a convex body in $\mathbb{R}^n$ and fix $\delta>0$. Let $P_i$ be a sequence of convex bodies  such that $P_i\to K$ in the Hausdorff metric. Then
\par
\noindent
(i) $P_i^\delta\to K^\delta$.
\par
\noindent
(ii) $I^\delta(P_i)\to I^\delta(K)$.
\end{proposition}
\vskip 2mm
\noindent
\begin{remark}
    In particular, this holds if $P_i$ is a sequence of polytopes and it  holds in every dimension $n\geq2$ and requires neither smoothness nor
strict convexity of $K$.
\end{remark}

\begin{proof}
For a convex body $A$ and $x \in A^C$, define
\[
 \Phi_A(x)
 =\vol_n\bigl([x,A]\setminus A\bigr)
 =\vol_n\bigl([x,A]\bigr)-\vol_n(A).
\]
The convex-hull operation is continuous in the Hausdorff metric: if
$A_i\to A$ and $x_i\to x$, then
\[
 [x_i,A_i]\longrightarrow
 [x,A].
\]
Volume is continuous under Hausdorff convergence of convex bodies.
Consequently,
\[
 \Phi_{A_i}(x_i)\longrightarrow\Phi_A(x).
\]
In particular, $\Phi_{P_i}\to\Phi_K$ uniformly on each fixed bounded
subset of $\R^n$.
\par
\noindent
Without loss of generality we assume that $0 \in\operatorname{int}(K)$.  For all sufficiently large $i$, the point
$0$ also belongs to $\operatorname{int}(P_i)$.  Along each ray
$0+tu$, $u\in S^{n-1}$, the function
\[
 t\longmapsto\Phi_A(0+tu)
\]
is continuous and strictly increasing after the ray leaves $A$.  Moreover,
the sets $P_i^\delta$ are uniformly bounded.  The uniform convergence of the
functions $\Phi_{P_i}$ therefore implies convergence of their
$\delta$-level sets, and hence
\[
 d_H(\partial P_i^\delta,\partial K^\delta)\longrightarrow0,
\]
which proves (i).
\par
\noindent
We next prove convergence of the centers of gravity.  Suppose that
\[
 x_i\in\partial P_i^\delta,
 \qquad x_i\longrightarrow x\in\partial K^\delta,
\]
and set
\[
 Q_i=[x_i,P_i],
 \qquad Q=[x,K].
\]
Then $Q_i\to Q$ in the Hausdorff metric.  First moments, just like volume,
are continuous under Hausdorff convergence of uniformly bounded convex
bodies.  Thus
\[
 \int_{Q_i}z\,dz\longrightarrow\int_Qz\,dz,
 \qquad
 \int_{P_i}z\,dz\longrightarrow\int_Kz\,dz.
\]
Since
\[
 \vol_n(Q_i\setminus P_i)=\vol_n(Q\setminus K)=\delta,
\]
we have, with $g_K(x)=g([x,K]\backslash K)$,
\begin{align*}
 g_{P_i}^\delta(x_i) &=\frac{1}{\delta} \left(\int_{Q_i}z\,dz-\int_{P_i}z\,dz\right)
 \longrightarrow \frac{1}{\delta} \left(\int_Qz\,dz-\int_Kz\,dz\right) =g_K^\delta(x).
\end{align*}
\par
\noindent
Finally, let $y_i\in I^\delta(P_i)$.  After passing to a subsequence, write
$y_i=g_{P_i}(x_i)$ with $x_i\to x\in\partial K^\delta$.  The preceding
argument gives $y_i\to g_K(x)\in I^\delta(K)$.  Conversely, for every
$x\in\partial K^\delta$, the convergence of the boundaries allows us to
choose $x_i\in\partial P_i^\delta$ such that $x_i\to x$, and then
$g_{P_i}(x_i)\to g_K(x)$.  These two assertions give
\[
 d_H\bigl(I^\delta(P_i),I^\delta(K)\bigr)\longrightarrow0.
\]
\end{proof}
\vskip 3mm
\noindent

\vskip 2mm
\subsection {Euclidean balls}
In this section we explicitly treat the case of  Euclidean balls $r B_2^n= B^n_2(0,r)$ with radius $r > 0$. This will show in particular  that  we cannot expect that  the gravitational illumination bodies are always contained in the given body not that they always contain the body.
\par
\noindent
For $x^\delta$ in the boundary of $(rB_2^n)^\delta$, we put  
$$A =\overline{ [x^\delta, rB_2^n] \backslash (rB_2^n)}.$$
We want to compute the center of gravity of $A$, which can be describes as follows.
All tangent lines to $rB_2^n$ through $x^\delta$ touch $rB_2^n$ in a hyperplane $H$ which cuts out an $(n-1)$-dimensional ball $(rB_2^n) \cap H$ of $rB_2^n$. Let $H^-$ be the halfspace 
in which $x^\delta$ sits. Let $B$ be the cap $B=rB_2^n \cap H^-$.
Then 
$A=\overline{[x^\delta, rB_2^n \cap H] \backslash B}$.
\vskip 2mm
\noindent
For the next proposition, we recall the Beta function and the  incomplete Beta function
$$
B(a,b) = \int _0^1 u^{a-1} (1-u)^{b-1} du, \hskip 10mm
B_t( a, b) = \int _0^t u^{a-1} (1-u)^{b-1} du.
$$ 
\vskip 3mm
\noindent
\begin{proposition} \label{prop1}
	For every $\delta > 0$, $I^\delta(rB_2^n)$ is a ball with radius
$$g_n (r, \delta) =\frac{r \  \frac{1}{n+1}\Delta^\frac{n+1}{2}\frac{(1-\frac{\Delta}{2r})^\frac{n+1}{2}}{(1-\frac{\Delta}{r})^2}}{\Delta^\frac{n+1}{2}\frac{(1-\frac{\Delta}{2r})^\frac{n+1}{2}}{(1-\frac{\Delta}{r})}-n \  r^\frac{n+1}{2}2^\frac{n-1}{2} B_{\frac{\Delta}{2r}}\left(\frac{n+1}{2}, \frac{n+1}{2}\right)},
$$
where $\Delta = \Delta(\delta)$ is the height of the cap $B$.
\end{proposition}
\vskip 2mm
\noindent
\begin{remark}
(i) In the case of low odd dimensions, we can compute $g_n(r, \delta)$ exactly. For instance, in dimensions $3$ and $5$ we get
\begin{align} \label{dim3}
g_3(r, \delta) &= r\frac{(1-\frac{\Delta}{2r})^2}{1-\frac{\Delta}{r}} = r + \frac{3 \delta}{4 r^2 \vol_2(B^2_2)}, \\
g_5 (r, \delta) &= r\frac{(1-\frac{\Delta}{2r})^3}{(1-\frac{\Delta}{r})(1-\frac{\Delta}{4r})} \nonumber.
\end{align}
In dimension $2$, we have
\begin{eqnarray} \label {dim2}
g_2(r, \delta) &=&   \frac { \frac r 3 (2\Delta)^{\frac 3 2} \frac{(1- \frac \Delta {2r})^{\frac 3 2}}{(1- \frac \Delta r)^2}}
{(2\Delta)^{\frac 3 2} \frac{(1- \frac \Delta {2r})^{\frac 3 2}}{1-\frac \Delta r} + r(2\Delta)^{\frac 1 2}(1-\frac \Delta r)(1-\frac \Delta {2r})^{\frac 1 2} - r^{\frac 3 2}\arccos(1-\frac \Delta r)} \nonumber \\
&=& r + \frac {1}{5} \left(\frac{9 \, \delta^2}{2\,  r}\right)^\frac{1}{3} + o( \delta^\frac{2}{3}).
\end{eqnarray}
\par
\noindent
(ii) It is convenient to rewrite $g_n(r, \delta)$. We put $t=\frac{\Delta}{2 r}$. Then $0 \leq t < \frac{1}{2}$ and
\begin{eqnarray} \label {beta}
g_n(r, t) &=& \frac{r}{(n+1) (1-2t)} \  \left( 1 - \frac{n}{2} \  (1-2t) \   \frac{B_t\left(\frac{n+1}{2}, \frac{n+1}{2}\right)}{t^\frac{n+1}{2} (1-t)^\frac{n+1}{2}} \right) ^{-1} \nonumber \\
&=& \frac{r}{(1-2t)} \  \left( n+1 - n \  (1-2t) \  \left( 1 + \sum _{k=0}^\infty \frac{B(\frac{n+3}{2}, k+1)}{B(n+1, k+1)} \  t^{k+1}\right) \right)^{-1} \nonumber\\
&=& \frac{r}{(1-2t)} \  \left( 1 +  2 n t - n \  (1-2t) \  \sum _{k=0}^\infty \frac{B(\frac{n+3}{2}, k+1)}{B(n+1, k+1)} \   t^{k+1} \right)^{-1}.
\end{eqnarray}
The second last identity follows by e.g., \cite{AbramowitzStegun}, p. 944.

\end{remark}
\par
\noindent
\begin{proof} 
	By symmetry,  all illumination bodies $(rB_2^n)^\delta$, $\delta > 0$, and consequently all bodies $I^\delta(rB_2^n)$ are again balls.
	Thus, it suffices to compute the distance of the  center of gravity of $A = \overline{[x^\delta,rB_2^n] \backslash (rB_2^n)}$ to the origin for one particular choice of $x^\delta$ in the boundary of $(rB_2^n)^\delta$.
	From now on we assume that $x^\delta$ lies on the positive $x_n$-axis.
	Since $A$ is symmetric about the $x_n$-axis, the center of gravity of $A$ lies on the same axis and $g_n(r, \delta) $ is equal to its $x_n$-coordinate.
	Letting $(g(A))_n$ to be the $x_n$-coordinate of the center of gravity $g(A)$,  then
	\begin{equation}
		\label{centerdecomp}
		g_n (r, \delta) = (g(A))_n = \frac 1 {\vol_n(A)} \{ (\vol_n(A)+\vol_n(B))(g(A \cup B))_n - \vol_n(B) (g(B))_n \}.
	\end{equation}
In terms of $r$ and $\Delta$, 
\begin{equation} 	\label{xdelta}
	\|x^\delta\|=\frac{r^2}{r-\Delta}.
\end{equation}
Hence, by a standard computation,
\begin{eqnarray*}
&&	(\vol_{n}(A)+\vol_{n}(B))(g(A\cup B))_n  \\
	&&=  \frac{\vol_{n-1}(B_2^{n-1})}{n(n+1)} \frac{(2r\Delta - \Delta^2)^{\frac{n-1} 2}}{(\|x^\delta\| - (r-\Delta))^{n-1}} \|x^\delta\| (\|x^\delta\| - (r-\Delta))^n \left[ 1 + n\frac{r-\Delta}{\|x^\delta\|} \right]   \\
	&&= \frac{\vol_{n-1}(B_2^{n-1})}{n(n+1)} (2r\Delta - \Delta^2)^{\frac{n-1} 2} (\|x^\delta\| - (r-\Delta)) (\|x^\delta\| + n(r-\Delta)).
\end{eqnarray*}
Plugging in expression \eqref{xdelta} for $\|x^\delta\|$ yields
\begin{align*}
	\left(\vol_n(A)+\vol_n(B)\right) (g(A \cup B))_n = \frac{\vol_{n-1}(B_2^{n-1})}{n(n+1)}(2r\Delta - \Delta^2)^{\frac{n+1} 2} \frac{r^2 + n(r-\Delta)^2}{(r-\Delta)^2} 
\end{align*} 
and similarly, 
\begin{eqnarray*}
&&\vol_n(B) (g(B))_n =\vol_{n-1}(B^{n-1}_2)\int_{r-\Delta}^r y(r^2-y^2)^\frac{n-1}{2}dy \\
&&= \frac{\vol_{n-1}(B_2^{n-1})}{n+1} r^{n+1} \left( 2\frac \Delta r - \left( \frac \Delta r \right)^2 \right)^{\frac{n+1} 2}
=\frac{\vol_{n-1}(B^{n-1}_2)}{n+1}(2r\Delta-\Delta^2)^\frac{n+1}{2}.
\end{eqnarray*}
Thus, the numerator in \eqref{centerdecomp} is equal to
\begin{equation*}
	\left(\vol_n(A)+\vol_n(B)\right)(g(A \cup B))_n - \vol_n(B) (g(B))_n = \frac{\vol_{n-1}(B_2^{n-1})}{n(n+1)} (2r\Delta - \Delta^2)^{\frac{n+1} 2} \frac{r^2}{(r-\Delta)^2}.
\end{equation*}
For the denominator we have $\vol_n(A)= \vol_n(A \cup B)  - \vol_n(B)$, where for the cone $A \cup B$
\begin{align*}
\vol_n(A\cup B)  &= \frac{1}{n}\vol_{n-1}(B^{n-1}_2)\  (2r\Delta-\Delta^2)^\frac{n-1}{2}\  (\|x^\delta\|-(r-\Delta)) \\
&= \frac{1}{n} \vol_{n-1}(B^{n-1}_2) \   \frac{(2r\Delta-\Delta^2)^\frac{n+1}{2}}{r-\Delta},
\end{align*}
and for the cap $B$
$$\vol_n(B)=\vol_{n-1}(B^{n-1}_2)\int_{r-\Delta}^r \left(r^2-y^2\right)^\frac{n-1}{2}dy.$$
The last integral can be rewritten in terms of an incomplete Beta function, by putting $u=\frac{r-y}{2r}$,
\begin{equation*}
	\int_{r-\Delta}^r (r^2-y^2)^{\frac{n-1} 2} dy = 2^n r^n \int_0^{\frac \Delta {2r}} u^{\frac{n-1} 2} (1-u)^{\frac{n-1} 2} du = (2r)^n B_{\frac{\Delta}{2r}} \left(\frac{n+1}{2}, \frac{n+1}{2} \right).
\end{equation*}
Collecting terms we get
\begin{align*}
	g_n (r, \Delta) &= \frac {  r \   \   \frac{1}{n+1}(2r\Delta)^{\frac{n+1} 2} \frac{(1- \frac \Delta {2r})^{\frac{n+1} 2}}{(1-\frac \Delta r)^2}} {(2 r \Delta)^{\frac{n+1} 2}  \frac{(1- \frac \Delta {2r})^{\frac{n+1} 2}}{1 - \frac \Delta r} - n \ 2^n\   r^{n +1}  \  B_{\frac{\Delta}{2r}} \left(\frac{n+1}{2}, \frac{n+1}{2} \right)}, 
\end{align*}
and cancelling out powers of 2 and $r$ yields the desired formula.
\end{proof}

\vskip 3mm
\noindent
The next lemma shows how $\delta$ and $\Delta$ are related.
\vskip 3mm
\begin{lemma} \label{deltaDelta} There are  constants $c_1(n)$ and $c_2(n)$ such that for $\delta$ sufficiently small, 
\begin{eqnarray*}
\hskip -30mm  && \frac{ \vol_{n-1}(B^{n-1}_2)}{n(n+1)} \left(2 \Delta \right)^\frac{n+1}{2} \  r^\frac{n-1}{2} \left( 1- c_1(n) \frac{\Delta }{r}\right)
 \leq  \  \delta  \\
&& \hskip 20mm \leq \frac{ \vol_{n-1}(B^{n-1}_2)}{n(n+1)}\left(2 \Delta \right)^\frac{n+1}{2} \  r^\frac{n-1}{2} \  \left( 1  + c_2(n)  \frac{\Delta}{r} \right).
\end{eqnarray*}
\end{lemma}
\vskip 2mm
\begin{proof}
The proof follows from  a formula in \cite{Werner94} using also (\ref{xdelta}).
\end{proof}
\vskip 3mm
\noindent
Recall that $t=\frac{\Delta}{2r}$ and observe that $t \rightarrow 0$ is equivalent to  $\Delta \rightarrow 0$ is equivalent to $\delta \rightarrow 0$. We will continue to write in short $\Delta$ for $\Delta(\delta)$.
\vskip 3mm
\begin{corollary} \label{Kor1}
The following statements hold true.
\par
\noindent
(i)   $\lim_{\Delta \rightarrow0} g_n(r, \delta) = r$. 
\vskip 2mm
\noindent
(ii) For all $0 \leq \Delta < r$,   $g_2(\delta) > r$ and $g_3(r, \delta) >  r$. 
\par
\noindent
For all $n \geq 4$ there is a $\Delta_0= \Delta(\delta_0)< r$ such that 
$$g_n(r, \delta_0) = r, $$ 
for all $\Delta \leq \Delta_0$, 
$$g_n(r, \delta) \leq r$$
and  for all $\Delta \geq \Delta_0$   we have that $$g_n(r, \delta) \geq r.$$
\end{corollary}
\vskip 2mm
\begin{proof}
(i) $\lim_{\Delta\rightarrow0} g_n(r, \delta) = \lim_{t \rightarrow0} g_n(r, t)$ and then the claim follows immediately by (\ref{beta}).
\vskip 2mm
\noindent
(ii) The cases $n=2$ and $n=3$ follow immediately from (\ref{dim2}) and (\ref{dim3}). 
For $n \geq 4$, we use (\ref{beta}) and instead of  the variable $\Delta$, we consider the variable $t$. Then
\begin{eqnarray} \label{equiv}
g_n(r, t) \leq r \,   \Longleftrightarrow \,  1 \leq  1-2t + 2nt(1-2t) - n (1-2t)^2 \  \sum _{k=0}^\infty \frac{B(\frac{n+3}{2}, k+1)}{B(n+1, k+1)} \   t^{k+1}, 
\end{eqnarray}
which is equivalent to 
\begin{eqnarray*} 
&&\hskip -5mm  0 \leq 2 \frac{n-3}{n+3} \  t  + 4 n t^2 \frac{n-7}{(n+3)(n+5)} - \\
&& \hskip 12mm   n \sum _{k=2}^\infty \left[\frac{B(\frac{n+3}{2}, k+1)}{B(n+1, k+1)}  - 4 \frac{B(\frac{n+3}{2}, k)}{B(n+1, k)} +
4 \frac{B(\frac{n+3}{2}, k-1)}{B(n+1, k-1)}  \right]\   t^{k+1},
\end{eqnarray*}
which holds for $t$ sufficiently small.
On the other hand, as $t \rightarrow \frac{1}{2}$, 
\begin{eqnarray*} 
 \frac{B_t\left(\frac{n+1}{2}, \frac{n+1}{2}\right)}{t^\frac{n+1}{2} (1-t)^\frac{n+1}{2}} \to   2 ^{n+1} \ \int _0^\frac{1}{2}  u^{\frac{n-1}{2}} (1-u)^{\frac{n-1}{2}}
\end{eqnarray*}
and 
$$ \frac{2}{n+1} \leq 2 ^{n+1} \ \int _0^\frac{1}{2}  u^{\frac{n-1}{2}} (1-u)^{\frac{n-1}{2}} \leq \frac{2^\frac{n+3}{2}}{n+1}$$
Thus, for $t \rightarrow \frac{1}{2}$, 
\begin{eqnarray*} 
g_n(r, t) &=& \frac{r}{(n+1) (1-2t)} \  \left( 1 - \frac{n}{2} \  (1-2t) \   \frac{B_t\left(\frac{n+1}{2}, \frac{n+1}{2}\right)}{t^\frac{n+1}{2} (1-t)^\frac{n+1}{2}} \right) ^{-1} 
\to \infty.
\end{eqnarray*} 
By this, (\ref{equiv}) and continuity, there is a $t_0=\frac{\Delta_0}{2r}$ such that $g_n(r, t_0) =r$.
\end{proof}
\vskip 2mm
\noindent
\begin{remark}
    The point (ii) of the corollary shows in particular  that in general we cannot expect that  the gravitational illumination bodies are always contained in the given body nor that they always contain the body.
\end{remark}

\subsection{Polytopes}

Next we derive a description of a boundary point in the gravitational illumination body. We use the following well known fact.
\par
\noindent
{\bf Fact.} The center of gravity of an $n$ dimensional cone of height $h$ lies on the line segment connecting the center of gravity of the base to the apex at a distance $\frac{h}{n+1}$ from the hyperplane containing the base.
\par
\noindent
Let $P$ be a convex polytope in $\mathbb{R}^n$. For $\delta>0$, let $P^\delta$ be the illumination body of $P$. For each $x^\delta\in \partial P^\delta$, $[x^\delta,P]\backslash P$ is a union of cones $C_i$, $1 \leq i \leq k$, 
$$
\overline{[x^\delta, P]\backslash P}=\cup_{i=1}^k C_i$$
where for $i=1,...,k$ , $C_i = [x^\delta, F_i]$ is the cone  with apex $x^\delta$ and base  $F_i$, with $F_i$ one of  the facets of $P$ 
that  is ``seen" from $x^\delta$. 
Using the above fact, it is easy to see that 
\begin{eqnarray*}
g(C_i)=\frac{x^\delta}{n+1}  + \left(1-\frac{1}{n+1}\right)g(F_i).
\end{eqnarray*}
Since $\overline{[x^\delta, P]\backslash P}=\cup_iC_i$,  and since $ \sum_{i=1}^k \vol_{n}(C_i) = \delta$,  we get the following formula for a point on the boundary of $I^\delta(P)$ 
\begin{align*} 
g\left([x^\delta, P]\backslash P\right)&=\frac{1}{\delta}\sum_{i=1}^k\operatorname{vol}_{n}(C_i)g(C_i)=\frac{1}{\delta} \sum_{i=1}^k\left[\vol_{n}(C_i) \left(\frac{x^\delta}{n+1}  + \left(1-\frac{1}{n+1}\right) g(F_i)\right)\right] \nonumber \\
&= \frac{x^\delta}{n+1} +  \left(1-\frac{1}{n+1}\right)  \sum_{i=1}^k  \frac{\vol_{n}(C_i)}{\delta}  \  g(F_i). 
\end{align*}
As $\sum_{i=1}^k  \frac{\vol_{n}(C_i)}{\delta} =1$,  the point $z(x^\delta) =  \sum_{i=1}^k  \frac{\vol_{n}(C_i)}{\delta}  \  g(F_i)$ is in $P$ and thus 
$g\left([x^\delta, P]\backslash P\right)=\frac{x^\delta}{n+1}  + \left(1-\frac{1}{n+1}\right)  z(x^\delta)$ is a convex combination of $x^\delta$ and  $z(x^\delta)$ which can also be written as
\begin{eqnarray}\label{eq:1}
g\left([x^\delta, P]\backslash P\right)&=& \frac{x^\delta}{n+1}  + \left(\frac{1}{n+1} \right)  \sum_{i=1} ^k  \frac{\vol_{n-1} (F_i) }{ \delta} \  \langle x^\delta -g(F_i), N_i \rangle  \  g(F_i).
\end{eqnarray}
where $N_i$ is the outer unit normal to $F_i$.
\vskip 2mm
\noindent
The following lemma establishes function notation to describe the boundary of $I^\delta (K)$ for a general convex body $K$. 

\begin{lemma} \label{LemmaTj} Let $P$ be a convex polytope in $\mathbb{R}^n$ and $\delta>0$ fixed. Let $G_j$ be a facet of $P^\delta$. Then for $x\in G_j$, the boundary of $I^\delta (P)$ is given (piecewise) by the maps 
$$T_j:G_j \rightarrow \partial I^\delta (P),$$
$$T_j(x):=\frac{x}{n+1}  + \left(\frac{1}{n+1} \right)  \sum_{i=1} ^k  \frac{\vol_{n-1} (F_i) }{ \delta} \  \langle x -g(F_i), N_i \rangle  \  g(F_i),$$
where the summation is taken over the facets of $P$ which are seen by $x$.
\vskip 2mm
\noindent
Similarly, define the map $T:\partial P^\delta \rightarrow \partial I^\delta (P)$, $T(x)=g(\overline{[x,P]\backslash P})$. $T$ is given explicitly by
\begin{equation}\label{T-see}
T(x):=\frac{x}{n+1}  + \left(\frac{1}{n+1} \right)  \sum_{i\in I_x}  \frac{\vol_{n-1} (F_i) }{ \delta} \  \max\{\langle x -g(F_i), N_i \rangle, 0 \} \  g(F_i),
\end{equation}
where $I_x=\{i: F_i \, \text{is a facet of}\,  P\}$.
\end{lemma}
\vskip 2mm
\begin{proof}
 The proof  follows from equation (\ref{eq:1}), noting that $x,y\in G_j$ if and only if $x$ and $y$ see the same set of facets $F_i\subset \partial P$, $i=1,...,k$.
\end{proof}
\vskip 3mm
\noindent
In the  next proposition, we give some properties of the maps defined in Lemma \ref{LemmaTj}.

\begin{proposition} \label{PropTj} Let $P$ be a convex polytope in $\mathbb{R}^n$ and $\delta>0$ fixed. Let $G_j$ be a facet of $P^\delta$, and the maps $T$ and $T_j$ defined as in Lemma \ref{LemmaTj}.

(i) $T_j$ preserves convex combinations,

(ii) If $x\in G_j\cap G_h$, then $T_j(x)=T_h(x)$.

(iii) $T$ is continuous.
\end{proposition}

\begin{proof} (i) Let $G_j$ be a facet of $P^\delta$. Then for all $x\in G_j$, $x$ sees the same set of facets $\{F_i\}_{i=1}^k$ of $P$. Take $x,y\in G_j$. Then by the convexity of $P^\delta$, $\lambda x+(1-\lambda)y\in G_j$. Then we have
\begin{eqnarray*}T_j(\lambda x+(1-\lambda)y)&=&\frac{\lambda x+(1-\lambda)y}{n+1}  \\
&+& \left(\frac{1}{n+1} \right)  \sum_{i=1} ^k  \frac{\vol_{n-1} (F_i) }{ \delta} \  \langle \lambda x+(1-\lambda)y -g(F_i), N_i \rangle  \  g(F_i)\\
&=&\frac{\lambda x+(1-\lambda)y}{n+1} + \left(\frac{1}{n+1} \right)  \sum_{i=1} ^k  \frac{\vol_{n-1} (F_i) }{ \delta} \\
&&\left( \langle \lambda x -g(F_i), N_i \rangle  +(1-\lambda) \langle  y -g(F_i), N_i \rangle\right)  \  g(F_i)\\
&=&\lambda T_j(x)+(1-\lambda)T_j(y).
\end{eqnarray*}
\vskip 2mm
\noindent
(ii) Let $x\in G_j\cap G_h$. Since $x\in G_j$, $x$ sees some set of facets $\{F_i\}_{i=1}^k$ of $P$ and forms a union of cones $\overline{[x,P]\backslash P}=\cup C_i$ which has volume $\delta$ and we have the map
$$T_j(x)=\frac{x}{n+1}  + \left(\frac{1}{n+1} \right)  \sum_{i=1} ^k  \frac{\vol_{n-1} (F_i) }{ \delta} \  \langle x -g(F_i), N_i \rangle  \  g(F_i).$$
Similarly, since $x\in G_h$, $x$ sees some other set of facets $\{F'_i\}_{i=1}^{k'}$ of $P$ and forms a union of cones $\overline{[x,P]\backslash P}=\cup C'_i$ with volume $\delta$ and we get the map
$$
T_h(x)=\frac{x}{n+1}  + \left(\frac{1}{n+1} \right)  \sum_{i=1} ^{k'}  \frac{\vol_{n-1} (F'_i) }{ \delta} \  \langle x -g(F'_i), N'_i \rangle  \  g(F'_i).$$
Consider any facet $F'\in \{F'_i\}\backslash\{F_i\}$ with normal vector $N$. Since, $\vol_n(\cup C_i)=\delta$, we have the volume of the cone $\vol_n([x,F'])=\frac{1}{n}\vol_{n-1}(F')\langle x-g(F'),N\rangle\leq 0$,  thus $\langle x-g(F'),N\rangle\leq 0$. But, since $x$ sees all the facets in $\{F'_i\}$, we have $\langle x-g(F'),N \rangle \geq 0$. So $\langle x-g(F'),N \rangle = 0$.
\par
\noindent
Likewise, for all facets $F\in \{F_i\}\backslash \{F'_i\}$with normal vector $N$,  $\langle x-g(F),N\rangle = 0$. Therefore $T_j(x)=T_h(x)$.
\vskip 2mm
\noindent
(iii) Follows from (i) and (ii) since each $T_j$ is shown to be an affine map restricted to the facet $G_j$.
\end{proof}
\vskip 3mm
\noindent
Now we consider polygons in $\R^2$.

\begin{proposition}\label{parallelpolygon}
Let $P$ be a convex polygon in $\R^2$. 
For every $\delta>0$, every
edge of $I^\delta(P)$ is parallel to the corresponding edge of
$P^\delta$.  

\end{proposition}

\vskip 2mm
\noindent
\begin{remark}
    Proposition \ref{parallelpolygon} is not true in general. We give a counter example in $\R^3$.
\end{remark}

\begin{example}\label{counterexampleparallel} 
Let $P\subset \R^3$ be the convex hull of the following vertices: 
$$P=[\{(0,0,0), (2,0,0),(0,2,0),(2,2,0),(0,0,-1), (6,0,-1),(0,2,-1),(2,2,-1)\}].$$
\end{example}
\par
\noindent
We do not have that all of the facets of $I^\delta(P)$ are parallel to the facets of $P^\delta$. The proof is given in the Appendix.

\begin{proof} [Proof of Proposition \ref{parallelpolygon}]

For $x\in\partial P^\delta$, put
\[
 D(x)=[x,P]\setminus P.
\]
Then $\vol_2(D(x))=\delta$ and  the centroid map is
\[
 g_P(x)
 =g(D(x))
 =\frac1\delta\int_{D(x)}y\,dy,
\]
and the boundary of the gravitational illumination body is traced by
\[
 \partial I^\delta(P)=g_P(\partial P^\delta).
\]
For vectors $a=(a_1,a_2)$ and $b=(b_1,b_2)$ in $\R^2$, we write
\(
 \det(a,b)=a_1b_2-a_2b_1.
\)
It's absolute value is the area of the
parallelogram spanned by $a$ and $b$.
\par
\noindent
Fix now an edge, or a subsegment of an edge, of $P^\delta$ on whose relative
interior the set of visible edges of $P$ is constant.  These visible
edges are consecutive.  Relabel their vertices as
\[
 v_0,v_1,\ldots,v_m,
\]
so that the visible edges are
\[
 [v_i,v_{i+1}],\qquad i=0,\ldots,m-1.
\]
Set
\[
 e_i=v_{i+1}-v_i,
 \qquad
 d=v_m-v_0=\sum_{i=0}^{m-1}e_i.
\]
For $x$ in the corresponding visibility region, let
\[
 \alpha_i(x)
 =\frac12\det(x-v_i,e_i).
\]
The orientation is chosen so that $\alpha_i(x)>0$.  Thus
$\alpha_i(x)$ is the area of the triangle
$[\{x,v_i,v_{i+1}\}]$.  The added region $D(x)$ is the union, with pairwise
disjoint interiors, of these triangles.  Hence
\begin{equation}\label{eq:area-sum}
 \vol_2(D(x))=\sum_{i=0}^{m-1}\alpha_i(x).
\end{equation}
We first determine the direction of the relevant edge of $P^\delta$.
For $t\in\R$, bilinearity of the determinant gives
\begin{eqnarray}\label{eq:alpha-variation}
 \alpha_i(x+td)
 &=&\frac12\det(x+td-v_i,e_i)\notag\\
& =&\alpha_i(x)+\frac{t}{2}\det(d,e_i).
\end{eqnarray}
Summing this  identity  over $i$ yields
\begin{align*}
 \sum_{i=0}^{m-1}\alpha_i(x+td)
 &=\sum_{i=0}^{m-1}\alpha_i(x)
   +\frac{t}{2}
     \det\left(d,\sum_{i=0}^{m-1}e_i\right)
 =\sum_{i=0}^{m-1}\alpha_i(x)
   +\frac{t}{2}\det(d,d)\\
 &=\sum_{i=0}^{m-1}\alpha_i(x).
\end{align*}
Thus the level line on which the added area equals $\delta$ is parallel
to $d$.  In particular, the edge of $P^\delta$ under consideration is
parallel to $d$.
\par
\noindent
The centroid of a triangle with vertices $x,v_i,v_{i+1}$ is
\[
 \frac{x+v_i+v_{i+1}}{3}.
\]
Using \eqref{eq:area-sum} and the fact that
$\sum_i\alpha_i(x)=\delta$ on $\partial P^\delta$, we obtain
\begin{align}
 g_P(x)
 &=\frac1{3\delta}
   \sum_{i=0}^{m-1}
   \alpha_i(x)(x+v_i+v_{i+1})\notag\\
 &=\frac{x}{3}
   +\frac1{3\delta}
    \sum_{i=0}^{m-1}\alpha_i(x)(v_i+v_{i+1}).
 \label{eq:centroid-map}
\end{align}
It follows from \eqref{eq:alpha-variation} and
\eqref{eq:centroid-map} that
\begin{equation}\label{eq:centroid-difference-preliminary}
 g_P(x+td)-g_P(x)
 =\frac{t}{3}d
  +\frac{t}{6\delta}
   \sum_{i=0}^{m-1}
   \det(d,e_i)(v_i+v_{i+1}).
\end{equation}
Let $A_{0,m}$ be the area enclosed by the polygonal chain
\[
 v_0,v_1,\ldots,v_m
\]
and the chord $[v_m,v_0]$.  We claim that
\begin{equation}\label{eq:vector-identity}
 \sum_{i=0}^{m-1}
 \det(d,e_i)(v_i+v_{i+1})
 =2A_{0,m}d.
\end{equation}
To verify this identity, translate and rotate the coordinate system so
that
\[
 v_0=(0,0),\qquad v_m=(L,0),\qquad d=(L,0),
\]
where $L>0$.  Write $v_i=(x_i,y_i)$.  Since $y_0=y_m=0$, the left-hand
side of \eqref{eq:vector-identity} becomes
\[
 L\sum_{i=0}^{m-1}
 (y_{i+1}-y_i)
 (x_i+x_{i+1},y_i+y_{i+1}).
\]
Its second component is
\[
 L\sum_{i=0}^{m-1}(y_{i+1}^2-y_i^2)=0.
\]
Its first component is
\begin{align*}
 L\sum_{i=0}^{m-1}
 (y_{i+1}-y_i)(x_i+x_{i+1})
 &=L\sum_{i=0}^{m-1}
 (x_i y_{i+1}-x_{i+1}y_i)\\
 &=2LA_{0,m}.
\end{align*}
The first equality also uses the telescoping identity
\[
 \sum_{i=0}^{m-1}
 (x_{i+1}y_{i+1}-x_i y_i)
 =x_my_m-x_0y_0=0,
\]
and the second equality is the shoelace formula.  Hence the vector is
$$
(2LA_{0,m},0)=2A_{0,m}d,
$$ 
proving
\eqref{eq:vector-identity}.
Substituting \eqref{eq:vector-identity} into
\eqref{eq:centroid-difference-preliminary}, we obtain the formula
\begin{equation}\label{eq:final-scaling}
 g_P(x+td)-g_P(x)
 =\frac{t}{3}
  \left(1+\frac{A_{0,m}}{\delta}\right)d.
\end{equation}
Because $P$ is convex, $A_{0,m}\geq0$.  The factor in
\eqref{eq:final-scaling} is therefore strictly positive.  It follows
that the centroid map sends the relevant edge segment of $P^\delta$ to
a nondegenerate segment parallel to it and with the same orientation.
\par
\noindent
This proves the assertion for every
$\delta>0$.

\end{proof}

\vskip 2mm
\noindent

\section {Proof of Theorem \ref{polytopes}}
\par
\noindent
(i) This follows immediately from Proposition \ref{parallelpolygon}. Indeed, for all $\delta$  the sides of $I^\delta(P)$ are parallel to the sides of $P^\delta$ so the interior angles of $I^\delta(P)$ and $P^\delta$ are equal. Convexity of $I^\delta(P)$ then follows from the convexity of $P^\delta$.
\vskip 2mm 
\noindent
(ii) follows from (i) by approximation by polytopes.
\vskip 3mm
\noindent
Before continuing the proof of Theorem \ref{polytopes} we need a further analysis of the boundary structure of $P^\delta$.
\vskip 3mm
\begin{lemma}
    Let $P$ be a convex polytope in $\R^n$ and fix $\delta>0$. Let $F^\delta$ be a facet of $P^\delta$.  Let  $x\in F^\delta$ be such that $x$  sees the set of $m$ facets of $P$, $\{F_i\}_{i=1}^m$. Let $N_i$ be the outer unit normal to the facet $F_i$ of $P$. Define the matrices, 
    $$F=
\begin{bmatrix}
\vol_{n-1}(F_1) &  & 0 \\
 & \ddots &  \\
0 &   & \vol_{n-1}(F_m)
\end{bmatrix}, G=
\begin{bmatrix}
 &  & \\
g(F_1) & \dots & g(F_m) \\
 &   & 
\end{bmatrix}, U=
\begin{bmatrix}
 &  & \\
N_1 & \dots & N_m \\
 &   & 
\end{bmatrix}\label{normalvectors}
$$
and $1=[1,...,1]^t\in\R^m$. Then
\par
    (i) $N_{F^\delta}=\sum_i^m \vol_{n-1}(F_i) N_i=UF1$ is a normal vector to $F^\delta\subset\partial P^\delta$, and
\par
    (ii) $N_{T(F^\delta)}=\sum_i^m \alpha_i N_i=U\alpha$ is a normal vector to $T(F^\delta)\subset\partial I^\delta(P)$, 
 \par
 \noindent   
    where $\alpha=(n+1)\delta (\delta I+FG^tU)^{-1}F1$ and where $T$ is as in (\ref{T-see}).
\end{lemma}
\vskip 2mm
\begin{proof}
    (i) For all $x\in F^\delta$,
    $$\delta=\sum \frac{1}{n}\vol_{n-1}(F_i)\langle x-g(F_i),N_i\rangle=\sum \frac{1}{n}\vol_{n-1}(F_i)\langle x,N_i\rangle-\sum \frac{1}{n}\vol_{n-1}(F_i)\langle g(F_i),N_i\rangle.$$
    Rearranging, we get
    $$\langle x, \sum \vol_{n-1}(F_i)N_i\rangle=n\delta + \sum \vol_{n-1}(F_i)\langle g(F_i),N_i\rangle.$$
    Thus we have found a vector $N_{F^\delta}$ such that $\langle x,N_{F^\delta}\rangle$ is constant for all $x\in F^\delta$.
\vskip 2mm
\noindent
    (ii) We write the map $T:F^\delta\to T(F^\delta)$ in matrix form
    \begin{align*}
      T(x)&=\frac{x}{n+1}+\frac{1}{(n+1)\delta}\sum \vol_{n-1}(F_i)\langle x-g(F_i),N_i\rangle g(F_i)\\
      &=\frac{x}{n+1}+\frac{1}{(n+1)\delta}\sum \vol_{n-1}(F_i)\langle x,N_i\rangle g(F_i)-\frac{1}{(n+1)\delta}\sum \vol_{n-1}(F_i)\langle g(F_i),N_i\rangle g(F_i)\\
      &=\frac{1}{(n+1)\delta}(\delta I+GFU^t)x-\frac{1}{(n+1)\delta}GF \text{diag}(U^tG),
    \end{align*}
    where $\text{diag}(U^tG)=(N_1^tg(F_1), ..., N_m^tg(F_m))^t$.
    Then,
    \begin{align*}
      \langle T(x), U\alpha\rangle &= \frac{1}{(n+1)\delta}\alpha^t U^t(\delta I+GFU^t)x-\frac{1}{(n+1)\delta}\alpha^t U^tGF \text{diag}(U^tG)\\
      &=\frac{1}{(n+1)\delta}\alpha^t (\delta I+U^tGF)U^tx-\frac{1}{(n+1)\delta}\alpha^t U^tGF \text{diag}(U^tG).
    \end{align*}
    Choosing $\alpha=(n+1)\delta (\delta I+FG^tU)^{-1}F1$, and noting from part (i) that $1^tF^tU^tx=n\delta+1^tF^t\text{diag}(U^tG)$, we get
    \begin{align*}
      \langle T(x), U\alpha\rangle &= 1^tF^tU^tx-\frac{1}{(n+1)\delta}\alpha^t U^tGF \text{diag}(U^tG)\\
      &=n\delta+1^tF^t\text{diag}(U^tG)-\frac{1}{(n+1)\delta}\alpha^t U^tGF \text{diag}(U^tG).
    \end{align*}
    Thus we have found a vector $N_{T(F^\delta)}$ such that $\langle T(x),N_{T(F^\delta)}\rangle$ is constant for all $T(x)\in T(F^\delta)$.
\end{proof}
\vskip 2mm
\noindent
\begin{remark}
    \par
\noindent
The matrix $(\delta I+FG^tU)$ is invertible as long as $ -\delta$ is not an eigenvalue of $FG^tU$.
\vskip 2mm
\noindent
\end{remark}
\vskip 3mm
\noindent
Now we continue the proof of Theorem \ref{polytopes}. We define the \emph{(convex) kernel} of a starshaped set, $S$,
$$\text{ker}(S)=\{x\in S:[x,y]\in S \text{ for all } y\in S\}.$$
That is, $\text{ker}(S)$ is the set of all points in $S$ from which all other points in $S$ are seen. It is clear from the definition $\text{ker}(S)$ is convex. Furthermore a set is starshaped if and only if $\text{ker}(S)\neq \emptyset$.

\begin{proof} [Proof of Theorem \ref{polytopes} (iii)]

    To show the set is starshaped, we show that $\text{ker}(I^\delta(P))\neq \emptyset$. Note that for nonconvex polytopes, a point is in the kernel if it is in the intersection of halfspaces defined as the affine span of the facets of the polytope.
   
    Since $g(P)=0$, we show that for all $T(x)\in \partial I^\delta(P)$ with outer normal $N_{T(x)}$ that we have $\langle T(x), N_{T(x)}\rangle>0$.
    In the proof of Lemma \ref{normalvectors}, we found
    \begin{align*}
      \langle T(x), N_{T(x)}\rangle&=n\delta+1^tF^t\text{diag}(U^tG)-\frac{1}{(n+1)\delta}\alpha^t U^tGF \text{diag}(U^tG)
    \end{align*}
    with $\alpha=(n+1)\delta (\delta I+FG^tU)^{-1}F1$. This simplifies to
    \begin{align*}
      \langle T(x), N_{T(x)}\rangle&=n\delta+1^tF^t\text{diag}(U^tG)-1^tF^t (\delta I+U^tGF)^{-1} U^tGF \text{diag}(U^tG)\\
      &=n\delta+1^tF^t\text{diag}(U^tG)-1^tF^t \frac{1}{\delta}\left(I+\frac{U^tGF}{\delta}\right)^{-1} U^tGF \text{diag}(U^tG)
    \end{align*}
    Then we write the Neumann expansion, for $\|\frac{U^tGF}{\delta}\|<1$, equivalently $\delta>\|U^tGF\|$ in the operator norm,    
    $$\left(I+\frac{U^tGF}{\delta}\right)^{-1}=I-\delta^{-1}U^tGF+\delta^{-2}(U^tGF)^2+o(\delta^{-2}), \qquad \delta\to\infty$$
    and we get
    \begin{align*}
      \langle T(x), N_{T(x)}\rangle&= n\delta+1^tF^t\text{diag}(U^tG)-1^tF^t\delta^{-1}U^tGF\text{diag}(U^tG)+o(\delta^{-1}), \qquad \delta\to\infty.
    \end{align*}
    Note $n\delta+1^tF^t\text{diag}(U^tG)>0$ since $F$ is diagonal with $F_{ii}>0$ and $\text{diag}(U^tG)$ has entries $\langle N_i,g(F_i)\rangle>0$ since $P$ is convex with $g(P)=0$. Then for $\delta$ sufficiently large $\langle T(x), N_{T(x)}\rangle>0$. Thus $0=g(P)\in \text{ker}(I^\delta(P))$ and therefore $I^\delta(P)$ is starshaped.
\end{proof}

\vskip 2mm
\noindent
\begin{proof} [Proof of Theorem \ref{polytopes} (iv)]
    Recall $T_j(x)=\frac{x}{n+1}  + \left(\frac{n}{n+1} \right)  \sum_{i=1} ^k  \frac{\vol_{n-1} (F_i) }{ n\delta}  \langle x -g(F_i), N_i \rangle  g(F_i)=\frac{x}{n+1}  + \left(\frac{n}{n+1} \right)  z(x)$, where $z(x)\in P$ as it is a convex combination of centers of gravity of facets of $P$. We have $\|z(x)\|\leq R$, for some constant $R$ since $P$ is bounded with $g(P)=0$. Then we get
    $$\lim_{\delta\to\infty}\frac{\|T_j(x)-\frac{x}{n+1}\|}{\|x\|}\leq \lim_{\delta\to\infty}\frac{R}{\|x\|}=0.$$
    Thus, for $\delta$ sufficiently large, the vertices of $I^\delta(P)$ are small perturbations of the vertices of $\frac{1}{n+1}P^\delta$, and therefore $I^\delta(P)$ is convex.
\end{proof}

\section {Proof of Theorem \ref{starshaped}}

\begin{proof} [Proof of Theorem \ref{starshaped}]

 Recall that $N_K(x)$ is  the outer unit normal to $\partial I^\delta(K)$ at
$x$.  For $x\in \partial K $ and $s>0$, set
\[
  y(x,s)=x+s \, N_K(x),
  \qquad
  V(x,s)= \vol_n([y(x,s),K]\setminus K)
\]
For fixed $x$, the function $s\mapsto \vol_n([y(x,s),K]\setminus K)$ is continuous and strictly
increasing near $s=0$, with $V(x,0)=0$.  Hence, for sufficiently small
$\delta>0$, there is a unique $s_\delta(x)>0$ satisfying
\begin{equation}\label{eq:volume-level}
  V\bigl(x,s_\delta(x)\bigr)=\delta.
\end{equation}
Every $y^\delta \in\partial K^\delta$ is of the form
$y^\delta=y(x,s_\delta(x))$ for a unique $x\in \partial K$: indeed, $x$ is the metric
projection of $y^\delta$ onto $K$.  Define
\[
  \gamma_\delta(x) = \frac{1}{\vol_n([y(x,s_\delta(x)),K]\setminus K)} \int _{[y(x,s_\delta (x)),K]\setminus K}  z dz, 
\hskip 5mm  x \in \partial K.
\]
Then
\begin{equation}\label{eq:gamma-image}
 \partial  I^\delta(K)=\gamma_\delta(\partial K).
\end{equation}
We first prove 
\begin{equation}\label{eq:C1-convergence}
  \gamma_\delta\longrightarrow Id_{\partial K}
  \qquad\text{in }C^1(\partial K, \R^n)
  \quad\text{as }\delta\to0.
\end{equation}
Fix $x\in \partial K$.  Choose orthonormal coordinates in which
\[
  x=0,
  \qquad
  N_K(x)=-e_n,
\]
and write the boundary of $K$ near $x$ as
\[
  y_x(z)=(z,f_x(z)),
  \qquad z\in\R^{n-1},
\]
where $K$ lies above the graph and
\[
  f_x(0)=0,
  \qquad
  \nabla f_x(0)=0.
\]
Put
\[
  q_x(z)=\langle\nabla f_x(z),z\rangle-f_x(z).
\]
For the apex $y(x,s)=-se_n$, the visible part of $\partial K$ is
parameterized by
\[
  D_{x,s}=\{z:q_x(z)\leq s\}.
\]
Indeed, since
\[
  N_K(y_x(z))
  =\frac{(\nabla f_x(z),-1)}
         {\sqrt{1+\lvert\nabla f_x(z)\rvert^2}},
\]
we have the identity
\begin{equation}\label{eq:visible-weight}
  \langle y(x,s)-y_x(z),N_K(y_x(z))\rangle
  \,d\mathcal H^{n-1}(y_x(z))
  =\bigl(s-q_x(z)\bigr)\,dz.
\end{equation}
Parameterizing $[y(x,s),K]\setminus K$, up to a set of measure zero, by
\[
  (t,z)\longmapsto y(x,s)+t \bigl(y_x(z)-y(x,s)\bigr),
  \qquad 0\leq t\leq1,\quad z\in D_{x,s},
\]
and using \eqref{eq:visible-weight}, we obtain
\begin{equation}\label{eq:local-volume}
  V(x,s)
  =\frac1n\int_{D_{x,s}}\bigl(s-q_x(z)\bigr)\,dz.
\end{equation}
Moreover,
\[
 \int_0^1
 \bigl(y(x,s)+t(y_x(z)-y(x,s))\bigr)t^{n-1}\,dt
 =\frac{y(x,s)+n y_x(z)}{n(n+1)}.
\]
Consequently, with
\[
  W(x,s)=\int_{D_{x,s}}\bigl(s-q_x(z)\bigr)\,dz,
\]
the centroid is given by the exact formula
\begin{equation}\label{eq:local-centroid}
  g(x,s)
  =\frac1{n+1}
   \left(
     y(x,s)
     +n\frac{\displaystyle
       \int_{D_{x,s}}(s-q_x(z))y_x(z)\,dz}
       {W(x,s)}
   \right).
\end{equation}
Because $K$ is of class $C^2_+$ and $\partial K$ is compact, all of the following
estimates are uniform in $x$.  We have
\begin{equation}\label{eq:Taylor}
  f_x(z)=\frac12\langle B_xz,z\rangle+o(\|z\|^2),
  \qquad
  q_x(z)=\frac12\langle B_xz,z\rangle+o(\|z\|^2),
\end{equation}
where $B_x=D^2f_x(0)$, the Hessian of $f_x$ at $0$,  is positive definite and its eigenvalues are bounded
above and bounded away from zero, independently of $x$.  It follows that
\begin{equation}\label{eq:domain-size}
  B_2^{n-1}(0,c\sqrt{s})
  \subset D_{x,s}
  \subset B_2^{n-1}(0,C\sqrt{s})
\end{equation}
and
\begin{equation}\label{eq:W-size}
  W(x,s)\asymp s^{\frac{n+1}{2}}.
\end{equation}
In particular,
\begin{equation}\label{eq:sdelta-size}
  s_\delta(x)\asymp\delta^{\frac{2}{n+1}}
\end{equation}
uniformly on $\partial K$.
Here, $h_1(\delta) \asymp h_2(\delta)$ means that there are constants $c, C$ such that $c h_1(\delta) \leq h_2(\delta) \leq C h_1(\delta)$ for all $\delta$ small enough.
In particular, $c, C$ and $\delta_0$, independent of $x$ and $\delta$ such that  for all $\delta < \delta_0$, for all $x \in \partial K$,  $c \delta^\frac{2}{n+1}  \leq s_\delta(x) \leq C \delta^\frac{2}{n+1} $.
Set
\[
  \overline z(x,s)
  =\frac1{W(x,s)}
    \int_{D_{x,s}}(s-q_x(z))z\,dz
\]
and
\[
  \overline f(x,s)
  =\frac1{W(x,s)}
    \int_{D_{x,s}}(s-q_x(z))f_x(z)\,dz.
\]
After the substitution $z=\sqrt{s}\,w$, the weighted domains converge,
uniformly in $x$, to
\[
  \left\{w:\frac12\langle B_xw,w\rangle\leq1\right\}
\]
with the even weight
\[
  1-\frac12\langle B_xw,w\rangle.
\]
The tangential first moment of the limiting weighted ellipsoid vanishes.
Together with \eqref{eq:Taylor}--\eqref{eq:W-size}, this gives
\begin{equation}\label{eq:moment-size}
  \overline z(x,s)=o(\sqrt{s}),
  \qquad
  \overline f(x,s)=O(s).
\end{equation}
Formula \eqref{eq:local-centroid} therefore implies
\begin{equation}\label{eq:C0-centroid}
  g(x,s)-x=o(\sqrt{s})
\end{equation}
uniformly in $x$.
\par
\noindent
It remains to control the derivatives.  On each of finitely many coordinate
neighborhoods of $\partial K$, the  graphs $f_x$ and the corresponding
orthonormal frames can be chosen to depend $C^1$-continuously on $x$.  Since
$f_x(0)=0$ and $\nabla f_x(0)=0$ for every $x$, differentiation with respect
to the boundary parameter gives
\begin{equation}\label{eq:p-derivatives}
  D_xf_x(z)=o(\|z\|),
  \qquad
  D_xq_x(z)=o(\|z\|),
\end{equation}
uniformly in $x$.
For small $s$, all the integrals may be taken over a fixed coordinate ball by
writing the integrands with the factor $(s-q_x)_+$.  Thus differentiation
produces no boundary term, since $s-q_x=0$ on $\partial D_{x,s}$.  From
\eqref{eq:domain-size} and \eqref{eq:p-derivatives},
\begin{equation}\label{eq:volume-derivatives}
  D_xW(x,s)=o(s^{n/2}),
  \qquad
  \partial_sW(x,s)\asymp s^{\frac{n-1}{2}}.
\end{equation}
Implicit differentiation of \eqref{eq:volume-level}, using
$V=W/n$, gives
\begin{equation}\label{eq:sdelta-derivative}
  D_xs_\delta(p)=o(\sqrt{s_\delta(x)}).
\end{equation}
Differentiating the two weighted first moments in the same way gives
\begin{equation}\label{eq:moment-p-derivatives}
  D_x\overline z(x,s)=o(1),
  \qquad
  D_x\overline f(x,s)=o(\sqrt{s}),
\end{equation}
and
\begin{equation}\label{eq:moment-s-derivatives}
  \partial_s\overline z(x,s)=o(s^{-1/2}),
  \qquad
  \partial_s\overline f(x,s)=O(1).
\end{equation}
For completeness, the estimates follow respectively from
\[
  \int_{D_{x,s}}D_xq_x(z)\,dz=o(s^{n/2}),
  \qquad
  \int_{D_{x,s}}D_xq_x(z)z\,dz
      =o(s^{\frac{n+1}{2}}),
\]
\[
  \int_{D_{x,s}}
  \bigl[-D_xq_x(z)f_x(z)+(s-q_x(z))D_x f_x(z)\bigr],dz
      =o(s^{\frac{n+2}{2}}),
\]
and, after scaling, from
\[
  \int_{D_{x,s}}z\,dz=o(s^{n/2}),
  \qquad
  \int_{D_{x,s}}f_x(z)\,dz=O(s^{\frac{n+1}{2}}).
\]
In a moving orthonormal frame, \eqref{eq:local-centroid} and
\eqref{eq:moment-size}--\eqref{eq:moment-s-derivatives} imply, at fixed $s$,
\[
  D_x\bigl(g(x,s)-x\bigr)=o(1),
  \qquad
  \partial_s\bigl(g(x,s)-x\bigr)=o(s^{-1/2}).
\]
Combining these estimates with \eqref{eq:sdelta-derivative} gives
\[
  D_x\bigl(g(x,s_\delta(x))-x\bigr)=o(1)
\]
uniformly in $x$.  Together with \eqref{eq:C0-centroid} and
\eqref{eq:sdelta-size}, this proves \eqref{eq:C1-convergence}.
\par
\noindent
We now finish the argument by radial projection.  Define
\[
  R_0:\partial K\longrightarrow \mathbb{S}^{n-1},
  \qquad
  R_0(x)=\frac{x}{\|x\|}.
\]
Because $K$ is convex and $0$ is in the interior of $K$, every ray from the origin meets
$\partial K$ exactly once.  Thus $R_0$ is a diffeomorphism.  Its differential is
nonsingular because
\[
  \langle x,N_K(x)\rangle=h_K(N_K(x))>0,
  \qquad x\in \partial K.
\]
For sufficiently small $\delta$, \eqref{eq:C1-convergence} allows us to
define
\[
  R_\delta:\partial K \longrightarrow \mathbb{S}^{n-1},
  \qquad
  R_\delta(x)=\frac{\gamma_\delta(x)}{\|\gamma_\delta(x)\|},
\]
and gives
\[
  R_\delta\longrightarrow R_0
  \qquad\text{in }C^1.
\]
Hence $R_\delta$ is a local diffeomorphism for all sufficiently small
$\delta$.
Let
\[
  r_0=\min_{x\in \partial K}\|x\|>0.
\]
For small $\delta$, we have
$\|\gamma_\delta(x)-x\|<r_0/2$ uniformly in $x$.  Therefore
\[
  H_t(x)
  =\frac{(1-t)x+t\gamma_\delta(x)}
         {\|(1-t)x+t\gamma_\delta(x)\|},
  \qquad 0\leq t\leq1,
\]
is a homotopy from $R_0$ to $R_\delta$.  It follows that
\[
  \deg R_\delta=\deg R_0=1, 
\]
where $\deg R$ denotes the topological degree of the map $R$, see e.g. \cite{Hatcher}.
A local diffeomorphism from the compact connected manifold $\partial K$ onto
$\mathbb{S}^{n-1}$ is a finite covering map.  The absolute value of its degree is
the number of sheets.  Hence $R_\delta$ has one sheet and is a global
diffeomorphism.
Consequently, every ray from the origin meets $\partial I^\delta(K)$ in exactly
one point.  There is therefore a positive $C^1$ function
$\rho_\delta:\mathbb{S}^{n-1}\to(0,\infty)$ such that
\[
  \partial I^\delta(K)
  =\{\rho_\delta(u)u:u\in \mathbb{S}^{n-1}\}.
\]
In particular, $\partial I^\delta(K)$ is an embedded hypersurface, and the
compact region it bounds is
\[
  I^\delta(K)
  =\{tu:u\in \mathbb{S}^{n-1},\ 0\leq t\leq\rho_\delta(u)\}.
\]
Thus, if $y\in I^\delta(K)$ and $0\leq\lambda\leq1$, then
$\lambda y\in I^\delta(K)$.  Hence $I^\delta(K)$ is star shaped with respect
to the origin.
\end{proof}
\vskip 3mm
\noindent
\subsection{Floating Bodies, Ulam Floating Bodies, Illumination Bodies and Gravitational Illumination Bodies}

In this subsection we explore the relations between floating bodies, Ulam floating bodies, illumination bodies and gravitational illumination bodies.
We start by introducing floating bodies and Ulam floating bodies.
\par
\noindent
Let $K$ be a convex body in $\mathbb{R}^n$ and let $\delta>0$. The (convex) floating body $K_\delta$ introduced independently in \cite{BaranyLarman1988} and \cite{SchuettWerner1990}
is defined as
\begin{equation}\label{floating body}
K_\delta= \bigcap _{\vol_n(H^- \cap K) \leq \delta} H^+.
\end{equation}
Let  $\mathds{1}_{K}$ be the characteristic function of $K$, and let 
$\mu$ be  the measure whose density with respect to Lebesgue measure
is $\delta^{-1}\mathds{1}_{K}$.Then the  Ulam floating body $M_\delta(K)$ is defined as \cite{HuangSlomka, HuangSlomkaWerner},
\begin{equation}\label{Ulam floating body}
M_\delta(K) 
=\bigcup_{\substack{0\le f\le1,\\
\int_{\R^{n}}f d{\mu}=1
}
}\left\{ \int_{\R^{n}}yf\left(y\right)d{\mu}\left(y\right)\right\},
\end{equation}
where the union is taken over all functions $0\le f\le1$ for which
$\int_{\R^{n}}f d{\mu=1}$ and $\int_{\R^{n}}yf\left(y\right) d{\mu}\left(y\right)$
exists.
\vskip 2mm
\noindent
The next corollary follows from the machinery developed in the proof of Theorem \ref{starshaped} and the definitions of $K_\delta$, $K^\delta$, $I^\delta(K)$ and $M_\delta(K)$.
It is well known \cite{HuangSlomkaWerner} that,  whenever $K_\delta$ is nonempty,
\begin{equation}\label{enthalt}
 K_\delta\subset M_\delta(K)\subset K\subset K^\delta.
\end{equation}
Moreover, the following inclusions hold.
\vskip 2mm 
\begin{corollary}\label{cor2}
Let $K$ be a convex body in $\mathbb{R}^n$ that is $C^2_+$.
\par
\noindent
(i) \(
 I^\delta(K)\subset K^\delta 
\) and $I^0(K)= K$. 
\par
\noindent
(ii) For $n=2$ and $\delta>0$ sufficiently small we have that  
$$ 
K_\delta \subset M_\delta(K) \subset K \subset I^\delta(K) \subset  K^\delta.
$$
For all $n \geq 4$ and for $\delta>0$ sufficiently small we have that 
$$
I^\delta(K) \subset K \subset K^\delta.
$$
(iii) For $n=3$ and $\delta>0$ sufficiently small we have that  
\[
K_\delta\subset I^\delta(K)\subset K^\delta.      
\]
\end{corollary}
\vskip 2mm
\begin{proof} 
\par
(i) Let $x\in\partial K^\delta$.  By  convexity of $K^\delta$, 
\[
 [x,K] \setminus K \subset  [x,K] \subset K^\delta.
\]
The barycentre of $[x,K] \setminus K$ thus belongs to $K^\delta$ which shows the first assertion of (i).
The second assertion of (i)  follows immediately from (\ref{eq:C0-centroid}). 
\vskip 2mm
\noindent
(ii) Because of (\ref{enthalt}) and (i), we only need to show that $K \subset I^\delta(K)$ for $n=2$ and $I^\delta(K)\subset K$ for $n\geq 4$. As in the proof of  Theorem (\ref{starshaped}), fix $x\in \partial K$ and choose orthonormal coordinates in which
\[
  x=0,
  \qquad
  N_K(x)=-e_n,
\]
with $f_x(z)$ and $q_x(z)$ defined as above. Let $B$ denote the osculating ball of $K$ at $x$ and define $f^B_x(z)$ and $q^B_x(z)$ analogously for the ball, that is, write the boundary of $B$ near $x$ as follows
\[
  y^B_x(z)=(z,f^B_x(z)),
  \qquad z\in\R^{n-1},
\]
where $B$ lies above the graph and
\[
  f^B_x(0)=0,
  \qquad
  \nabla f^B_x(0)=0.
\]
We define
\[
  q^B_x(z)=\langle\nabla f^B_x(z),z\rangle-f^B_x(z).
\]
Since $K$ and its osculating ball $B$ agree to second order at $x$, we have by (\ref{eq:Taylor}),
\[
f_x(z)-f^B_x(z)=o(\|z\|^2),
\qquad
q_x(z)-q^B_x(z)=o(\|z\|^2).
\]
On $D_{x,s}$ we have $\|z\|=O(\sqrt{s})$, hence
\[
f_x(z)-f^B_x(z)=o(s),
\qquad
q_x(z)-q^B_x(z)=o(s).
\]
Consequently, with the visible part of $\partial B$ parametrized by
\[
  D^B_{x,s}=\{z:q^B_x(z)\leq s\}
\]
we get
\[
  W^B(x,s)=\int_{D^B_{x,s}}\bigl(s-q^B_x(z)\bigr)\,dz\asymp s^{\frac{n+1}{2}}.
\]
Define the weighted integral for the ball $B$ analogously to $\overline{f}(x,s)$,
\[
  \overline f^B(x,s)
  =\frac1{W^B(x,s)}
    \int_{D^B_{x,s}}(s-q^B_x(z))f^B_x(z)\,dz.
\]
Then, after the scaling $z=\sqrt{s}w$, both weighted integrals converge uniformly to the same weighted ellipsoid determined by the common second fundamental form at $x$. Since $f_x(z)-f^B_x(z)=o(s)$ and $q_x(z)-q^B_x(z)=o(s)$, the difference of the corresponding normalized moments is $o(s)$. That is,
\[
\overline{f}(x,s)-\overline{f}^B(x,s)=o(s)
\]
uniformly in $x$.
From (\ref{eq:local-centroid}),
\[
g(x,s)
=
\frac1{n+1}
(n\overline{z}(x,s),-s+n\overline{f}(x,s)),
\]
it follows, with $g^B(x,s)$ the corresponding centroid for the ball $B$ and with $N_K(x)=N_B(x)$, that
\[
\langle g(x,s)-g^B(x,s),N_K(x)\rangle
=
-\frac n{n+1}
(\overline f(x,s)-\overline f^B(x,s))
=o(s).
\]
Using now (\ref{eq:sdelta-size}),
\[
\langle g(x,s_\delta(x))-g^B(x,s_\delta(x)),N_K(x)\rangle=o(\delta^\frac{2}{n+1}).
\]
From this and (\ref{beta}), for the ball we have
\[
\langle g^B(x,s_\delta(x))-x,N_K(x)\rangle
=
-\frac{1}{2}\frac{n-3}{n+3}\left(\frac{n(n+1)}{\vol_{n-1}(B_2^{n-1})}\right)^\frac{2}{n+1}r(x)^{-\frac{n-1}{n+1}}\delta^\frac{2}{n+1}+o(\delta^\frac{2}{n+1}),
\]
so we conclude for $K$ that
\begin{equation}\label{centroidexpansion}
    \langle g(x,s_\delta(x))-x,N_K(x)\rangle
=
-\frac{1}{2}\frac{n-3}{n+3}\left(\frac{n(n+1)}{\vol_{n-1}(B_2^{n-1})}\right)^\frac{2}{n+1}r(x)^{-\frac{n-1}{n+1}}\delta^\frac{2}{n+1}+o(\delta^\frac{2}{n+1}),
\end{equation}
where $r(x)$ is the radius of the osculating ball $B$ at $x$. Since the estimate $o(\delta^\frac{2}{n+1})$ holds uniformly in $x$, for $\delta$ sufficiently small, we get $K\subset I^\delta(K)$ for $n=2$ and $ I^\delta(K) \subset K$ for $n \geq 4$.
\vskip 2mm
\noindent
(iii) We only need to show that $K_\delta\subset I^\delta(K)$. We use the
normal parametrization $x\mapsto \gamma_\delta(x)$ from the proof of Theorem  \ref{starshaped}.
In dimension $3$, the  barycentre has zero inward normal
coordinate. For example, its first normal moment is proportional to
\[
 \int_0^1 r\left[r^4-(2r-1)^2\right]\,dr=0.
\]
Thus, uniformly in $x\in\partial K$,
\begin{equation}\label{eq:three-normal1}
 \langle \gamma_\delta(x)-x, N_K(x) \rangle = o(s_\delta(x))=o(\sqrt\delta),
\end{equation}
where \(s_\delta(x)\) is determined by
\[
\operatorname{vol}_n\!\left(
[x+s_\delta(x) u,K] \setminus K
\right)=\delta
\]
 as in the proof of Theorem  \ref{starshaped} and $u=N_K(x)$.
The cap and added-volume formulas become
\begin{equation*}
 \delta=\frac{\pi}{\sqrt{\kappa(x)}}c_\delta(u)^2(1+o(1)), \quad \quad
 \delta=\frac{\pi}{3\sqrt{\kappa(x)}}s_\delta(x)^2(1+o(1)),
\end{equation*}
where $c_\delta=h_K(u)-a_\delta(u)$ is the height of the cap
$$K\cap\{y:\langle y,u\rangle\ge a_\delta(u)\}$$
defining the floating body $K_\delta$.
So

\begin{equation}\label{eq:three-cap-height}
 c_\delta(u)=\frac{\kappa(x)^{1/4}}{\sqrt\pi}\sqrt\delta
             +o(\sqrt\delta).
\end{equation}
Since $\kappa$ has a positive minimum, \eqref{eq:three-normal1} and
\eqref{eq:three-cap-height} imply
\[
 \langle \gamma_\delta(x), u \rangle >h_K(u)-c_\delta(u)=a_\delta(u)
\]
for all $x$ and sufficiently small $\delta$. The same strict inequality
holds throughout the segment from $x$ to $\gamma_\delta(x)$.
Every such segment therefore avoids $K_\delta$, because the latter
lies in the halfspace $\langle y, u \rangle \le a_\delta(u)$. The straight-line
homotopy from $\partial K$ to the barycentre surface avoids every
point of $K_\delta$. Its degree remains $1$ about those points, so the
surface encloses $K_\delta$. This proves the  inclusion.

\end{proof}

\vskip 2mm
\noindent
\begin{remark} \label{Beispiele}
\par
\noindent
We cannot expect in general complete inclusion chains as in (ii) in dimension $n \geq 3$.  
\par
\noindent
1. For instance, the Ulam floating body need not be contained in the gravitational illumination body. Indeed,  
\[
I^\delta(B_2^{12})\subset M_\delta(B_2^{12}).
\]
This follows from Proposition \ref{prop1} and Lemma \ref{deltaDelta}  and the computations for the Ulam floating body of Euclidean balls in \cite{HuangSlomkaWerner}, which give
\begin{equation}\label{eq:three-radii}
 r_M (\delta) =1-B_n \Delta +O(\Delta^2),\qquad
        r_I(\delta)=1-A_n \Delta +O(\Delta ^2), 
\end{equation} 
where  $\Delta=\Delta(\delta)$ is as in Proposition \ref{prop1}, $r_M$ and $r_I$ are the radii of $M_\delta(B_2^{n})$ and $I^\delta(B_2^{n})$  respectively,  and 
\begin{equation} \label{AnBn}
A_n=\frac{n-3}{n+3},  \quad \quad
 B_n=\frac{n+1}{(n+3) n^{2/(n+1)}}. 
\end{equation}
\vskip 2mm
\noindent
2. Also, we cannot expect that $K_\delta\subset I^\delta(K)$ in general. We have 
\[
 I^\delta(B_2^{24})\subset(B_2^{24})_\delta.
\]
This follows from (\ref{eq:three-radii}) and (\ref{AnBn}) and the fact that the radius $r_F$ of the floating body is $r_F(\delta) = 1- \frac{\Delta}{n^{2/(n+1)}} +O(\Delta ^2)$ \cite{JenkinsonWerner}.
\vskip 2mm
\noindent
3. It follows from Proposition \ref{prop1} that in dimension $3$
\begin{equation*}\label{eq:three-ball-body}
I^\delta(B_2^3)
       =\left(1+\frac{3\delta}{4\pi}\right)B_2^3
       \supset B_2^3.
\end{equation*}
On the other hand the convex body 
\begin{equation*}\label{eq:explicit-K}
K=\left\{(x_1,x_2,x_3):
s-s^2+s^3+x_3^2\le1,\quad s=x_1^2+x_2^2
\right\}
\end{equation*}
provides an example of a $C^2_+$ convex body in $\mathbb{R}^3$ for which we have that $K$ is {\bf not} contained in $I^\delta(K)$ and $I^\delta(K)$ is {\bf not} contained in $K$.
We sketch the argument.
\par
\noindent
It is easy to see  that $K$ is $C^2_+$ and strictly convex.
We describe the boundary of $K$ in the first octant  by 
\begin{equation*}
 f(s)=(1- s + s^2- s^3)^{\frac{1}{2}}, 
\end{equation*}
where $s=(x_1^2+x_2^2)$. 
Let $p=(0,0,1)$. Let $t>0$ and let $p_t=(0,0,1+t)$. We compute $\delta=\delta(t)$ such that $p_t \in \partial K^\delta$. 
Then  the barycenter of $[p_t,K]\setminus K$ is
$$
(0,0,
1-\frac9{4\pi}\delta+o(\delta))
$$
and it belongs to the interior of $K$. On the other hand, it also belongs to $\partial I^\delta(K)$, thus it cannot be that $K\subset I^\delta(K)$. 
\par
\noindent
Let $q=(1,0,0)$.  Let $q^\delta$ be such that $q^\delta \in K^\delta$ with $\delta=\delta(t)$ given as above.
Then  the barycenter of $[q^\delta,K]\setminus K$ is
$$
(1+ \frac{15}{8\sqrt2\,\pi}\delta+o(\delta) , 0,0)
$$
and it belongs to the interior of the complement of $K$. On the other hand, it also belongs to $\partial I^\delta(K)$, thus it cannot be that $I^\delta(K) \subset K$. 
\end{remark}

\vskip 3mm
\noindent
We conclude this subsection with the inclusion chain for the Euclidean ball.
For $K=B_2^n$, with sufficiently small positive $\delta$
depending on $n$, the following inclusions hold. The proofs follow from computations of the radii in \cite{HuangSlomkaWerner, JenkinsonWerner, SchuettWerner1990,Werner94}.
\begin{center}
\small
\begin{tabular}{@{}ll@{}}
\toprule
Dimension & Inclusion chain for the Euclidean ball $K=B_2^n$\\
\midrule
$n=2,3$
 & $K_\delta\subset M_\delta(K)\subset K
       \subset I^\delta(K)\subset K^\delta$\\[3pt]
$4\le n\le11$
 & $K_\delta\subset M_\delta(K)\subset I^\delta(K)
       \subset K\subset K^\delta$\\[3pt]
$12\le n\le23$
 & $K_\delta\subset I^\delta(K)\subset M_\delta(K)
       \subset K\subset K^\delta$\\[3pt]
$n\ge24$
 & $I^\delta(K)\subset K_\delta\subset M_\delta(K)
       \subset K\subset K^\delta$\\
\bottomrule
\end{tabular}
\end{center}

\vskip 4mm

\section {Proof of Theorems \ref{main} and  \ref{main3} }

In what follows we shall write $H(x,u)$ for the hyperplane in $\mathbb{R}^n$ through $x\in \mathbb{R}^n$ with normal $u\in\partial B^n_2$ and $H(x,u)^-$ for the halfspace bounded by $H(x,u)$ into which $u$ points. In addition, we recall that $N_K(x)$ denotes the outer unit normal vector at $x\in\partial K$ of a convex body $K$. 
\par
\noindent
We will  also need several lemmas. 
 The first one can be found for convex bodies in  e.g., \cite{SchuettWerner1990} or in \cite{Werner94}. That it also holds for star-shaped bodies can be found in e.g.
\cite{WernerYe:2008}.
\vskip 2mm
\begin{lemma} \label{difference}
 Let $K$ and $L$ be starshaped bodies in $\mathbb{R}^{n}$
such that $0\in \text{int}(L)$ and $0\in \text{int}(K)$. When $L \subset K$, 
\[
\vol_{n}(K) - \vol_{n}(L) =\frac{1}{n}\int_{\partial K}\langle x,N_K(x)\rangle\left(1-\left(\frac{\|x_L\|}{\|x\|
}\right)^{n}\right) d{\mu_{K}\left(x\right)}.
\]
More generally, 
\[d_s(K, L)=\frac{1}{n}\int_{\partial K}\langle x,N_K(x)\rangle\left|1-\left(\frac{\|x_L\|}{\|x\|
        }\right)^{n}\right| d{\mu_{K}\left(x\right)}
\]
where $x_{L}$ is the unique point in the intersection of $\partial L$ with the ray $\{\lambda x:\lambda>0, x\in \partial K\}$. 
\end{lemma}
\vskip 3mm
\noindent
When  $K$ is  $C^2_+$,  the
Blaschke rolling theorem (see, e.g.,  \cite{SchneiderBook}) states  that  there exists a ball with radius $r_0$ that rolls freely in $K$:
For all $x \in \partial K$,
\begin{equation}\label{BlaschkeRolling}
B^n_2(x-r _0N_K(x), r_0) \subset K.
\end{equation}
We also recall that  $g^\delta(x)$ is the unique point in the intersection $\partial  I^\delta(K)  \cap \left\{ \lambda x : \lambda \geq 0 \right\}$.
\vskip 2mm
\noindent
\begin{lemma} \label{bounded}
Let $K$ be a convex body in $\mathbb{R}^n$  that is $C^2_+$ and such that $0 \in  \text{int}(K)$. 
Then there exists $\delta_0= \delta_0(K)$ such that for all $\delta \leq \delta_0$, for all $x \in \partial K$ the functions 
\[
 \frac{\langle x,N_K(x)\rangle}{ n \  \delta^{2/(n+1)}} \ \left|1-\left(\frac{\|g^\delta(x)\|}{\|x\|
}\right)^{n}\right| 
\]
are uniformly (in $\delta$)  bounded by a  function that is integrable on $\partial K$.
\par
\noindent
When $n=2$, they are bounded by $\gamma_2 \,  r_0 ^{-\frac{n-1}{n+1}}$, where $\gamma_2 $ is a constant. 
\par
\noindent
When $n \geq 4$, they are bounded by $\gamma(n) \,  r_0 ^{-\frac{2n}{n+1}}$ where $\gamma(n)$ is a constant that depends on $n$ only.
\par
\noindent
When $n=3$, and $K$ is $C^4_+$, then we have for $\delta \leq \delta_0$ that 
\begin{equation*}
 \frac{\langle x,N_K(x)\rangle}{3\,  \delta^\frac{1}{2}} \ \left|1-\left(\frac{\|g^\delta(x)\|}{\|x\|
}\right)^{3}\right| \leq  \gamma_3 \delta^\frac{1}{2},  
\end{equation*}
where $\gamma_3$ is an absolute constant constant.
\end{lemma}
\vskip 2mm
\begin{proof}
\par
\noindent
We first treat the case $n=2$.
\par
\noindent
Let $\{x^\delta \} = \{ \lambda \ x: \lambda \geq 0\} \cap \partial K^\delta $. By Corollary \ref{cor2} (ii), there is $\delta _0$ such that for all $\delta \leq \delta_0$,  $\|x^\delta\| \geq \|g^\delta(x)\|$ and thus 
\begin{eqnarray*}
&&\frac{1}{n} \  \frac{\langle x,N_K(x)\rangle}{\delta^\frac{2}{n+1}}  \left|1-\left(\frac{\|g^\delta(x)\|}{\|x\|}\right)^{n}\right|  =  
\frac{1}{n} \  \frac{\langle x,N_K(x)\rangle}{\delta^\frac{2}{n+1}}  \left[\left(\frac{\|g^\delta(x)\|}{\|x\|}\right)^{n} -1 \right] \\
&& 
\leq  \frac{1}{n}  \  \frac{\langle x,N_K(x)\rangle}{\delta^\frac{2}{n+1}}  \left [\left(\frac{\|x^\delta\|}{\|x\|}\right)^n -1 \right]
\end{eqnarray*}
and it follows from \cite{Werner94}, Lemma 2, that the quantity is bounded above by $\gamma_2 r_0 ^{-\frac{n-1}{n+1}}$.
\vskip 2mm
\noindent
Now we  treat the case $n \geq 4$. 
\par
\noindent
By Corollary \ref{cor2} (ii), there is $\delta_0$ such that  for all $\delta \leq \delta_0$, $I^\delta(K) \subset K$.  Thus 
$$
\left|1-\left(\frac{\|g^\delta(x)\|}{\|x\|}\right)^{n}\right| = 1-\left(\frac{\|g^\delta(x)\|}{\|x\|}\right)^{n}.
$$
By the colinearity of $x$ and $g^\delta(x)$, we have that $\left(\frac{\|g^\delta(x)\|}{\|x\|}\right)^{n} = \left(1- \frac{\|x-g^\delta(x)\|}{\|x\|} \right)^n$ and thus 
\begin{align}\label{1UNGL}
\frac{1}{n}  \langle x,N_K(x)\rangle \left(1-\left(\frac{\|g^\delta(x)\|}{\|x\|}\right)^{n}\right) \leq  \left\langle \frac{x}{\|x\|},N_K(x)\right\rangle  \|x-g^\delta(x)\| \leq \|x-g^\delta(x)\|.
\end{align}
Let $\delta_0$ be such that for all $x \in \partial K$, $\|x - g^{\delta_0}(x)\| \leq \frac{r_0}{2}$, where $r_0$ is the radius of the Blaschke rolling ball. Then we can view $g^{\delta}(x)$  for all
$\delta \leq \delta_0$  as an element of $\partial I^{\tilde{\delta}}(r_0 B^n_2)$ with $\tilde{\delta} \leq \delta$. 
By (\ref{beta}), we get with $\Delta=\Delta(\tilde{\delta})$and  a constant $c$, 
\begin{eqnarray*}
g_n \left(r_0, \frac{\Delta}{2r_0}\right)&=&\frac{r}{(1-\frac{\Delta}{r_0})} \  \left( 1 +   n \frac{\Delta}{r_0} - n \  \left(1-\frac{\Delta}{r_0}\right) \  \sum _{k=0}^\infty \frac{B(\frac{n+3}{2}, k+1)}{B(n+1, k+1)} \,  \left(\frac{ \Delta}{2r_0}\right)^{k+1} \right)^{-1}\\
&=& r_0\left( 1- \frac{n-3}{n+3} \frac{\Delta}{r_0} + \frac{3(n^2+15)}{(n+3)^2(n+5)} \left(\frac{\Delta}{r_0}\right)^2 + c  \left(\frac{\Delta}{r_0}\right)^3 \right), 
\end{eqnarray*}
which implies that 
\begin{eqnarray} \label{r+g}
r_0-g_n \left(r_0, \frac{\Delta}{2r_0}\right) &=& \frac{n-3}{n+3} \frac{\Delta}{r_0} - \frac{3(n^2+15)}{(n+3)^2(n+5)} \left(\frac{\Delta}{r_0}\right)^2 - 
c  \left(\frac{\Delta}{r_0}\right)^3 \leq \frac{n-3}{n+3} \frac{\Delta}{r_0} \nonumber\\
&\leq& \frac{n-3}{n+3} \left(\frac{n(n+1)}{\vol_{n-1}(B^{n-1}_2)}\right)^\frac{2}{n+1} \, \frac{\tilde{\delta}^\frac{2}{n+1}}{r_0^\frac{2n}{n+1}}.
\end{eqnarray}
In the last inequality we have also used Lemma \ref{deltaDelta}. We put $\langle \frac{x}{\|x\|},N_K(x)\rangle = \cos \theta$. As without loss of generality 
we can assume that $0 \in \text{int} (K)$, there exists $\beta >0$ such that 
\begin{equation}\label{beta2}
B^n_2(0,\beta) \subset K \subset B^n_2\left(0,\frac{1}{\beta}\right).
\end{equation}
This means that $\cos \theta \geq \beta^2$. By the law of cosine,
$$
g_n \left(r_0, \frac{\Delta}{2r_0}\right)^2= r_0^2 +\|x-g^\delta(x)\|^2 - 2 \cos \theta \, r_0 \|x-g^\delta(x)\|.
$$
Solving for $\|x-g^\delta(x)\|$, we get 
$$
\|x-g^\delta(x)\| \leq \frac{2}{\beta^2} \, \left(r_0 - g_n \left(r_0, \frac{\Delta}{2r_0}\right)\right).
$$
With this and  (\ref{r+g}), 
\begin{eqnarray*}
\frac{1}{n}  \frac{\langle x,N_K(x)\rangle }{\delta^\frac{2}{n+1}} \, \left(1-\left(\frac{\|g^\delta(x)\|}{\|x\|}\right)^{n}\right) &\leq& \frac{2}{\beta^2} 
\, \frac{(n-3)}{(n+3)} \left(\frac{n(n+1)}{|B^{n-1}_2|}\right)^\frac{2}{n+1} \,         r_0^{-\frac{2n}{n+1}}  \frac{\tilde{\delta}^\frac{2}{n+1}} {{\delta}^\frac{2}{n+1}} \\
&\leq&\frac{2}{\beta^2} \, \frac{(n-3)}{(n+3)} \left(\frac{n(n+1)}{\vol_{n-1}(B^{n-1}_2)}\right)^\frac{2}{n+1} \,         r_0^{-\frac{2n}{n+1}} .
\end{eqnarray*}
\vskip 2mm
\noindent
Finally, we treat the case $n=3$.
\par
\noindent
We first consider the case that  $g^\delta(x)$ is in the interior of the complement  of $K$,  i.e., 
$g^\delta(x) \in \text{int} (K^C)$. 
Then we get  with an absolute  constant $c$
\begin{eqnarray*}
\frac{1}{3} \, \langle x,N_K(x)\rangle \left|1-\left(\frac{\|g^\delta(x)\|}{\|x\|}\right)^{3}\right|  &=& 
\frac{1}{3} \,  \langle x,N_K(x)\rangle  \left[\left(\frac{\|g^\delta(x)\|}{\|x\|}\right)^{3} -1 \right] \\
&\leq&  c \, \langle \frac{x}{\|x\|}, N_K(x) \rangle \|g^\delta (x) - x\|.
\end{eqnarray*}
As  $g^\delta(x) \in \text{int} ((B^3_2(x-r_0 N_K(x), r_0))^C)$,  by (\ref{dim3}), $g^\delta(x)$ can be viewed, with a shift, as  an element of $\partial I^{\tilde{\delta}}(B^3_2(r_0))$
with $\tilde{\delta} \leq \delta$ and then 
\begin{eqnarray} \label{bounded-dim3}
\frac{1}{3} \, \langle x,N_K(x)\rangle \left|1-\left(\frac{\|g^\delta(x)\|}{\|x\|}\right)^{3}\right| 
&\leq&  c \, \langle \frac{x}{\|x\|}, N_K(x) \rangle \|g^\delta (x) - x\| = \langle g_3(r_0, \tilde{\delta}) , e_3 \rangle - r_0 \nonumber \\
&=&c\,  \frac{\tilde{\delta}}{4 \pi r_0^2} \leq c\,  \frac{{\delta}}{4 \pi r_0^2} .
\end{eqnarray}
Now we treat  the case that  $g^\delta(x) \in \text{int} (K)$. We have 
\begin{eqnarray*}
\frac{1}{3} \, \langle x,N_K(x)\rangle \left|1-\left(\frac{\|g^\delta(x)\|}{\|x\|}\right)^{3}\right|  &=& 
\frac{1}{3} \,  \langle x,N_K(x)\rangle  \left[1-\left(\frac{\|g^\delta(x)\|}{\|x\|}\right)^{3}  \right] \\
&\leq&   \langle \frac{x}{\|x\|}, N_K(x) \rangle \|g^\delta (x) - x\|.
\end{eqnarray*}
Let \(y^\delta \in  \partial K^\delta\) be the  point such that
\(
g^\delta(x)\) is the barycenter of the convex hull $y^\delta$ and $K$ without $K$, i.e. the barycenter of 
$[K, y^\delta]\setminus K$. Let $y$ be the closest point in $K$ to $y^\delta$. As $K$ is convex, 
$y$ is unique and we have that 
\[
y^\delta=y+t N_K(y),
\]
where $ t=\operatorname{dist}(y^\delta,K)$.
We estimate $\delta$ by the volume of the cone $[y^\delta, B^3_2(y- r_0 N_K(y), r_0)]$  intersected with the upper halfspace given by the tangent hyperplane to $\partial B^3_2(y- r_0 N_K(y), r_0)$ in $y$, 
and get
\begin{equation}\label{delta-Lemma5}
\delta\ge  \frac{\pi r_0}{6 M} \quad t^2 = c_1 \, t^2, 
\end{equation}
where $M$ is an absolute constant.
Moreover we claim that there is $t_0$ such that for all $t \leq t_0$, 
\[\big|\langle g^\delta(x)-y,N_K(y)\rangle\big|
\le c_2\,  t^2.
\]
To show this, we use coordinates such that \(y=0\),  \(N_K(y)=-e_3\) and  \(y^\delta =-t e_3\). We write locally, i.e., in a neighborhood $U(0)$,  the boundary of $\partial K$ as \((v,u(v))\), such that 
for all $v \in U(0)$, 
\[
u(v)=Q(v)+P(v)+O(\|v\|^4),
\]
with $Q(v)=\tfrac12v^T D^2 u  v$, where $D^2 u$ is the Hessian of $u$  and  \(P\) is homogeneous of degree $3$. Denoting $D u$ the gradient of $u$, we put
\[
W_t(v)=t+u(v)-\langle v,D u(v)\rangle +O(\|v\|^4)
=t - Q(v)-2P(v)+O(\|v\|^4).\]
The last equality holds as $\langle v,D u(v)\rangle = 2 Q(v) + 3 P(v)  +O(\|v\|^4)$, as $Q$ is homogeneous of degree $2$ and $P$ is homogeneous of degree $3$.
We also set
\[ D_t=\{v \in \mathbb{R}^2: W_t (v) >0\}
\]
and choose $t_0$ such that for all $t \leq t_0$, $D_t \subset U(0)$. By the  \(C^4_+\) assumption and compactness of \(\partial K\) we can choose the same \(t_0\)
uniformly for all $y \in \partial K$. Define $\delta_0= c_1 \, t_0^2$. We set
$$
F_t(v,s) = (s v, -t + s(t+u(v))) \quad \quad v \in D_t, \quad 0 \leq s \leq 1.
$$
Then
$$
[-t e_3, K] \setminus K= \{F_t(v,s):  v \in \overline {D_t}, 0 \leq s \leq 1\}
$$
and the $e_3$-coordinate $b_3(t)$ of the barycenter $b(t)$  of $[-t e_3, K] \setminus K$ is 
\begin{eqnarray}\label{b3-Lemma5}
\langle g^\delta(x)-y,N_K(y)\rangle &=& b_3(t) =\frac{ \int_{[-t e_3, K] \setminus K} z dz dv}{\int_{[-t e_3, K] \setminus K} dz dv} = 
\frac{\int_{D_t} \int_{s=0}^1 (-t +s(t+u(v)) s^2 W_t(v) ds dv}{\int_{D_t} \int_{s=0}^1 s^2 W_t(v) ds dv} \nonumber \\
&=&
\frac{3 \int_{D_t}(t+u(v)) W_t(v)dv}{4\int_{D_t} W_t(v) dv}  -t.
\end{eqnarray}
In the above we have also used that the Jacobian $\det D F_t(v,s) = s^2 W_t(v)$. We make a change of variable and put 
$$v= t^\frac{1}{2} w, \quad E_t=\{w:\sqrt t\,w\in D_t\}.
$$
Since \(Q\) has degree two and \(P\) has degree three,
\[
Q(\sqrt t\,w)=tQ(w),\qquad
P(\sqrt t\,w)=t^{3/2}P(w).
\]
Thus
\begin{equation*}
u(\sqrt t\,w)
 =t\left[Q(w)+\sqrt tP(w)+O(t)\right],
\end{equation*}
and
\begin{equation*}
W_t(\sqrt t\,w)
 =t\left[1-Q(w)-2\sqrt tP(w)+O(t)\right].
\end{equation*}
Hence
\begin{equation}\label{Wt-Lemma5}
\int_{D_t}W_t\dd v
= t^2\left(\int_{E_t} \left[1-Q(w)-2\sqrt tP(w)+O(t)\right] \dd w \right).
\end{equation}
Let now  $E= \{w: Q(w) <1\}$. 
Now note that the difference between
\(W_t(\sqrt t\,w)/t\) and \(1-Q(w)\) is \(O(\sqrt t)\). Therefore  we have for the symmetric difference $\Delta$, 
\begin{equation*}\label{eq:shell-inclusion}
(E_t \Delta E) 
\subset\left\{w:|Q(w)-1|\le C\sqrt t\right\}
\end{equation*}
and consequently $\vol_2(E_t \Delta E) \leq O(\sqrt t)$. 
On the set $E_t \Delta E$, the maps  \(1-Q\) and \(Q(1-Q)\) are \(O(\sqrt t)\). Thus
\begin{eqnarray*}
\int_{E_t}  \left[1-Q(w)-2\sqrt tP(w) \right]  \dd w
 &=&\int_E  \left[1-Q(w)-2\sqrt tP(w)+O(t)\right] \dd w +O(t) \\
 &= &\int_E  \left[1-Q(w)\right]  \dd w +O(t) + \sqrt{t} \int_E  P(w) \dd w \\
 &=&\int_E  \left[1-Q(w)\right]  \dd w +O(t).
\end{eqnarray*}  
The last identity holds as integration of the odd function $P$ over the symmetric ellipse $E$ is equal to $0$.
Thus, with (\ref{Wt-Lemma5}), 
\begin{eqnarray*}
\int_{D_t}W_t\dd v
&=& t^2 \left( \int_E  \left[1-Q(w)\right]  \dd w +O(t) \right) = t^2\left( \frac{2\pi}{\sqrt{\det D^2 u}}\int_0^1(1-s)\dd s +O(t) \right)\\
 &=& t^2 \left(\frac{\pi}{\sqrt{\det D^2 u}} +O(t) \right).
\end{eqnarray*}
Similarly we show that 
\begin{eqnarray*}
\int_{D_t}u\,W_t\dd v &=& t^3\left[\int_{E} \left[  Q(1-Q) \right]\dd w +O(t) \right] = 
t^3\left( \frac{2\pi}{\sqrt{\det D^2 u}}\int_0^1 s(1-s)\dd s +O(t) \right) \\
 &=& t^3 \left(\frac{\pi}{3 \sqrt{\det D^2 u}} +O(t) \right).
\end{eqnarray*}
Then we get with (\ref{b3-Lemma5}),
\begin{align*}
    b_3(t)&=\frac{3\int_{D_t}(t+u(v))W_t(v)dv}{4\int_{D_t}W_t(v)dv}-t\\
    &=\frac{3}{4}\frac{t^3\bigg[\int_E[1-Q(w)]dw+\int_E[Q(w)(1-Q(w))]dw+O(t)\bigg]}{t^2[\int_E[1-Q(w)]dw+O(t)]}-t\\
    &=\frac{3}{4}\frac{t^3\bigg[\frac{\pi}{\sqrt{\det D^2u}}+\frac{\pi}{3\sqrt{\det D^2u}}+O(t)\bigg]}{t^2[\frac{\pi}{\sqrt{\det D^2u}}+O(t)]}-t
   =\frac{3t}{4}\left[\frac{4}{3}+O(t)\right]-t
   =O(t^2).
\end{align*}
\par
\noindent
As $g^\delta(x) \in [0,x]$, there is $0\leq s <1$ such that  \(g^\delta(x)=sx\). By (\ref{beta2}), we have that $\langle  y, N_K(y)\rangle \geq \beta$ 
and it holds that $\langle x,  N_K(y)\rangle \leq \langle y,  N_K(y)\rangle$. 
Then we get with (\ref{delta-Lemma5}), 
\begin{eqnarray*}
\beta(1-s)
\le \langle y,  N_K(y)\rangle  - \langle sx,N_K(y)\rangle &\leq& \big|\langle y- g^\delta(x),N_K(y)\rangle\big|
= |b_3(t)|\leq c_2\, t^2 \\ 
&\leq& \frac{c_2}{c_1}\, \delta.
\end{eqnarray*}
Therefore
\begin{equation}\label{dim3-innen}
\|g^\delta(x)-x\|
=(1-s)\|x\|
\le\frac{c_2}{\beta^2 c_1}\,\delta.
\end{equation}
\end{proof}

\vskip 3mm
\noindent
The next lemma treats the dimension $3$ case with a different power of $\delta$.
\vskip 3mm
\noindent
\begin{lemma} \label{dim3-bounded}
Let $K$ be a convex body in $\mathbb{R}^3$ that is $C^4_+$ and such that $0 \in  \text{int}(K)$. 
Then there is $\delta_0$ such that  for all $\delta \leq \delta_0$, for all $x \in \partial K$, 
\begin{equation*}\label{bounded:n=3}
 \frac{\langle x,N_K(x)\rangle}{3\,  \delta} \ \left|1-\left(\frac{\|g^\delta(x)\|}{\|x\|
}\right)^{3}\right| \leq  c_3 
\end{equation*}
where  $c_3$ is a constant independent of $\delta$.
\end{lemma}
\vskip 2mm
\noindent
\begin{proof}
Again, we first consider the case that  $g^\delta(x)$ is in the interior of the complement  of $K$,  i.e., 
$g^\delta(x) \in \text{int} (K^C)$.  Then it follows from (\ref{bounded-dim3}), that 
$$
\frac{1}{3 \delta } \, \langle x,N_K(x)\rangle \left|1-\left(\frac{\|g^\delta(x)\|}{\|x\|}\right)^{3}\right| \leq   \frac{c}{4 \pi r_0^2} .
$$
Next we treat  the case that  $g^\delta(x) \in \text{int} (K)$.  Then it follows from
(\ref{dim3-innen}) that 
$$
\frac{1}{3 \delta } \, \langle x,N_K(x)\rangle \left|1-\left(\frac{\|g^\delta(x)\|}{\|x\|}\right)^{3}\right|
\le\frac{c_2}{\beta^2 c_1}.
$$
\end{proof}
\vskip 3mm
\noindent
The concept of the Dupin indicatrix is a useful tool to describe the local shape of the boundary of a convex body. 
By the assumptions  on Theorem \ref{main}, the indicatrix of Dupin is an ellipsoid everywhere on $\partial K$.
Since the quantities considered in Theorem \ref{main} are affine invariant,  we can assume that the indicatrix is a Euclidean ball.
\vskip 2mm
\noindent
The next result can be found in e.g.,  \cite{SchuettWerner1990}, under even weaker assumptions on $K$.
\begin{lemma} \label{Dupin}
\noindent 
Let $K$ be a convex body in $\mathbb{R}^n$ that is $C^2_+$.
Let  $x \in \partial K$ be such that the indicatrix of Dupin exists.
We can assume that the indicatrix of Dupin is a Euclidean ball.   Then locally at $x$, $\partial K$ is approximated by a Euclidean ball, i.e., 
for every  $\varepsilon >0$ there exists $\Delta_\varepsilon>0$ such that for all $\Delta\leq \Delta_\varepsilon$,
\begin{eqnarray*}
&& \hskip -15mm B^n_2\big(x-(1-\varepsilon) r N_K(x), (1-\varepsilon) r \big) \cap H (x-\Delta N_K(x), N_K(x) ) ^-  \\ 
&& \hskip 5mm \subseteq K  \cap  H (x-\Delta N_K(x), N_K(x)) ^-  \\
&& \hskip 5mm \subseteq B^n_2\big(x-(1+\varepsilon) r N_K(x), (1+\varepsilon) r\big) \cap  H (x-\Delta N_K(x), N_K(x)) ^-, 
\end{eqnarray*}
where $r=r(x)= \kappa(x)^{-\frac{1}{n-1}}$.
\end{lemma}

\vskip 3mm
\noindent
The next lemma addresses the case $n\neq 3$.
\vskip 2mm
\begin{lemma} \label{limit}
Let $K$ be a convex body in $\mathbb{R}^n$ that is $C^2_+$. Then   we have for all $x \in \partial K$, 
\[
\frac{1}{n}  \lim _{\delta \rightarrow 0} \frac{\langle x,N_K(x)\rangle \left(1-\left(\frac{\|g^\delta(x)\|}{\|x\|}\right)^{n}\right) } {\delta^{2/(n+1)}} = \alpha_n\  \kappa(x) ^\frac{1}{n+1}, 
\]
where $\alpha _n= \frac{1}{2} \  \frac{n-3}{n+3} \left( \frac{n(n+1)}{\vol_{n-1}(B^{n-1}_2)}\right)^\frac{2}{n+1}$. 
\end{lemma}
\vskip 3mm
\noindent
\begin{remark}
Note that for $n=3$ the right hand side of Lemma \ref{limit} is $0$. Lemma \ref{limit3} shows that 
the left hand side of the expression   in the lemma is  also equal to $0$.
\end{remark}

\vskip 2mm
\noindent
\begin{proof}
Without loss of generality we can assume that $0 \in \text{int} (K)$. 
By the assumptions on $K$, 
the indicatrix exists everywhere on $\partial K$ and  is an ellipsoid.  By affine invariance we can assume that the indicatrix is a Euclidean ball. Let $r=r(x) =  \kappa(x)^{-\frac{1}{n-1}}$ and let  $B^n_2\big(x- r N_K(x)),  r \big)$ be the ball locally approximating $K$, given by Lemma \ref{Dupin}.  Let $ \varepsilon >0$ be given.
By  Lemma \ref{Dupin} there exists $\Delta_\varepsilon>0$ such that for all $\Delta\leq \Delta_\varepsilon$,
\begin{eqnarray*}
&& \hskip -15mm B^n_2\big(x-(1-\varepsilon) r N_K(x), (1-\varepsilon) r \big) \cap H (x-\Delta N_K(x), N_K(x) ) ^- \\
&& \hskip 5mm \subseteq K  \cap  H (x-\Delta N_K(x), N_K(x)) ^- \\
&&\hskip 5mm  \subseteq  B^n_2\big(x-(1+\varepsilon) r N_K(x), (1+\varepsilon) r\big) \cap  H (x-\Delta N_K(x), N_K(x)) ^-.
\end{eqnarray*}
We choose in addition $\Delta  \leq \min\{ \Delta_\varepsilon, (1-\varepsilon) r/2 \}$.
\par
\noindent
In the proof we consider two cases, $g^\delta(x) \in \text{int} (K)$ and $g^\delta(x) \in \text{int} (K^C)$, i.e. $g^\delta(x)$  is in the interior of the complement of $K$.
It can happen that $g^\delta(x) \in \partial K$, see Corollary \ref{Kor1} (ii). Then we will pass to a new $\delta$ that is strictly smaller and then $g^\delta(x) \in \text{int} (K)$
or $g^\delta(x) \in \text{int} (K^C)$.
\vskip 2mm
\noindent
\underline {Case 1:} \hskip 2mm $g^\delta(x) \in \text{int} (K)$. 
\par
\noindent
By Corollary \ref{cor2}, this happens for $n \geq 4$ but it can happen for $n=3$ as well, see Example 3 in Remark \ref{Beispiele}.
Then 
$$
\left|1-\left(\frac{\|g^\delta(x)\|}{\|x\|}\right)^{n}\right| = 1-\left(\frac{\|g^\delta(x)\|}{\|x\|}\right)^{n}.
$$
We have that $\left(\frac{\|g^\delta(x)\|}{\|x\|}\right)^{n} = \left(1- \frac{\|x-g^\delta(x)\|}{\|x\|} \right)^n$,  and thus 
\begin{eqnarray}\label{nx1}
\frac{1}{n}  \langle x,N_K(x)\rangle \left(1-\left|\frac{\|g^\delta(x)\|}{\|x\|}\right|^{n}\right) \leq  \left\langle \frac{x}{\|x\|},N_K(x)\right\rangle  \  \|x-g^\delta(x)\|.
\end{eqnarray}
By Corollary  \ref{cor2} (i),  $g^\delta(x) \rightarrow x$ as $\delta \rightarrow 0$. Therefore we get  for $\delta$ small enough that with a constant $c$, 
\begin{eqnarray}\label{nx2}
\frac{1}{n}  \langle x,N_K(x)\rangle \left(1-\left|\frac{\|g^\delta(x)\|}{\|x\|}\right|^{n}\right) \geq  \left\langle \frac{x}{\|x\|},N_K(x)\right\rangle  \  \|x-g^\delta(x)\| \  \left(1-c\  \frac{\|x-g^\delta(x)\|}{\|x\|}\right).
\end{eqnarray}
By definition of $g^\delta(x)$ there exists $x^\delta$ such that 
$$
 \vol_n\left([x^\delta, K] \setminus K\right) =  \delta  \hskip 3mm \text{and} \hskip 3mm g^\delta(x)= \frac{1}{\delta}  \  \int_{[x^\delta, K] \setminus K}y \   dy.
$$
Moreover  we choose $\delta$ so small that  the following conditions hold:
\begin{eqnarray}  \label{conditions}
&& \left \langle \frac{x}{\|x\|},N(x)\right\rangle  \  \|x-g^\delta(x)\| \leq \Delta/2, \text{  which we can do by Corollary  \ref{cor2} (i)} \nonumber \\
&&\left[x^\delta, B^n_2\big(x-(1-\varepsilon) r N_K(x), (1-\varepsilon) r \big)\right ]  \subset H (x-\Delta N_K(x), N_K(x)) ^-,  \nonumber \\
&&\left[x^\delta, B^n_2\big(x-(1+\varepsilon) r N_K(x), (1+\varepsilon) r \big)\right]   \subset H (x-\Delta N_K(x), N_K(x)) ^-.
\end{eqnarray}
The last two can be achieved as $x^\delta \rightarrow x$ by \cite{Werner94}. 
Then
\begin{eqnarray} \label{delta1}
\delta &=&  \vol_n\left([x^\delta, K] \setminus K\right) \leq \Bigg[  \vol_n\left([x^\delta, B^n_2(x-(1+\varepsilon)r N_K(x),(1+\varepsilon)r )] \cap H_1^- \right)  \nonumber\\
 && \hskip 28mm  - \  
 \vol_n\left( B^n_2(x-(1-\varepsilon)r N_K(x),(1-\varepsilon) r)  \cap H_1^- \right) \Bigg]
\end{eqnarray}
and
\begin{eqnarray} \label{delta2}
\delta &=& \vol_n\left([x^\delta, K] \setminus K\right)  \geq \Bigg[  \vol_n\left([x^\delta, B^n_2(x-(1-\varepsilon)r N_K(x),(1-\varepsilon)r )] \cap H_2^- \right)  \nonumber\\
 &&  \hskip 28mm - \
 \vol_n\left(B^n_2(x-(1+\varepsilon)r N_K(x),(1+\varepsilon)r ) \cap H_2^- \right) \Bigg],
\end{eqnarray}
where $H_1$ is the hyperplane of intersection of  the convex hull of $x^\delta$ with $B^n_2(x-(1+\varepsilon)r N_K(x),(1+\varepsilon) r )$
and $H_2$ is the hyperplane of intersection of  the convex hull of $x^\delta$ with $B^n_2(x-(1-\varepsilon)r N_K(x),(1-\varepsilon) r )$.
This means that up to an error of the order of $\varepsilon$, $x^\delta$ is on the boundary of the illumination body  $\left(B^n_2(x- r N_K(x), r )\right)^{\delta}$.
\par
\noindent
Without loss of generality we assume now that $N_K(x) = e_n= (0, \dots, 0, 1)$ and that the origin is at $x- (1+\varepsilon) r N_K(x)$. Then
\begin{equation}\label{Ausdruck1}
\left\langle \frac{x}{\|x\|},N_K(x)\right\rangle  \  \|x-g^\delta(x)\| = (1+\varepsilon)  r - \langle g^\delta(x), e_n \rangle.
\end{equation}
We write in short $B^n_2((1+\varepsilon)r)= B^n_2(0,(1+\varepsilon)r )$, and $ B^n_2((1-\varepsilon)r)= B^n_2(2\varepsilon r e_n,(1-\varepsilon)r )$, and 
keep the notations $H_1$ and $H_2$ for the hyperplanes  of above. We get
\begin{eqnarray} \label{center1}
 \langle g^\delta(x), e_n \rangle  &=& \frac{1}{\delta}  \  \int_{[x^\delta, K] \setminus K} \langle y, e_n \rangle \  dy = \frac{1}{\delta}  \  \int_{[x^\delta, K] \setminus K}  y_n \  dy \nonumber\\
& \leq &\frac{1} {\delta}  \   
\int_{[x^\delta, B^n_2((1+\varepsilon)r)]\cap H_1^- \setminus (B^n_2((1-\varepsilon)r) \cap H_1^-)}  y_n \  dy.
\end{eqnarray}
Let $z^{\delta} \in \partial \left(B^n_2(x- r N_K(x), r )\right)^{\delta}$ be such that 
$z^{\delta} = \{\lambda e_n: \lambda >0\} \cap \left(B^n_2(x- r N_K(x), r )\right)^{\delta}$.  
As  up to an error of the order of $\varepsilon$, $x^\delta$ is on the boundary of the illumination body  $\left(B^n_2(x- r N(x), r )\right)^{\delta}$, 
we 
then obtain, continuing after (\ref{center1}), that
\begin{eqnarray} \label{center1,1}
 \langle g^\delta(x), e_n \rangle  &\leq&  \frac{1}{\delta}  \  \int_{[z^{\delta}, B^n_2((1+\varepsilon)r)]\cap H_1^- \setminus (B^n_2((1-\varepsilon)r) \cap H_1^-)}  y_n \  dy \nonumber \\
&\leq& \frac{1}{\delta}  \  \int_{[z^{\delta}, B^n_2((1+\varepsilon)r)] \setminus (B^n_2((1+\varepsilon)r)}  y_n \  dy = 
   g_n ((1+  \varepsilon)r, \delta )  \\
  &=&  (1+ \varepsilon) r \Big[ 1- \frac{1}{2} \  \frac{n-3}{n+3} \left( \frac{n(n+1)}{\vol_{n-1}(B^{n-1}_2)}\right)^\frac{2}{n+1}  \  \frac{\delta^\frac{2}{n+1}}{(1 + \varepsilon) r} \ \left((1+\varepsilon) r\right) ^{-\frac{n-1}{n+1} } 
  \nonumber\\
  &+&   c_1(n) \   \delta^\frac{4}{n+1}\Big]\nonumber,
\end{eqnarray}
where in the second last and the last equality we have used Proposition \ref{prop1},  
(\ref{beta}), the definition of $t=\frac{\Delta}{r}$  and Lemma \ref{deltaDelta} to get that 
\begin{eqnarray} \label{r-delta}
g_n(r, \delta) = r \left[ 1- \frac{1}{2} \  \frac{n-3}{n+3} \left( \frac{n(n+1)}{\vol_{n-1}(B^{n-1}_2)}\right)^\frac{2}{n+1}  \  \frac{\delta^\frac{2}{n+1}}{r} \  r^{-\frac{n-1}{n+1} }
+ c_1(n)\   \delta^\frac{4}{n+1} \right].
\end{eqnarray}
We also assume that $n\neq 3$.
\par
\noindent
Hence, by (\ref{Ausdruck1}), 
\begin{eqnarray*}
\left\langle \frac{x}{\|x\|},N_K(x)\right\rangle  \  \|x-g^\delta(x)\|  &\geq& 
\frac{1}{2} \  \frac{n-3}{n+3} \left( \frac{n(n+1)}{\vol_{n-1}(B^{n-1}_2)}\right)^\frac{2}{n+1}  \  \frac{\delta^\frac{2}{n+1}}{((1+\varepsilon) r )^{\frac{n-1}{n+1} }}  \\
&-&(1+ \varepsilon) \, r \, c_1(n)\,   \delta^\frac{4}{n+1}.
\end{eqnarray*}
Therefore, with (\ref{nx2}),
\begin{eqnarray}\label{unten}
&&\hskip -10mm  \frac{1}{n} \frac{ \langle x,N_K(x)\rangle}{\delta^\frac{2}{n+1}} \left(1-\left|\frac{\|g^\delta(x)\|}{\|x\|}\right|^{n}\right) \geq  \nonumber \\
&&\left(1-c\  \frac{\|x-g^\delta(x)\|}{\|x\|}\right) \Bigg[ \frac{1}{(1+\varepsilon)^\frac{n-1}{n+1}}\ 
\frac{n-3}{2(n+3)} \left( \frac{n(n+1)}{\vol_{n-1}(B^{n-1}_2)}\right)^\frac{2}{n+1}  \  \kappa^\frac{1}{n+1}  \nonumber \\
&& \hskip 50mm-(1+ \varepsilon) \, r \, c_1(n)\,   \delta^\frac{2}{n+1}  \Bigg].
\end{eqnarray}
\par
\noindent
Next we consider the estimate from above. 
\par
\noindent
We still assume  that $N_K(x) = e_n= (0, \dots, 0, 1)$  but that now the origin is at $x- (1-\varepsilon) r N_K(x)$. We now write 
 in short $B^n_2 ((1+\varepsilon) r)= B^n_2(-2 \varepsilon r,(1+\varepsilon)r )$, and $ B^n_2((1-\varepsilon)r)= B^n_2(0,(1-\varepsilon)r )$,
Then
\begin{equation}\label{Ausdruck22}
\left\langle \frac{x}{\|x\|},N_K(x)\right\rangle  \  \|x-g^\delta(x)\| = (1-\varepsilon)  r - \langle g^\delta(x), e_n \rangle.
\end{equation}
\begin{eqnarray}\label{center2}
  \langle g^\delta(x), e_n \rangle  &=& \frac{1}{\delta}  \  \int_{[x^\delta, K] \setminus K}y_n \   dy    
\geq \frac{ 1}{\delta}  \    \int_{[x^\delta, B^n_2((1-\varepsilon)r)]\cap H_2^- \setminus (B^n_2((1+\varepsilon)r) \cap H_2^-)}  y_n \  dy  \nonumber\\
&\geq&   \cos \theta \ g_n((1-\varepsilon) r, \delta),
\end{eqnarray}
where $\theta$ is the angle between $x^\delta$ and $e_n$. Hence, 
$$
\cos \theta = \frac{(1-\varepsilon)^2r^2 +\|x^\delta\|^2 - \|x - x^\delta\|^2}{2 (1-\varepsilon) r \|x^\delta\|}.
$$
By  (\ref{delta1}) and (\ref{delta2}),  $x^\delta$ is, up to an error of $\varepsilon$,  on the boundary of the illumination body  $\left(B^n_2(x- (1-\varepsilon) r N(x), r )\right)^{\delta}$.
Therefore, by \cite{JenkinsonWerner, Werner94}
$$
\|x^\delta\| = 
(1-\varepsilon) r\left(1+ \frac{d_n}{ ((1-\varepsilon) r)^\frac{2n}{n+1} } \delta^\frac{2}{n+1} + c \    \delta^\frac{4}{n+1}\right),
$$ 
where $d_n=\frac{1}{2} \left( \frac{n(n+1)}{ \vol_{n-1}(B^{n-1}_2)}\right)^\frac{2}{n+1}$, i.e., $\alpha_n = \frac{n-3}{n+3} \, d_n$.  Thus we get with $c=c(n)$
\begin{eqnarray}\label{costhetaest}
\cos \theta &=& \frac{(1-\varepsilon)^2 r^2 +(1-\varepsilon)^2 (1-\varepsilon)^2r^2\left(1+ \frac{d_n}{ ((1-\varepsilon)r)^\frac{2n}{n+1} } \delta^\frac{2}{n+1} + c \   \delta^\frac{4}{n+1}\right)^2 - \|x - x^\delta\|^2}{2 (1-\varepsilon)^2r^2\left(1+ \frac{d_n}{ ((1-\varepsilon)r)^\frac{2n}{n+1} } \delta^\frac{2}{n+1} + c \   \delta^\frac{4}{n+1}\right)}\nonumber \\
&=& 1- \frac{\|x - x^\delta\|^2}{2 (1-\varepsilon)^2r^2}  + \frac{d_n \   \|x - x^\delta\|^2 \   \delta^\frac{2}{n+1}}{2 ((1-\varepsilon)r)^{\frac{4n+2}{n+1}}} + c \, \delta^\frac{4}{n+1}\nonumber \\
&\geq&  1- \frac{\|x - x^\delta\|^2}{2 (1-\varepsilon)^2r^2}.
\end{eqnarray}
We use  this,  (\ref{center2}), (\ref{r-delta})  and (\ref{Ausdruck22})  and get that with a (new) constant $c$ 
\begin{eqnarray*}
 \left\langle \frac{x}{\|x\|},N_K(x)\right\rangle  \  \|x-g^\delta(x)\|  &\leq&
  \frac{\|x-x^\delta\|^2}{2(1-\varepsilon) r}+ 
  \alpha_n \, \delta^\frac{2}{n+1} ((1-\varepsilon)) r )^{-\frac{n-1}{n+1}}\\
  &+& \delta^\frac{4}{n+1} \left( c(1-\varepsilon) r  +\frac{\|x-x^\delta\|^2}{2(1-\varepsilon) r} \right)
\end{eqnarray*}
Again, we assume that $n\neq 3$. Therefore, with (\ref{nx1}), 
\begin{eqnarray}\label{above22}
\frac{1}{n} \frac{ \langle x,N(x)\rangle} {\delta^\frac{2}{n+1}} \  \left(1-\left|\frac{\|g^\delta(x)\|}{\|x\|}\right|^{n}\right) &\leq & 
  \frac{\|x-x^\delta\|^2}{2(1-\varepsilon) r\delta^\frac{2}{n+1}} + 
\alpha_n   \frac{ \kappa(x) ^{\frac{1}{n+1} }}{ (1-\varepsilon)^\frac{n-1}{n+1} }  \nonumber \\
&+ &
 \delta^\frac{2}{n+1}\left(\frac{\|x - x^\delta\|^2}{2 (1-\varepsilon)r }  +  c (1-\varepsilon) r \right).
\end{eqnarray}
 As $x^\delta \in \partial K^\delta$, we have, using the Blaschke rolling ball (\ref{BlaschkeRolling}), 
\begin{eqnarray*}
\delta &=&
\vol_n\left([x^\delta, K] \setminus K  \right)
 \geq \vol_n\left([x^\delta,  B^n_2\big(x- r_0 N_K(x),  r_0\big)]  \cap  H^-\left(x, N_K(x)\right)\right)\\
 &&\geq  \vol_n\left(C\right),
\end{eqnarray*}
where  $C=[z^\delta,  B^n_2\big(x- r_0 N_K(x),  r_0\big)]  \cap  H^-\left(x, N_K(x)\right)$ is the cone with apex 
$$z^\delta= H(x^\delta, N_K(x)) \cap \left\{ x+ \lambda N_K(x): \lambda \geq 0 \right\}.
$$  
We compute the volume of $C$, where $\alpha$ is the angle between $x^\delta-x$ and $N_K(x)$. For $\delta$ sufficiently small,
\begin{eqnarray}\label{volumeC}
\vol_n\left(C\right) &=&  r_0^{n-1} \frac{ \vol_{n-1}\left(B^{n-1}_2\right)}{n}  \  \frac{(\|x-x^\delta\|\cos\alpha )^\frac{n+1}{2}} {(\|x-x^\delta\|\cos\alpha +2r_0)^\frac{n-1}{2}} \nonumber\\
&\geq&  r_0^{n-1} \frac{ \vol_{n-1}\left(B^{n-1}_2\right)}{n}  \  \frac{(\|x-x^\delta\|\cos\alpha )^\frac{n+1}{2}} {(3r_0)^\frac{n-1}{2}}\nonumber
\end{eqnarray}
and therefore 
\begin{eqnarray*}
\delta^\frac{2}{n+1}  \geq     r_0^\frac{2(n-1)}{n+1} \   \left(\frac{ \vol_{n-1}\left(B^{n-1}_2\right)}{n}\right)^\frac{2}{n+1}   \  \frac{\|x-x^\delta\|\cos\alpha  } {(3r_0)^\frac{n-1}{n+1}} .
\end{eqnarray*}
Equivalently,
\begin{eqnarray*}
\left(\frac{ n}{\vol_{n-1}\left(B^{n-1}_2\right)}\right)^\frac{2}{n+1}\left(\frac{3 }{r_0}\right)^\frac{n-1}{n+1}  \frac{\delta^\frac{2}{n+1}}{\cos\alpha} \geq \|x-x^\delta\|.
\end{eqnarray*}
It follows from this and (\ref{above22}) that
\begin{eqnarray}\label{above2}
& \frac{1}{n} \frac{ \langle x,N(x)\rangle} {\delta^\frac{2}{n+1}} \  \left(1-\left|\frac{\|g^\delta(x)\|}{\|x\|}\right|^{n}\right) \leq  
  \left(\frac{ n}{\vol_{n-1}\left(B^{n-1}_2\right)}\right)^\frac{4}{n+1}\left(\frac{3 }{r_0}\right)^\frac{2(n-1)}{n+1}  \frac{\delta^\frac{2}{n+1}}{\cos^2\alpha}\frac{1}{2(1-\varepsilon) r}
  \nonumber\\
  &+ \kappa(x) ^{\frac{1}{n+1} } \frac{1}{2}\ \frac{n-3}{n+3}\left( \frac{n(n+1)}{\vol_{n-1}(B^{n-1}_2)}\right)^\frac{2}{n+1}
\left( 1-   \frac{\|x - x^\delta\|^2}{2 (1-\varepsilon) ^2 r^2} \right)  +  \left(\frac{\|x - x^\delta\|^2}{2 (1-\varepsilon)r }  +  c (1-\varepsilon) r \right) \delta^\frac{2}{n+1} . \nonumber  \\
\end{eqnarray}
With (\ref{unten}) and (\ref{above22}), we get that for $g^\delta(x) \in \text{int} (K)$,
\[
\frac{1}{n}  \lim _{\delta \rightarrow 0} \frac{\langle x,N(x)\rangle \left(1-\left(\frac{\|g^\delta(x)\|}{\|x\|}\right)^{n}\right) } {\delta^{2/(n+1)}} = \alpha_n\  \kappa(x) ^\frac{1}{n+1}, 
\]
where $\alpha _n= \frac{1}{2} \  \frac{n-3}{n+3} \left( \frac{n(n+1)}{\vol_{n-1}(B^{n-1}_2)}\right)^\frac{2}{n+1}$. 
\vskip 2mm
\noindent
\underline{Case 2}: $g^\delta(x) \in \text{int} (K^C)$. 
\par
\noindent
This happens when $n=2$, but it can happen for $n=3$ as  well, see Example 3 in Remark \ref{Beispiele}.
Therefore, 
$$
\left|1-\left(\frac{\|g^\delta(x)\|}{\|x\|}\right)^{n}\right| = \left(\frac{\|g^\delta(x)\|}{\|x\|}\right)^{n} -1.
$$
We have that $\left|\frac{\|g^\delta(x)\|}{\|x\|}\right|^{n} = \left(\frac{\|x-g^\delta(x)\|}{\|x\|} + 1 \right)^n$,  which we use to estimate 
\begin{eqnarray}\label{2nx1}
\frac{1}{n}  \langle x,N_K(x)\rangle \left(\left|\frac{\|g^\delta(x)\|}{\|x\|}\right|^{n} -1 \right) \geq  \left\langle \frac{x}{\|x\|},N_K(x)\right\rangle  \  \|x-g^\delta(x)\|.
\end{eqnarray}
By Corollary  \ref{cor2} (i),  $g^\delta(x) \rightarrow x$, as $\delta \rightarrow 0$. Therefore we get  for $\delta$ small enough that with a constant $c$, 
\begin{eqnarray}\label{2nx2}
\frac{1}{n}  \langle x,N_K(x)\rangle \left(\left|\frac{\|g^\delta(x)\|}{\|x\|}\right|^{n} -1 \right) \leq  \left\langle \frac{x}{\|x\|},N_K(x)\right\rangle  \  \|x-g^\delta(x)\| \  \left(1+c\  \frac{\|x-g^\delta(x)\|}{\|x\|}\right).
\end{eqnarray}
We again choose $\delta$ so small that the conditions  (\ref{conditions}) hold.  Likewise, the estimates (\ref{delta1}), (\ref{delta2}), (\ref{center1}), and  (\ref{center2}) also hold. Hence,
up to an error of the order of $\varepsilon$, $x^\delta$ is on the boundary of $\left(B^n_2(x- r N_K(x), r )\right)^{\delta}$.  We continue the proof for the case $g^\delta(x) \in \text{int} (K^C)$ in the same  way as in the case $g^\delta(x) \in \text{int} (K)$  and arrive again at
\begin{equation}\label{nicht3}
\frac{1}{n}  \lim _{\delta \rightarrow 0} \frac{\langle x,N(x)\rangle \left(1-\left(\frac{\|g^\delta(x)\|}{\|x\|}\right)^{n}\right) } {\delta^{2/(n+1)}} = \alpha_n\  \kappa(x) ^\frac{1}{n+1}, 
\end{equation}
where $\alpha _n= \frac{1}{2} \  \frac{n-3}{n+3} \left( \frac{n(n+1)}{\vol_{n-1}(B^{n-1}_2)}\right)^\frac{2}{n+1}$. 
\end{proof}
\vskip 3mm
\noindent
Finally, noting that $\alpha_3=0$, we determine the asymptotics  in the case $n=3$. 
\vskip 2mm
\begin{lemma}\label{limit3}
Let $K$ be a convex body in $\mathbb{R}^3$ that is $C^4_+$. For all $x \in \partial K$ the limiting behavior in the case $n=3$ is
\begin{eqnarray} \label{limes:dim3}
\frac{1}{3}  \lim _{\delta \rightarrow 0}  \frac{ \langle x,N_K(x)\rangle} {\delta} \  \left|1-\left(\frac{\|g^\delta(x)\|}{\|x\|}\right)^{3}\right|
=  \frac{\kappa(x)^{1/4}}{12 \pi} |J(x) -9\,  H_{\mathrm{aff}}(x)|.
\end{eqnarray}
\end{lemma}
\vskip 2mm
\noindent
\begin{proof}
\par
\noindent
Let \(K\subset\R^3\) be a convex body of class \(C^4_+\), and fix
\(x\in\partial K\).
When $n=3$ it can happen that $g^\delta(x) \in \text{int} (K)$ or that $g^\delta(x) \in \text{int} (K^C)$.
When $g^\delta(x) \in \text{int} (K)$, 
we have that for $\delta$ small enough 
\begin{eqnarray}\label{Nx1}
&&\hskip -5mm \left\langle \frac{x}{\|x\|},N_K(x)\right\rangle  \  \|x-g^\delta(x)\| \  \left(1-\frac{1}{\beta}\, \|x-g^\delta(x)\| + \frac{\beta^2}{3} \|x-g^\delta(x)\|^2 \right) \nonumber \\
&& \hskip 10mm  \leq \frac{1}{3}  \langle x,N_K(x)\rangle \left(1-\left|\frac{\|g^\delta(x)\|}{\|x\|}\right|^{3}\right) \nonumber \leq \\
&&\hskip -5mm   \left\langle \frac{x}{\|x\|},N_K(x)\right\rangle  \  \|x-g^\delta(x)\| \, \left(1-\beta \, \|x-g^\delta(x)\| + \frac{1}{3 \beta^2} 
  \|x-g^\delta(x)\|^2\right),
\end{eqnarray}
where $\beta$ is as in (\ref{beta2}).
When $g^\delta(x) \in \text{int} (K^C)$, 
we have that for $\delta$ small enough 
\begin{eqnarray}\label{Nx2}
&&\hskip -5mm \left\langle \frac{x}{\|x\|},N_K(x)\right\rangle  \  \|x-g^\delta(x)\| \  \left(1 + \beta\, \|x-g^\delta(x)\| + \frac{\beta^2}{3} \|x-g^\delta(x)\|^2 \right) \nonumber \\
&& \hskip 10mm  \leq \frac{1}{3}  \langle x,N_K(x)\rangle \left(\left|\frac{\|g^\delta(x)\|}{\|x\|}\right|^{3}-1\right) \nonumber \leq \\
&&\hskip -5mm   \left\langle \frac{x}{\|x\|},N_K(x)\right\rangle  \  \|x-g^\delta(x)\| \, \left(1+ \frac{1}{\beta} \, \|x-g^\delta(x)\| + \frac{1}{3 \beta^2} 
  \|x-g^\delta(x)\|^2 \right).
\end{eqnarray}
Let $\phi$ be the angle between $x$ and its outer normal $N_K(x)$. Then $\left\langle \frac{x}{\|x\|},N_K(x)\right\rangle = \cos \phi$. We keep $\phi$  fixed and 
assume without loss of generality that $x=0$ and that the tangent space of $K$  in $x=0$ coincides with $\mathbb{R}^2$. 
We describe $\partial K$ in a neighborhood $B^2_2(0,\rho)$ of $x=0$ by a convex $C^4_+$ function 
$$
f : B^2_2(0,\rho) \rightarrow \mathbb{R^+}
$$
with $f(0)=0$, $N_K(0)=(\grad f(0), -1) = -e_3$ and we can assume that the principal directions of curvature at $x=0$ coincide with $e_1$ and $e_2$. The Hessian of $f$ at $0$ is 
\begin{equation} \label{Hessische}
\nabla^2 f (0) =\begin{pmatrix} \kappa_1 & 0\\
               0 & \kappa_2 
\end{pmatrix},
\end{equation}
where we write in short $\kappa_i =  \kappa_i(0)$,  for the strictly positive principal curvatures $\kappa_i$,  $i=1,2$.
Then $\partial K$  locally around $0$ is described by $(y, f(y))$ where $y=(y_1, y_2)$, with
$$
f(y) = \frac{1}{2} \left( \kappa_1 y_1^2 +\kappa_2 y_2^2 \right) + O(\|y\|^3).
$$
We make a change of variable and put $v_i = \sqrt{\kappa_i} y_i$,  $i=1,2$,  and then, still calling the function describing 
$\partial K$  locally around $0$ by $f$, we now have that $\partial K$ is given by  $(v, f(v))$ where $v=(v_1, v_2)$, with
\begin{equation}\label{v-coordinates}
f(v) = \frac{1}{2} \left(  v_1^2 + v_2^2 \right) + O(\|v\|^3)
\end{equation}
and $\nabla^2 f (0) = I_2$, 
the $2 \times 2$ identity matrix.
Note that this change of variables brings in the Jacobian 
\begin{equation} \label{Jacob}
d y_1 d y_2 d y_3 = \frac{d v_1}{\sqrt{\kappa_1}} \frac{d v_2}{\sqrt{\kappa_2}} dy_3 = \frac{d v_1 d v_2 d y_3} {\sqrt{\kappa}}.
\end{equation}
\par
\noindent
Let $t \in \mathbb {R}$, $t>0$ and consider 
\begin{equation}\label{Kegel}
C_t= [-t e_3, K] \setminus K.
\end{equation}
We want to first compute the corresponding $\| x - g^\delta (x) \|= \|g^\delta(0)\|$ in that case.
We put  $f_1 =\frac{\partial f}{\partial v_1}$, $f_2 =\frac{\partial f}{\partial v_2}$, $D f(v)= (f_1(v), f_2(v))$  and set, denoting the gradient again by $D$, 
$$
W_t(v) = t + f(v) - \langle v, D f(v)\rangle.
$$
Geometrically, $f(v) - \langle v, D f(v)\rangle$ is the point of intersection of the tangent plane to $\partial K$ in $(v, f(v))$. Thus 
$W_t(v)$ is this point  shifted by $-t$  and  \(W_t(v)\) is the signed vertical distance from the tip \((0,0,-t)\) to that intersection  and 
\(W_t(v)\geq 0\) is the set of all those  \((v,f(v))\) that are  ``visible" from the tip $-t e_3$.
We set 
$$
D_t= \{v \in \mathbb{R}^2: W_t(v) >0\}.
$$ We choose
$t$ so small that $D_t \subset B^2_2(0,\rho)$.
We now parametrize the line segments from $-t e_3$ to that part of $\partial K$ that is ``visible" from $-t e_3$ as follows,
\begin{equation}
 F_t(v,s)
=(sv,-t+s(t+f(v))),\qquad 0\leq s \leq 1, \qquad v \in D_t.
\end{equation}
Then 
\[
C_t= [-te_3,K]\setminus K = \left\{ F_t(v,s):v \in\overline{D_t},\ 0\le s\leq 1\right \}
\]
We compute the volume and the center of gravity of this set. We observe that with $D$ denoting the derivative with respect to $v_1,v_2, s$
$$
\det D F_t (v,s) = \det \begin{pmatrix} s & 0 & s f_1(v) \\
               0 & s & s f_2(v)\\ 
v_1 & v_2 & t +f(v)
\end{pmatrix} = s^2 W_t(v).
$$
Thus, also taking into account the Jacobian  (\ref{Jacob}) of the earlier change of variables
\begin{equation} \label{3-volume}
\vol_3(C_t) = \frac{1}{\sqrt{\kappa} } \int_{D_t} \int_{s=0}^1 s^2 W_t(v) ds dv = \frac{1}{3 \sqrt{\kappa}} \int_{D_t} W_t(v)  dv.
\end{equation}
We denote by $b(t)  = (b_1(t), b_2(t) , b_3(t))$, $i=1,2,3$ the barycenter of $C_t$. Then for $i=1,2$,
\begin{equation} \label{baryi}
b _i(t)=\frac{1}{\vol_3(C_t)}\int_{D_t} \int_{s=0}^1 s v_i s^2 W_t(v) ds dv = \frac{3}{4} \frac{\int_{D_t}  v_i W_t(v)  dv}{\int_{D_t}  W_t(v)  dv}
\end{equation} 
and 
\begin{eqnarray} \label{bary3}
b_3(t) &=&\frac{1}{\vol_3(C_t)}\int_{D_t} \int_{s=0}^1 (-(1-s) t  +s f(v)) s^2 W_t(v) ds dv \nonumber \\
&=& -\frac{t}{4} + \frac{3}{4} \frac{\int_{D_t}  f(v) W_t(v)  dv}{\int_{D_t}  W_t(v)  dv}.
\end{eqnarray}
We next expand $f$ into a Taylor series around $0$.
We write  in short  $T_{ijk}=f_{ijk}(0)$, $i,j,k\in\{1,2\}$.  Then 
\begin{eqnarray*}
T[v,v,v] &=& D^3 f(0)[v,v,v]= \sum_{i,j,k=1}^2 T_{ijk} v_iv_jv_k \nonumber\\
&=&  f_{111}(0) v_1^3 + 3 f_{112} (0) v_1^2v_2+3 f_{122}(0) v_1v_2^2 + f_{222}(0)v_2^3
\end{eqnarray*}
and similarly
\begin{eqnarray*}
H[v,v,v,v] &=& D^4 f(0)[v,v,v,v] =\sum_{i,j,k,l=1}^2 H_{ijkl} v_iv_jv_k v_l \\
&=& f_{1111}(0)v_1^4 + 4 f_{1112}(0)\,v_1^3v_2 +6 f_{1122}(0)\,v_1^2v_2^2
+4 f_{1222}(0)\,v_1v_2^3 +f_{2222}(0)\,v_2^4.
\end{eqnarray*}
Then
\begin{eqnarray*}
f(v) = \frac{1}{2} \|v\| + \frac{1}{6} T[v,v,v]  + \frac{1}{24} H[v,v,v,v] + o(\|v\|^4).
\end{eqnarray*}
\par
\noindent
To compute the above integrals we pass to polar coordinates for $v$, 
\begin{equation*}\label{eq:polar}
v=r \theta,\qquad \theta \in S^1, 
\qquad \dd v=r\dd r\dd\theta.
\end{equation*}
Thus by Taylor's theorem we have  uniformly in \(\theta\),
\begin{equation*}
f(r \theta) = \frac{1}{2} r^2 + \frac{r^3}{6} D^3 f(0)[\theta,\theta,\theta]+ \frac{r^4}{24} D^4 f (0)[\theta,\theta,\theta,\theta] +o(r^4).
\end{equation*}
Abbreviate the  Taylor coefficients by
\begin{equation}\label{eq:alpha-beta}
\alpha(\theta)=\frac16D^3f(0)[\theta, \theta, \theta],\qquad
\beta(\theta)=\frac1{24}D^4f(0)[\theta,\theta, \theta,\theta].
\end{equation}
 So 
\begin{equation}\label{eq:u-polar}
f(r \theta)
=\frac12r^2+\alpha(\theta)r^3+\beta(\theta)r^4+o(r^4).
\end{equation} 
By chain rule 
$$
\frac{\partial}{\partial r} f(r \theta) = f_1(r\theta) \theta_1 +  f_2(r\theta) \theta_2 = \langle \theta, Df(r\theta) \rangle,  
$$
and hence  $r \frac{\partial}{\partial r} f(r \theta) = \langle r\theta, Df(r\theta) \rangle$. 
Consequently,
\begin{equation}\label{eq:W-polar}
W_t(r \theta) =  t + f(r \theta) - \langle r \theta, Df(r \theta) \rangle 
= t-\frac12r^2-2\alpha(\theta)r^3-3\beta(\theta)r^4+o(r^4).
\end{equation}
We determine the domain of integration of the above integrals.
For each \(\theta\), let \(R_t(\theta)>0\) be such that 
\(W_t(R_t \theta)=0\). Strict convexity of $f$  implies
\(\partial_r W_t(r \theta)=-r\partial_{rr}f(r \theta)<0\)
for small positive \(r\), so the small positive root is unique and
\[
D_t=\{r \theta:0\le r<R_t(\theta),\ 0\le\theta<2\pi\}.
\]
Solving \eqref{eq:W-polar} for that root gives
\begin{equation}\label{eq:radius}
R_t(\theta)=\sqrt{2t}-4\alpha(\theta)t
+\sqrt2\bigl(20\alpha(\theta)^2-6\beta(\theta)\bigr)t^{3/2}
+o(t^{3/2}).
\end{equation}
Substituting (\ref{eq:u-polar}), \eqref{eq:W-polar}, and
\eqref{eq:radius} into the three integrals (\ref{3-volume}), (\ref{baryi}) and  (\ref{bary3}) and integrating over
\(0\le r\le R_t(\theta)\), gives
\begin{align}
&\int_{D_t}W_t(v) dv =\int_{\theta =0}^{2 \pi} \int _{r=0}^{R_t(\theta)} W_t(r \theta) rdr d\theta
=\pi t^2+t^3\int_0^{2\pi}(16\alpha^2-4\beta)\dd\theta
+o(t^3),\label{eq:I0-angular} 
\\
&\int_{D_t} v W_t(v) dv = \int_{\theta =0}^{2 \pi} \int _{r=0}^{R_t(\theta)} r \theta W_t(r \theta) rdr d\theta
=-\frac83t^3\int_0^{2\pi}\theta\alpha\dd\theta
+o(t^3),\label{eq:Iv-angular}
\\
&\int_{D_t}f(v) W_t(v) dv = \int_{\theta =0}^{2 \pi} \int _{r=0}^{R_t(\theta)} f(r \theta) W_t(r \theta) rdr d\theta
=\frac\pi3t^3+t^4\int_0^{2\pi}
\left(12\alpha^2-\frac83\beta\right)\dd\theta
+o(t^4).\label{eq:Iu-angular}
\end{align}
We have $\theta = (\theta_1, \theta_2) =(\cos\theta, \sin\theta)$ and compute 
\begin{eqnarray*}
\int_0^{2\pi} \theta\,  \alpha(\theta) \dd\theta &=& \frac16 \int_0^{2\pi} \theta\,  D^3f(0)[\theta, \theta, \theta] \dd\theta \\
&=& \frac16 \int_0^{2\pi} \theta 
\bigl(T_{111}\cos^3\theta
+3T_{112}\cos^2\theta\sin\theta
+3T_{122}\cos\theta\sin^2\theta+T_{222}\sin^3\theta\bigr) \dd \theta\\
&=&\frac\pi8 (T_{111}+T_{122}, T_{112}+T_{222}) = \frac\pi8 (f_{111}(0)+f_{122}(0), f_{112}(0)+f_{222}(0))\\
&=&\frac\pi8 \tr_{1,2} D^3f (0), 
\end{eqnarray*}
i.e., it is   a contraction of  the third derivative tensor $T=D^3f(0)$ over two indices, i.e. the  partial trace of $\tr_{1,2} D^3f (0)$ of $D^3f(0)$.
\begin{eqnarray*}
\int_0^{2\pi}\beta(\theta) \dd\theta &=& \frac1{24}\int_0^{2\pi} D^4f(0)[\theta,\theta, \theta,\theta] \\
&=&\frac1{24}\int_0^{2\pi} (f_{1111}(0)\cos^4 \theta + 4 f_{1112}(0)\,\cos^3\theta \sin \theta +6 f_{1122}(0)\, \cos^2\theta \sin^2\theta \\
&+&4 f_{1222}(0)\,\cos\theta \sin^3 \theta +f_{2222}(0)\, \sin^4\theta )\, \dd\theta\\
&=&\frac\pi{32} (f_{1111}(0)+2f_{1122}(0)+f_{2222}(0))\\
&=&\frac\pi{32} (\Delta_L^2 f)(0), 
\end{eqnarray*}
where $\Delta_L$ is the Laplacian and $\Delta_L^2  = \Delta_L \circ \Delta_L$ is the composition.
\begin{eqnarray*}
\int_0^{2\pi}\alpha(\theta)^2 \dd\theta &=& \frac{1}{36}  \int_0^{2\pi} \left( D^3f(0)[\theta, \theta, \theta] \right)^2 \dd\theta \\ 
&=&  \frac{1}{36}  \int_0^{2\pi} \bigl(T_{111}\cos^3\theta
+3T_{112}\cos^2\theta\sin\theta
+3T_{122}\cos\theta\sin^2\theta+T_{222}\sin^3\theta\bigr)^2 \dd \theta\\ 
&=& \frac\pi{288}\bigl(5 T_{111}^2+9T_{112}^2+9T_{122}^2+5T_{222}^2 +6T_{111}T_{122}+6T_{112}T_{222}\bigr)\\
&=&\frac\pi{144}\|D^3 f(0)\|_{F}^2
+\frac\pi{96} \| \tr_{1,2} D^3f (0)\|^2, 
\end{eqnarray*}
where 
\begin{equation}\label{Frobenius}
\|D^3 f(0)\|_{F} = \norm T_{F} =\left(\sum _{i,j,k=1}^2 T_{ijk}^2 \right)^\frac{1}{2} = \left(T_{111}^2 +3 T_{112}^2 + 3 T_{122}^2 +T_{222}^2\right)^\frac{1}{2}
\end{equation}
is the Frobenius norm or Hilbert Schmidt norm of the third derivative tensor $T=D^3 f(0)$.
We sometimes write in short 
\begin{equation}\label{PT}
\tau = \tr_{1,2} D^3f (0).
\end{equation}
Thus it follows with (\ref{eq:I0-angular}), (\ref{eq:Iv-angular}) and (\ref{eq:Iu-angular}) that 
\begin{equation} \label{Wt}
\int_{D_t}W_t(v) dv=  \pi t^2+\pi t^3\left(
-\frac {(\Delta_L^2 f)(0)}{8}+\frac19\norm T_{F}^2+\frac16\norm\tau^2
\right) +o(t^3) =  \pi t^2 +o(t^2)
\end{equation}
and consequently
\begin{equation}\label{volCt}
\vol_3(C_t) = \frac{1}{3 \sqrt{\kappa}} \left(\pi t^2+\pi t^3\left(
-\frac {(\Delta_L^2 f)(0)}{8}+\frac19\norm T_{F}^2+\frac16\norm\tau^2
\right)\right) +o(t^3) =  \frac{\pi t^2}{3 \sqrt{\kappa}} +o(t^2).
\end{equation}
By (\ref{baryi}) and (\ref{Wt}), the first and second coordinate of $b(t)$ are, with $i=1,2$
\begin{equation*}
b_i (t) = \frac{3}{4} \frac{\int_{D_t}  v_i W_t(v)  dv}{\int_{D_t}  W_t(v)  dv} = -\frac{t}{4} \tr_{1,2} D^3f (0) +o(t) ,\label{eq:Iv-jets}
\end{equation*}
and by (\ref{bary3}) and (\ref{Wt}), the third coordinate of $g^\delta(0)$ is
\begin{equation*}
b_3(t)=  -\frac{t}{4} + \frac{3}{4} \frac{\int_{D_t}  f(v) W_t(v)  dv}{\int_{D_t}  W_t(v)  dv}=
t^2\left(-\frac{(\Delta_L^2 f)(0)}{32}
+\frac{5}{144}\norm T_{ F}^2
+\frac5{96}\norm\tau^2\right)+o(t^2).
\end{equation*}
We passed to polar coordinates to compute the integrals. Now we go back to the $v$-variable.
We denote by  \(d_K\)  the signed distance, positive inside \(K\) from a point $z$ to $K$
\[
d_K(z)=
\begin{cases}
\operatorname{dist}(z,\partial K),&z\in K,\\
-\operatorname{dist}(z,\partial K),&z\notin K.
\end{cases}
\]
We have to consider the signed distance as $b(t)$ can be in the interior of $K$ or in the interior of the complement of $K$.
We  compute the centroid’s signed distance from the boundary.  In our fixed $v$-coordinates, its second-order expansion is
\begin{eqnarray*}
d_K(b(t)) &=& b_3(t)-\frac12\bigl(b_1(t)^2+b_2(t)^2\bigr)+o(t^2)\\
&= &\frac{t^2}{288} \left(-9  (\Delta_L^2 f)(0)  +10 \|T\|_{F}^2+6\|\tau\|^2\right) +o(t^2).
\end{eqnarray*}
When \(\vol_3(C_t) =\delta\), (\ref{volCt}) implies
\[
t^2=\frac{3\sqrt\kappa}{\pi}\delta+o(\delta).
\]
Consequently,
\begin{equation*}
d_K(b(t))=\frac{\sqrt\kappa}{96\pi}
\left(-9 (\Delta_L^2 f)(0)+10\|T\|_{F}^2+6\|\tau\|^2\right)\,\delta+o(\delta).
\end{equation*}
The calculation has uniform remainders for boundary footpoints near \(0\), and its coefficient depends continuously on the footpoint. The footpoint of the tip generating the radial point \(g^\delta(0)\) tends to \(0\). Therefore
\[
d_K(g^\delta(0))=\left(-9 (\Delta_L^2 f)(0)+10\|T\|_{F}^2+6\|\tau\|^2\right)\, \delta+o(\delta).
\]
Finally, along the fixed radial line through \(0\), normal displacement and radial displacement satisfy
\[
| d_K(g^\delta(0))| =  \frac{\langle x,N_K(x)\rangle}{\|x\|} \| g^\delta(0)\|+O(\|g^\delta(0)\|^2).
\]
Since \(\frac{\langle x,N_K(x)\rangle}{\|x\|}>0\), we get
\[
\frac{\langle x,N_K(x)\rangle}{\|x\|} \|g^\delta(0)\|=
\frac{\sqrt\kappa}{96\pi}
\left|-9 (\Delta_L^2 f)(0) +10\|D^3 f(0)\|_{F}^2+6\|\tr_{1,2} D^3 f(0) \|^2\right|
\delta+o(\delta).
\]
\vskip 2mm
\noindent
We show next that 
\begin{equation} \label{Pick+Haff}
\left|-9 (\Delta_L^2 f)(0) +10\|D^3 f(0)\|_{F}^2+6\|\tr_{1,2} D^3 f(0) \|^2\right| = \frac{8}{\kappa(0)^\frac{1}{4}} |J(0) -9 H_{aff}(0)|.
\end{equation}
In order to do that, we describe,  as above in (\ref{v-coordinates}), 
$\partial K$  locally around $0$ by  $(v, f(v))$ where for $v=(v_1, v_2)$ 
\begin{equation*}
f(v) = \frac{1}{2} \left(  v_1^2 + v_2^2 \right) + O(\|v\|^3)
\end{equation*}
and $D^2 f (0) = I_2$, 
the $2 \times 2$ identity matrix. To shorten notation in what follows we will write  
 \[A_0=D^2f(0), \quad 
A=D^2f(v), \quad (a^{ij})=(a^{ij})_{1 \leq i,j\leq (n-1)}(v) =A(v) ^{-1}.\] 
We also put
\[
\phi(v)=\kappa_K(0)^{1/4}\bigl(\det D^2f(v)\bigr)^{1/4}.
\]
The affine mean curvature \(2 H_{\mathrm{aff}}\) is the trace of the shape operator. In local coordinates this is expressed as follows, see e.g. \cite{Dillenetal}.
The shape operator is the matrix $S_{\mathrm{aff}}$ with entires in row $k$ and column $i$, $1 \leq k,i \leq 2$,
\[
[S_{\mathrm{aff}}]_{ki}(v)=
\sum_{j=1}^2\frac{\partial}{\partial v_i} \!\left(a^{kj}\frac{\partial}{\partial v_j}\phi\right)(v)
\]
and thus 
\[
[S_{\mathrm{aff}}]_{ki}(0)
=\frac{\partial}{\partial v_i} \, \frac{\partial}{\partial v_k} \phi(v) \big|_{v=0}
-\sum_{j=1}^2 f_{ikj}(0)\,\frac{\partial}{\partial v_j} \phi(v) \big |_{v=0}.
\]
Recalling that  \(T=D^3f(0)\), \(\tau_i=\sum_j T_{ijj}\), and the definition of $\phi$, we get
\[
H_{\mathrm{aff}} (0)
=\frac{\kappa(0)^{1/4}}8
\left(\Delta_L^2f(0)-\|T\|_{\mathrm F}^2-\frac34\|\tau\|^2\right).
\]
For the cubic tensor \(\mathcal C=\mathcal C_0\), choose an \(h_0\)-orthonormal basis \(e_1, e_{2}\) of \(T_0(\partial K)\). Then, independent of the chosen \(h_0\)-orthonormal basis, 
\[
\|\mathcal C\|_{h_0}
=\left(\sum_{i,j,k=1}^{2}\mathcal C(e_i,e_j,e_k)^2\right)^{1/2} = \left(\sum_{i,j,k=1}^{2}\mathcal C_{ijk}^2\right)^{1/2}
\]
with
\[
\mathcal C_{ijk}
=\phi(0)^{-1}\left[
T_{ijk}-\frac14
\bigl(\tau_i\delta_{jk}+\tau_j\delta_{ik}+\tau_k\delta_{ij}\bigr)
\right].
\]
The norm $\| \cdot\|_{h_{0}}$  is the same as the Frobenius norm $\| \cdot\|_{F}$. Since \(h^{-1}(0)=\phi(0)I\), 
\[
J=\frac{\kappa(0)^{1/4}}8
\left(\|T\|_{\mathrm F}^2-\frac34 \| \tau\|^2\right).
\]
Subtracting nine times the formula for \(H_{\mathrm{aff}}\) from the formula for \(J\) now gives the signed identity
\[
\frac8{\kappa(0)^{1/4}}(J-9H_{\mathrm{aff}})
=-9\Delta_L^2f(0)+10\|T\|_{\mathrm F}^2+6\|\tau\|^2.
\]
Taking absolute values gives (\ref{Pick+Haff}).

\end{proof}
\vskip 5mm
\noindent
The proof of Theorem \ref{main}  and Theorem \ref{main3} follow immediately from these lemmas.
\par
\begin{proof}[Proof of Theorem \ref{main}]
By the assumptions of the theorem, Lemmas \ref{difference}, \ref{bounded} and Lebesgue's Dominated Convergence Theorem we get that 
\begin{eqnarray*}
		\lim_{\delta \to 0} \frac{\vol_{n}(K) - \vol_{n}(I^\delta(K))}{\delta^\frac{2}{n+1}} &=&  
		\frac{1}{n} \int_{\partial K} \lim _{\delta \rightarrow 0} \frac{\langle x,N_K(x)\rangle \left(1-\left|\frac{\|g^\delta(x)\|}{\|x\|}\right|^{n}\right) } {\delta^{2/(n+1)}} d\mu_K (x).\\
\end{eqnarray*}
Lemma \ref{limit}  finishes the proof. 
We want to remark that by Lemma \ref{bounded} integration and limit can also be interchanged in the case $n=3$ when we consider $\delta^\frac{1}{2}$ 
and thus by Lemma \ref{limit} the left hand side of the theorem is also $0$ in this case.
\end{proof}
\vskip 3mm
\noindent
\begin{proof}[Proof of Theorem \ref{main3}]
The proof of Theorem \ref{main3} is the same. Now we use Lemmas \ref{difference},  \ref{dim3-bounded},  Lebesgue's Dominated Convergence Theorem
and Lemma \ref{limit3}.
\end{proof}

\vskip 5mm
\section{Appendix}

\subsection{Proof of Example \ref{counterexampleparallel}}

\begin{proof}
Let $F_1$ be the facet of $P$ with normal vector $N_1=\langle 0,0,1 \rangle$,
$$F_1=[\{(0,0,0), (2,0,0),(0,2,0),(2,2,0)\}]$$
and let $F_2$ be the facet of $P$ with normal vector $N_2=\langle 0,-1,0 \rangle$,
$$F_2=[\{(0,0,0), (2,0,0),(0,0,-1),(6,0,-1)\}].$$

\begin{figure}[ht]
\centering
    \includegraphics[width=5in]{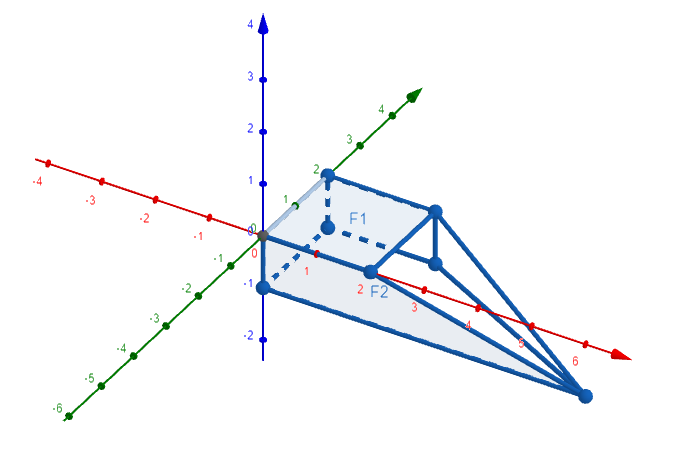}
    \caption{The polytope $P$ in Example \ref{counterexampleparallel} is shown with facets $F_1$ and $F_2$ highlighted.}
    \label{figcounterexampleparallel}
\end{figure}

Take $\delta$ small enough so each vertex of $P^\delta$ sees only one facet of $P$. We consider the facet of $P^\delta$ for which a point $x^\delta\in \partial P^\delta$ sees only $F_1$ and $F_2$. We look at the corresponding facet of $I^\delta(P)$. A normal to this facet of $P^\delta$ is given by
$$N=\vol_{n-1}(F_1)N_1+\vol_{n-1}(F_2)N_2=\langle 0, -4, 4\rangle.$$

To compute a normal $N^*$ to the corresponding facet of $I^\delta(P)$, we compute three vertices and the cross product of the vectors joining them.

We take the vertex $(0,0,0)\in P$. Then there is a facet of $P^\delta$ parallel to $F_1$ at a distance $\frac{3\delta}{4}$, with vertex $(0,0,\frac{3\delta}{4})$. The corresponding vertex of $I^\delta(P)$ is given by
$$v_1=\frac{(0,0,\frac{3\delta}{4})}{n+1}+\left(1-\frac{1}{n+1}\right)g(F_1)=\left(0,0,\frac{3\delta}{16}\right)+\frac{3}{4}(1,1,0)=\left(\frac{3}{4},\frac{3}{4},\frac{3\delta}{16}\right).$$

Similarly, there is a facet of $P^\delta$ parallel to $F_2$ at a distance $\frac{3\delta}{4}$, with vertex $(0,-\frac{3\delta}{4},0)$. The corresponding vertex of $I^\delta(P)$ is given by:
$$v_2=\frac{(0,-\frac{3\delta}{4},0)}{n+1}+\left(1-\frac{1}{n+1}\right)g(F_2)$$
$$=\left(0,-\frac{3\delta}{16},0\right)+\frac{3}{4}\left(\frac{13}{6},0,-\frac{7}{12}\right)=\left(\frac{39}{24},-\frac{3\delta}{16}, -\frac{7}{16}\right).$$
 $F_2$ is a trapezoid, so we compute its center of gravity by decomposing it into a $2\times 1$ rectangle and a triangle with base 4 and height 1, $$g(F_2)=\frac{2(1,0,-\frac{1}{2})+2(\frac{10}{3},0,-\frac{2}{3})}{4}=\left(\frac{13}{6},0,-\frac{7}{12}\right).$$

Next, consider the vertex $(2,0,0)\in P$. There is a facet of $P^\delta$ parallel to $F_1$ at a distance $\frac{3\delta}{4}$. We compute the intersection of the hyperplane containing $F_1$ and the neighboring hyperplanes using similar triangles. Consider the same triangle with base 4 and height 1 as above,
$$\frac{1+\frac{3\delta}{4}}{1}=\frac{4+x}{4}$$
$$\implies x=3\delta$$
and thus the vertex of $P^\delta$ is $(2-3\delta,0,\frac{3\delta}{4})$. Then the corresponding vertex of $I^\delta(P)$ is given by
$$v_3=\frac{(2-3\delta,0,\frac{3\delta}{4})}{4}+\frac{3}{4}(1,1,0)=\left(\frac{5}{4}-\frac{3\delta}{4},\frac{3}{4},\frac{3\delta}{16}\right).$$

Lastly, we compute two vectors lying in the facet of $I^\delta(P)$ under consideration and compute their cross product,

$$t_1=v_2-v_1=\left(\frac{7}{8}, -\frac{3}{4}-\frac{3\delta}{16},-\frac{7}{16}-\frac{3\delta}{16}\right)$$

$$t_2=v_3-v_1=\left( \frac{1}{2}-\frac{3\delta}{4},0,0 \right)$$

$$N^*=\left( 0, \left(-\frac{7}{16}-\frac{3\delta}{16}\right)\left(\frac{1}{2}-\frac{3\delta}{4}\right), \left(\frac{3}{4}+\frac{3\delta}{16}\right)\left(\frac{1}{2}-\frac{3\delta}{4}\right) \right).$$
We factor out common terms to find a normal  
$$N^*=(0,-7-3\delta, 12+3\delta  ).$$
Note, the normal $N^*$ depends on $\delta$ and is not parallel to the normal $N$ to $P^\delta$.
\end{proof}

\vskip 3mm
\noindent
Furthermore, we do not have convexity in general.

\subsection{Counterexamples to Convexity}

\begin{example}\label{counterexample1convex}
    Let $P\subset \R^3$ be the simplex with a vertex cut off defined as the convex hull:
    $$P=[(0,-2,0), (2,0,-2), (-2,0,-2), (0,1,0), (1,1,-1), (-1,1,-1)]$$
    $I^\delta (P)$ is not convex for $\delta<\frac{1}{3}$.
\end{example}

\begin{figure}[ht]
\centering
    \includegraphics[width=5in]{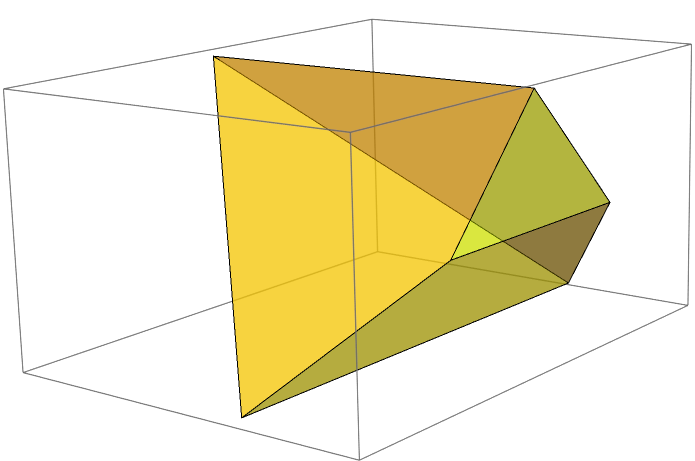}
    \caption{The polytope $P$ in Example \ref{counterexample1convex}.}
    \label{figcounterexample1convex}
\end{figure}

\begin{proof}
    In the case $\delta$ is sufficiently small, note that for each facet $F_i\subset \partial P$, there exists $x\in \partial P^\delta$ such that $x$ sees only $F_i$ and each vertex of $P^\delta$ lies in these facets on the line, $aff(F_j\cap F_k)$, defined by the intersection of facets of $P$.
    
    Let $v_1\in\partial P^\delta$ be the vertex on the line $l(t)=(0,3t,0)$ which adds on $\delta$ to the facet $F_1=[(0,1,0),(1,1,-1), (-1,1,-1)]$. The area of $F_1$ is 1 and the height of the pyramid with base $F_1$ and vertex $v_1$ is $\langle(0,3t-1,0),(0,1,0)\rangle=3t-1$. Then
    $$\delta=\frac{1}{3}(1)(3t-1)\implies t=\frac{3\delta+1}{3}\implies v_1=(0,3\delta+1,0).$$
    Note, the planes containing the faces of $P$ intersect at $(0,2,0)$, thus this description only applies for $3\delta+1<2\implies \delta<\frac{1}{3}$.
    
    Let $v_2\in\partial P^\delta$ be the vertex on the line $l(t)=(1-t,1,t-1)$ which adds on $\delta$ to the facet $F_2=[(0,1,0), (-1,1,-1), (0,-2,0),(-2,0,-2)]$. The area of $F_2$ is $\frac{1}{2}(4)(\sqrt{8})-\frac{1}{2}(1)(\sqrt{2})=\frac{7\sqrt{2}}{2}$ and the height of the pyramid with base $F_2$ and vertex $v_2$ is $\langle(1-t,0,t-1),(\frac{-1}{\sqrt{2}},0,\frac{1}{\sqrt{2}})\rangle=\sqrt{2}(t-1)$. Then
    $$\delta=\frac{1}{3}\left(\frac{7\sqrt{2}}{2}\right)(\sqrt{2}(t-1))\implies t=\frac{3\delta}{7}+1\implies v_2=\left(-\frac{3\delta}{7},1,\frac{3\delta}{7}\right).$$
    Let $v_3\in\partial P^\delta$ be the vertex which adds on $\delta$ to the facet $F_3=[(0,1,0), (1,1,-1), (0,-2,0),(2,0,-2)]$. By symmetry,
    $$v_3=\left(\frac{3\delta}{7},1,\frac{3\delta}{7}\right).$$
    Let $v_4\in\partial P^\delta$ be the vertex on the line $l(t)=(2-2t,-2t,-2+2t)$ which adds on $\delta$ to the facet $F_2$. The height of the pyramid with base $F_2$ and vertex $v_4$ is $\langle(2-2t,2-2t,-2+2t),(\frac{-1}{\sqrt{2}},0,\frac{1}{\sqrt{2}})\rangle=2\sqrt{2}(t-1)$. Then
    $$\delta=\frac{1}{3}\left(\frac{7\sqrt{2}}{2}\right)(2\sqrt{2}(t-1))\implies t=\frac{3\delta}{14}+1\implies v_4=\left(-\frac{3\delta}{7},-\frac{3\delta}{7}-2,\frac{3\delta}{7}\right).$$    
    Let $v_5\in\partial P^\delta$ be the vertex which adds on $\delta$ to the facet $F_3$. By symmetry,
    $$v_5=\left(\frac{3\delta}{7},-\frac{3\delta}{7}-2,\frac{3\delta}{7}\right).$$   
    Let $v_6\in\partial P^\delta$ be the vertex on the line $l(t)=(2-4t,0,-2)$ which adds on $\delta$ to the facet $F_2$. The height of the pyramid with base $F_2$ and vertex $v_6$ is $\langle(4-4t,0,0),(\frac{-1}{\sqrt{2}},0,\frac{1}{\sqrt{2}})\rangle=2\sqrt{2} (t-1)$. Then
    $$\delta=\frac{1}{3}\left(\frac{7\sqrt{2}}{2}\right)(2\sqrt{2} (t-1))\implies t=\frac{3\delta}{14}+1\implies v_6=\left(-\frac{6\delta}{7}-2,0,-2\right).$$

    The vertices $v_2, v_3,v_4$ and $v_5$ define the facet of $P^\delta$ containing all $x$ that see only $F_2\cup F_3$. We compute the corresponding vertices $T(v_2), T(v_3),T(v_4), T(v_5)\in \partial I^\delta (P)$, which define the corresponding facet, $G$, of $I^\delta(P)$. We also compute $T(v_1)$ and $T(v_6)$ and show they lie on opposite sides of the tangent plane $aff(G)$.

    Note in the case $x$ sees only one facet $F_i$, the map $T$ simplifies to $T(x)=\frac{x}{n+1}+\frac{n}{n+1}g(F_i)$.

    Since $F_1$ is a triangle, averaging the vertices yields
    $$g(F_1)=\left(0,1,-\frac{2}{3}\right).$$

    We compute the center of gravity of $F_2$ as the weighted difference of two triangles,
    $$g(F_2)=\frac{4\sqrt{2}(-\frac{2}{3},0,-\frac{2}{3})-\frac{\sqrt{2}}{2}(-\frac{1}{3},\frac{4}{3},-\frac{1}{3})}{\frac{7\sqrt{2}}{2}}=\left(-\frac{15}{21},\frac{4}{21},-\frac{15}{21}\right).$$

    Thus we get, using symmetry to compute $T(v_3)$ and $T(v_5)$,
    $$T(v_1)=\frac{(0,3\delta+1,0)}{4}+\frac{3}{4}(0,1,-\frac{2}{3})=\left(0,\frac{3\delta}{4}+1,-\frac{1}{2}\right)$$
    $$T(v_2)=\frac{1}{4}\left(-\frac{3\delta}{7},1,\frac{3\delta}{7}\right)+\frac{3}{4}\left(-\frac{15}{21},\frac{4}{21},-\frac{15}{21}\right)=\left(-\frac{3\delta+15}{28},\frac{8}{7},\frac{3\delta-15}{28} \right)$$
    $$T(v_3)=\left(\frac{3\delta+15}{28},\frac{8}{7},\frac{3\delta-15}{28} \right)$$    
    $$T(v_4)=\frac{1}{4}\left(-\frac{3\delta}{7},-\frac{3\delta}{7}-2,\frac{3\delta}{7}\right)+\frac{3}{4}\left(-\frac{15}{21},\frac{4}{21},-\frac{15}{21}\right)=\left(-\frac{3\delta+15}{28},-\frac{3\delta+10}{28},\frac{3\delta-15}{28} \right)$$
    $$T(v_5)=\left(\frac{3\delta+15}{28},-\frac{3\delta+10}{28},\frac{3\delta-15}{28} \right)$$
    $$T(v_6)=\frac{1}{4}\left(-\frac{6\delta}{7}-2,0,-2\right)+\frac{3}{4}\left(-\frac{15}{21},\frac{4}{21},-\frac{15}{21}\right)=\left(-\frac{6\delta+29}{28},\frac{1}{7},-\frac{29}{28} \right)$$
   The outer unit normal to the facet, $G$, defined by $T(v_2), T(v_3),T(v_4),$ and $T(v_5)$ is $N=(0,0,1)$. Considering $T(v_1)$ we want to show $-\frac{1}{2}>\frac{3\delta-15}{28}$, which holds for $\delta<\frac{1}{3}$. On the other hand, considering $T(v_6)$, we want to show $-\frac{29}{28}<\frac{3\delta-15}{28}\iff \delta>-\frac{14}{3}$ but we have $\delta>0$.
\end{proof}

\begin{example}
    Let $P\subset \R^3$ be the truncated cube with vertices that are all permutations of $(\pm 2,\pm 2,\pm 1)$. Then $I^\delta(P)$ is not convex for $\delta<44/3$.
\end{example}

\begin{figure}[ht]
    \centering
    \includegraphics[width=.3\textwidth]{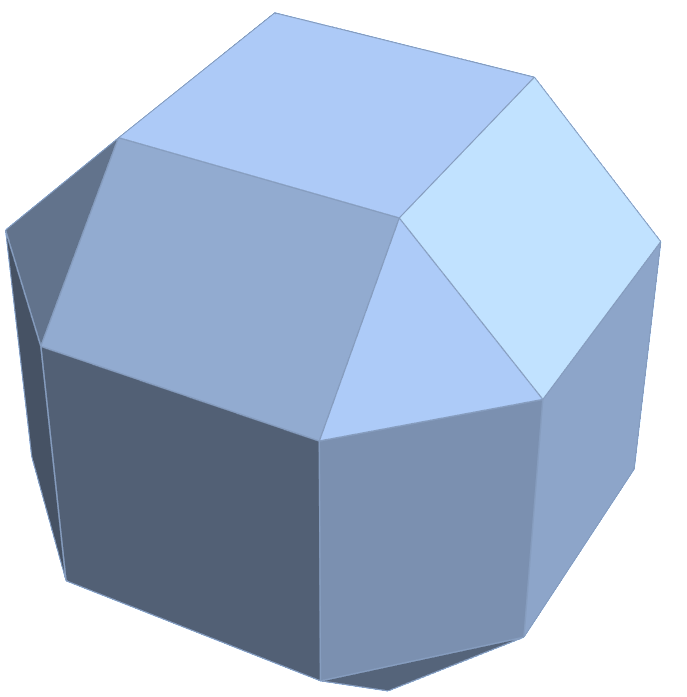}
    \hfill
    \includegraphics[width=.3\textwidth]{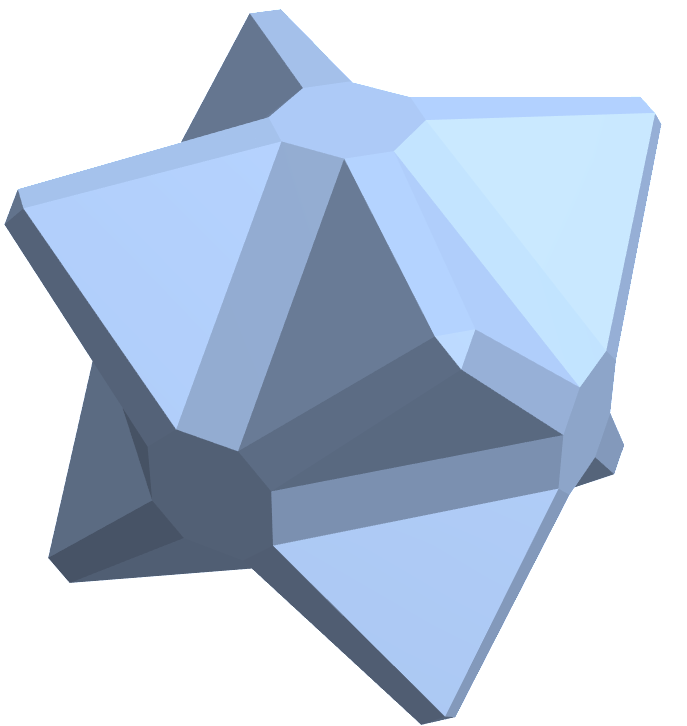}
    \hfill
    \includegraphics[width=.3\textwidth]{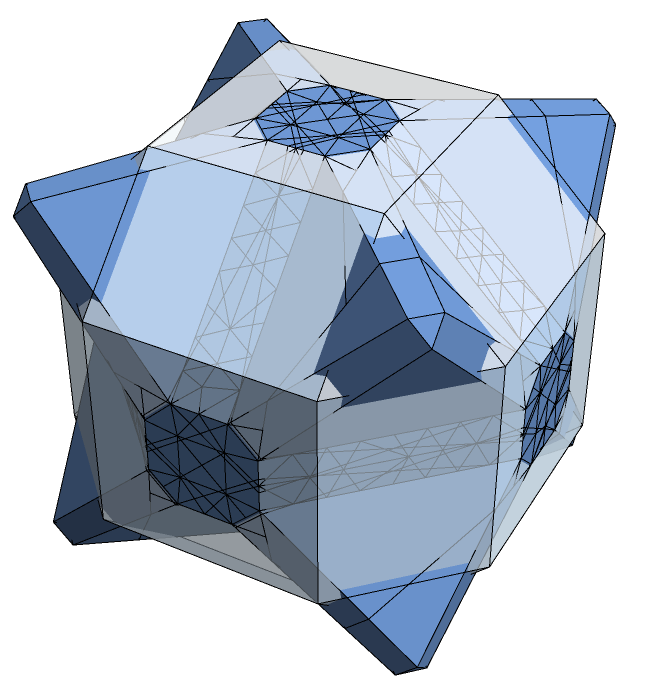}
    \caption{For the truncated cube, $P$ (left), the gravitational illumination body, $I^\delta(P)$ (middle) is not convex. Plotted together (right), we see neither $I^\delta(P)\subset P$ nor $P\subset I^\delta(P)$.}
    \label{FigureGIBtruncubeexample}
\end{figure}

\vskip 3mm
\noindent

\begin{proof}
    Note $P$ has 24 vertices. The faces of $P$ are one of two shapes: 1) octagons with normals that are permutations of $N=(\pm 1,0,0)$, area given by a $4\times 4$ square minus 4 triangles with base=height=1, giving area=14, and with center of gravity equal to the normal vector. The defining support planes then have the form $H=\{x|\langle x,N\rangle=2\}$. 2) triangles with normals that are permutations of $N=\frac{1}{\sqrt{3}}(\pm 1, \pm 1, \pm 1)$, area=$\sqrt{3}/2$, and center of gravity $(\frac{5}{3},\frac{5}{3},\frac{5}{3})=\frac{5}{\sqrt{3}}N$. The defining support planes then have the form $H=\{x|\langle x,N\rangle=5\}$.

    Then for $\delta$ sufficiently small, since each of the 24 vertices of $P$ are at the intersection of 3 planes, $P^\delta$ has 72 vertices.
    
    For $\delta<14/3$, 48 of the vertices are the permutations of $v=(\pm(1-\frac{3\delta}{14}),\pm(2+\frac{3\delta}{14}),\pm 2)$. One can verify that each such vertex is indeed at the intersection of two defining support planes of $P$ and that $[v,P]\backslash P$ is a pyramid with octagonal base and volume $\delta$. For example, take the vertex $v=(1-\frac{3\delta}{14},2+\frac{3\delta}{14},2)$. Then $\langle v,(0,0,1)\rangle=2$ and $\langle v,(1,1,1)\rangle=5$ and the volume of the pryamid is
    $$V=\frac{1}{3}(14)h=\frac{1}{3}(14)\langle v-(0,2,0),(0,1,0)\rangle=\frac{1}{3}(14)\left(\frac{3\delta}{14}\right)=\delta.$$
    
    Additionally requiring $\delta<1/6$, the remaining 24 vertices are the permutations of $v=(\pm(1+6\delta), \pm 2, \pm 2)$. Again one verifies $v$ is at the intersection of two of the defining planes of $P$ and that $[v,P]\backslash P$ is a pyramid with triangular base and volume $\delta$. For example, take the vertex $v=(1+6\delta, 2, 2)$. Then $\langle v,(0,0,1)\rangle=2$ and $\langle v,(0,1,0)\rangle=2$ and the volume of the pryamid is
    $$V=\frac{1}{3}\left(\frac{\sqrt{3}}{2}\right)h=\frac{1}{3}\left(\frac{\sqrt{3}}{2}\right)\left\langle v-\left(\frac{5}{3},\frac{5}{3},\frac{5}{3}\right),\frac{1}{\sqrt{3}}(1,1,1)\right\rangle=\frac{1}{3}\left(\frac{\sqrt{3}}{2}\right)\left(\frac{6\delta}{\sqrt{3}}\right)=\delta.$$

    Then the vertices of $I^\delta(P)$ for $0<\delta<1/6$ are all permutations of $(\pm(\frac{1}{4}-\frac{3\delta}{56}),\pm(\frac{1}{2}),\pm(2+\frac{3\delta}{56}))$ and $(\pm(\frac{3}{2}+\frac{3\delta}{2}),\pm\frac{7}{4},\pm \frac{7}{4})$. One can verify using the fact that in the case $x$ sees only one facet $F_i$, the map $T$ simplifies to $T(x)=\frac{x}{n+1}+\frac{n}{n+1}g(F_i)$, and choosing the appropriate center of gravity $g(F_i)$.

    Finally, we have that the 4 vertices $w_1=(\frac{1}{4}-\frac{3\delta}{56},\frac{1}{2},2+\frac{3\delta}{56}), w_2=(\frac{1}{4}-\frac{3\delta}{56},2+\frac{3\delta}{56},\frac{1}{2}), w_3=(-\frac{1}{4}+\frac{3\delta}{56},\frac{1}{2},2+\frac{3\delta}{56}),$ and $w_4=(-\frac{1}{4}+\frac{3\delta}{56},2+\frac{3\delta}{56},\frac{1}{2})$ define a facet of $I^\delta(P)$ with normal $N=\frac{1}{\sqrt{2}}(0,1,1)$. And, the vertices $v_1=(\frac{3}{2}+\frac{3\delta}{2},\frac{7}{4},\frac{7}{4})$ and $v_2=(-\frac{3}{2}-\frac{3\delta}{2},-\frac{7}{4},-\frac{7}{4})$ lie on opposite sides of the plane defined by the facet:
    $$\langle v_2,N\rangle<\langle w_i,N\rangle<\langle v_1,N\rangle$$
    since
    $$\langle w_i,N\rangle=\frac{1}{4\sqrt{2}}\left(10+\frac{3\delta}{14}\right)$$
    $$\langle v_1,N\rangle=\frac{7}{2\sqrt{2}}$$
    $$\langle v_2,N\rangle=-\frac{7}{2\sqrt{2}}$$
    Thus, $I^\delta(P)$ is not convex for $\delta<1/6$. We continue the computation to find the value of $\delta_0$ for which $I^\delta(K)$ is convex with $\delta>\delta_0$.

    For $\delta<1/6$, the 24 vertices of $P^\delta$ that lie in the intersection of two supporting planes of $P$ with normals that are permutations of $(\pm1,0,0)$ see only a triangular face of $P$. For $\delta>1/6$, these 24 vertices additionally see an octagonal face. Thus a new description is needed. We get for $1/6<\delta<44/3$, these 24 vertices are all permutations of $(\pm 2, \pm 2, \pm \frac{57+6\delta}{29})$. Again one verifies $v$ is at the intersection of two of the defining planes of $P$ and that $[v,P]\backslash P$ is union of pyramids with volume $\delta$. For example, take the vertex $v=( 2, 2, \frac{57+6\delta}{29})$. Then $\langle v,(1,0,0)\rangle=2$ and $\langle v,(0,1,0)\rangle=2$ and the volume of the pryamid is
    \begin{align*}
        V&=\frac{1}{3}\left(\frac{\sqrt{3}}{2}\right)h_1+\frac{1}{3}(14)h_2\\
        &=\frac{1}{3}\left(\frac{\sqrt{3}}{2}\right)\left\langle v-\left(\frac{5}{3},\frac{5}{3},\frac{5}{3}\right),\frac{1}{\sqrt{3}}(1,1,1)\right\rangle+\frac{1}{3}(14)\langle v-(0,0,2),(0,0,1)\rangle\\
        &=\delta.
    \end{align*}

    Then, for $v=( 2, 2, \frac{57+6\delta}{29})$, the corresponding vertex of $I^\delta(P)$ is
    \begin{align*}
        g(v)&=\frac{v}{n+1}+\frac{1}{(n+1)\delta}\left(\vol_{n-1}(F_1)\langle v-g(F_1),N_1\rangle g(F_1)+\vol_{n-1}(F_2)\langle v-g(F_2),N_2\rangle g(F_2) \right)\\
        &=\frac{v}{4}+\frac{1}{4\delta}\left(\frac{\sqrt{3}}{2}\left\langle v-\left(\frac{5}{3},\frac{5}{3},\frac{5}{3}\right),\frac{1}{\sqrt{3}}(1,1,1)\right\rangle\left(\frac{5}{3},\frac{5}{3},\frac{5}{3}\right)+14\langle v-(0,0,2),(0,0,1)\rangle(0,0,2)\right)\\
        &=\left(\frac{189\delta+70}{348\delta},\frac{189\delta+70}{348\delta}, \frac{-7+345\delta+9\delta^2}{174\delta}\right)
    \end{align*}
    
    and by symmetry we get that the 24 new vertices are all permutations of points of the form $\left(\pm\frac{189\delta+70}{348\delta},\pm\frac{189\delta+70}{348\delta}, \pm\frac{-7+345\delta+9\delta^2}{174\delta}\right)$. Then letting $v_1=\left(\frac{189\delta+70}{348\delta},\frac{189\delta+70}{348\delta},\frac{-7+345\delta+9\delta^2}{174\delta}\right)$ and $v_2=\left(-\frac{189\delta+70}{348\delta},-\frac{189\delta+70}{348\delta},-\frac{-7+345\delta+9\delta^2}{174\delta}\right)$ and considering again the plane defined by $w_1,w_2, w_3,$ and $w_4$, we have:
    $$\langle v_2,N\rangle<\langle w_i,N\rangle<\langle v_1,N\rangle$$
    since
    $$\langle w_i,N\rangle=\frac{1}{4\sqrt{2}}\left(10+\frac{3\delta}{14}\right)$$
    $$\langle v_1,N\rangle=\frac{56+879\delta+18\delta^2}{348\sqrt{2}\delta}$$
    $$\langle v_2,N\rangle=-\frac{56+879\delta+18\delta^2}{348\sqrt{2}\delta}$$
    Thus, we get that $I^\delta(P)$ is not convex for $1/6<\delta<14/3$.

    For $\delta<14/3$, the 48 vertices of $P^\delta$ that lie in the intersection of one supporting plane of $P$ with normal that is a permutation of $(\pm1,0,0)$ and a second supporting plane of $P$ with normal that is a permutation of $\frac{1}{\sqrt{3}}(\pm 1,\pm 1,\pm 1)$ see only an octagonal face of $P$. For $\delta>14/3$, these 48 vertices additionally see a triangular face. Thus a new description is needed. We get for $14/3<\delta<44/3$, these 48 vertices are permutations of $\left(\pm\left(1-\frac{3\delta+1}{15}\right),\pm2,\pm\left(2+\frac{3\delta+1}{15}\right)\right)$. One can verify that vertices of this form add on a set of volume $\delta$. Taking for example $v=\left(1-\frac{3\delta+1}{15},2,2+\frac{3\delta+1}{15}\right)$, we get
        \begin{align*}
        V&=\frac{1}{3}\left(\frac{\sqrt{3}}{2}\right)h_1+\frac{1}{3}(14)h_2\\
        &=\frac{1}{3}\left(\frac{\sqrt{3}}{2}\right)\left\langle v-\left(-\frac{5}{3},\frac{5}{3},\frac{5}{3}\right),\frac{1}{\sqrt{3}}(-1,1,1)\right\rangle+\frac{1}{3}(14)\langle v-(0,0,2),(0,0,1)\rangle\\
        &=\delta.
    \end{align*}

    Then, for $v=\left(1-\frac{3\delta+1}{15},2,2+\frac{3\delta+1}{15}\right)$, the corresponding vertex of $I^\delta(P)$ is
    \begin{align*}
        g(v)&=\frac{v}{n+1}+\frac{1}{(n+1)\delta}\left(\vol_{n-1}(F_1)\langle v-g(F_1),N_1\rangle g(F_1)+\vol_{n-1}(F_2)\langle v-g(F_2),N_2\rangle g(F_2) \right)\\
        &=\frac{v}{4}+\frac{1}{4\delta}\left(\frac{\sqrt{3}}{2}\left\langle v-\left(-\frac{5}{3},\frac{5}{3},\frac{5}{3}\right),\frac{1}{\sqrt{3}}(-1,1,1)\right\rangle\left(-\frac{5}{3},\frac{5}{3},\frac{5}{3}\right)+14\langle v-(0,0,2),(0,0,1)\rangle(0,0,2)\right)\\
        &=\left(\frac{70+27\delta-9\delta^2}{180\delta},\frac{21\delta-14}{36\delta}, \frac{14+360\delta+9\delta^2}{180\delta}\right)
    \end{align*}
    and by symmetry, the corresponding 48 vertices of $I^\delta(P)$ are permutations of the points $\left(\pm\frac{70+27\delta-9\delta^2}{180\delta},\pm\frac{21\delta-14}{36\delta}, \pm\frac{14+360\delta+9\delta^2}{180\delta}\right)$. Letting $w_i=\left(\pm\frac{70+27\delta-9\delta^2}{180\delta},\frac{21\delta-14}{36\delta}, \frac{14+360\delta+9\delta^2}{180\delta}\right)$ and permuting the $y$ and $z$ components yields four points that define a facet of $I^\delta(P)$ with normal $N=\frac{1}{\sqrt{2}}(0,1,1)$. For $v_1=\left(\frac{189\delta+70}{348\delta},\frac{189\delta+70}{348\delta},\frac{-7+345\delta+9\delta^2}{174\delta}\right)$ and $v_2=\left(-\frac{189\delta+70}{348\delta},-\frac{189\delta+70}{348\delta},-\frac{-7+345\delta+9\delta^2}{174\delta}\right)$, we have

        $$\langle v_2,N\rangle<\langle w_i,N\rangle<\langle v_1,N\rangle$$
    since
    $$\langle w_i,N\rangle=\frac{9\delta^2+465\delta-56}{180\sqrt{2}\delta}$$
    $$\langle v_1,N\rangle=\frac{56+879\delta+18\delta^2}{348\sqrt{2}\delta}$$
    $$\langle v_2,N\rangle=-\frac{56+879\delta+18\delta^2}{348\sqrt{2}\delta}$$
    Thus, we get that $I^\delta(P)$ is not convex for $14/3<\delta<44/3$.
\end{proof}

\vskip 3mm
 \noindent
    \textbf{Data availability statement}.
    Data sharing not applicable to this article, as no datasets were generated or analyzed during the current study. 
\vskip 2mm
 \noindent  
    {\bf Conflict of interest}. The corresponding author declares, on behalf of all authors, that there are no conflicts of interest.
\vskip 2mm
 \noindent    
{\bf Statement on the use of AI}.
The  mathematical ideas, the research direction and  the
first version of the paper  predate this use of AI, see \cite{Glasgo2020}. 
During the recent preparation of the manuscript, the authors used ChatGPT (OpenAI) to assist with
the verification of some mathematical arguments of some lemmas.
AI tools also assisted with minor editing.
The authors take full responsibility for the mathematical content and the final
manuscript.

\vskip 10mm

\noindent
Victor Glasgo\\
{\small Department of Mathematics}\\
{\small Kennesaw State University}\\
{\small Kennesaw, Georgia 30144, U. S. A.}\\
{\small \tt vglasgo@kennesaw.edu}\\ \\
\vskip 3mm
\noindent
Elisabeth M. Werner\\
{\small Department of Mathematics \\
{\small Case Western Reserve University \\
{\small Cleveland, Ohio 44106, U. S. A. \\
{\small \tt elisabeth.werner@case.edu}

\end{document}